\documentclass{article}
\usepackage{manuscriptStyleFile} 

\usepackage{authblk}
\title{Image reconstructions using sparse dictionary representations and implicit, non-negative mappings}
\author[1]{Elizabeth Newman}
\author[1]{Jack Michael Solomon}
\author[1]{Matthias Chung}

\affil[1]{\centering Department of Mathematics, Emory University, Atlanta, GA, USA $\big\{$\email{elizabeth.newman},~\email{jack.michael.solomon},~\email{matthias.chung}$\big\}$\email{@emory.edu}}
\date{\today}

\begin{document}


\maketitle

\begin{abstract}
    Many imaging science tasks can be modeled as a discrete linear inverse problem.  Solving linear inverse problems is often challenging, with ill-conditioned operators and potentially non-unique solutions.  Embedding prior knowledge, such as smoothness, into the solution can overcome these challenges. In this work, we encode prior knowledge using a non-negative patch dictionary, which effectively learns a basis from a training set of natural images.   In this dictionary basis, we desire solutions that are non-negative and sparse (i.e., contain many zero entries).  With these constraints, standard methods for solving discrete linear inverse problems are not directly applicable.  One such approach is the modified residual norm steepest descent (MRNSD), which produces non-negative solutions but does not induce sparsity.  In this paper, we provide two methods based on MRNSD that promote sparsity.  In our first method, we add an $\ell_1$-regularization term with a new, optimal step size. In our second method, we propose a new non-negative, sparsity-promoting mapping of the solution.  We compare the performance of our proposed methods on a number of numerical experiments, including deblurring, image completion, computer tomography, and superresolution. Our results show that these methods effectively solve discrete linear inverse problems with non-negativity and sparsity constraints. 
\end{abstract}

\section{Introduction}
\label{sec:introduction}

\emph{Inverse problems} are an important class of mathematical problems with the aim to re- or uncover unobserved information from noisy data paired with an underlying mathematical or statistical model. Inverse problems are ubiquitous in modern science, from data analysis and machine learning to many imaging applications, such as image deblurring, computer tomography, and super-resolution \cite{hansen2010discrete}.

In this work, we consider the common imaging problem that is given as a discrete linear inverse problem
\begin{equation}\label{eq:inverseproblem}
    \bfb = \bfC\bfy_{\rm true} + \bfeta,
\end{equation}
where $\bfy_{\rm true} \in \calY \subset \bbR^n_+$ is the desired solution (a vectorized representation of the image is some image space), $\bfC \in \bbR^{m\times n}$ is a given forward imaging process and $\bfb\in\bbR^m$ contains observed measurements. We assume these observations $\bfb$ are corrupted by additive Gaussian noise $\bfeta \sim \calN({\bf0},\sigma^2\bfI_m)$. 
In most cases, the underlying operator $\bfC$ is ill-conditioned (i.e., the solution is sensitive to small perturbations in $\bfb$) or underdetermined ($m < n$) such that the solution is not unique. 
In most of these common settings, prior knowledge is required for selecting a meaningful solution. 
One approach to embedding prior knowledge into the solution is to cast~\eqref{eq:inverseproblem} as a \emph{variational inverse problem} of the form~\cite{ChungGazzola2023:hybrid, Hansen1998:rankDeficient}
\begin{equation}\label{eq:variationalip}
    \min_{\bfy\in\calY} \quad  \thf \norm[2]{\bfC\bfy-\bfb}^2 + \mu \calR(\bfy).    
\end{equation}
The first term measures the data fidelity between the model prediction, $\bfC\bfy$, and observations, $\bfb$, with discrepancy measured in $\ell^2$-norm. For $\mu = 0$, this corresponds to a maximum likelihood estimator with Gaussian noise distribution.  The second term, $\calR(\mdot)$, denotes an appropriate regularizer that promotes desirably properties in the solution, $\bfy$, (e.g., small norm, smoothness). 
The regularization parameter $\mu > 0$ balances the weight of the data fidelity and regularization terms, and is often chosen heuristically or based on noise estimates~\cite{Golub1979:gcv}.

Regularization for imaging problems is motivated by the fact that the \emph{natural image space} $\calY \subset \bbR^n_+$ is minuscule compared to the set of all images of dimension $\bbR^n_+$. For instance, drawing a random sample from $\bbR^n_+$ will most likely look like noise \cite{brunton2019data}. In other words, natural images exhibit a high degree of redundancy and can be efficiently captured within a low-dimensional manifold. Hence, the features and the structure of natural images must be included in an inverse problem process.One common approach is to utilize an edge-preserving \emph{total variation} (TV) regularization $\calR(\bfy) = \norm[p]{\bfL\bfy}$. Here, $\bfL\bfy$ is a finite difference approximation of the directional pixel gradients, and $\norm[p]{\bfy}$ denotes the $\ell^p$-norm~\cite{Osher2003:TV}.   
This type of regularization has shown tremendous success in various applications. However, disadvantages of this approach include the potential over-smoothing of the reconstructed image (apart from the edges) and potential inconsistent intensity shifts in $\bfy$, \cite{chan2005image,strong2003edge,elad2010sparse}. 

An entirely different approach to regularize is to represent the solution to~\eqref{eq:variationalip} in terms of a pre-trained image data ``basis'' or \emph{dictionary}~\cite{SoltaniAH15}.  Dictionaries of natural images or image elements provide an informative prior and incorporate expected image features (e.g., edges, smoothness) into $\bfy$ by directly.  To be precise, we represent $\bfy$ as a linear combination of vectors (atoms) $\bfg_k \in \bbR^n_+$ for $k = 1,\ldots, K$, i.e., $\bfy = \bfG\bfx$, where $\bfG = [\bfg_1,\ldots, \bfg_K]$ is referred to as the \emph{dictionary} and $\bfx\in \bbR^K_+$ are appropriate non-negative selection coefficients. If $K>n$ or even $K\gg n$, we say the dictionary is \emph{overcomplete}.  As a result, the coefficient vector is not unique and larger than the original image.  In this setting, we desire coefficients that exhibit sparsity, i.e., many elements in $\bfx$ are being zero, to ensure we do not increase the storage of our solution. 

To select such a suitable and sparse coefficient vector $\bfx$ (and therefore a suitable $\bfy$), we ideally minimize
\begin{equation}\label{eq:sparsereconstruction}
    \min_{\bfx\geq {\bf0}} \quad \norm[0]{\bfx} \quad \mbox{subject to} \quad \thf\norm[2]{\bfC\bfG\bfx -\bfb}^2 \leq \delta^2, 
\end{equation}
where $\norm[0]{\bfx}$ is not a proper norm and is defined as the cardinality of non-zero elements in $\bfx$, while $\delta$ represents the given noise level. Solving \eqref{eq:sparsereconstruction} is NP-hard and efficient approximation approaches need to be utilized, \cite{elad2010sparse,tovsic2011dictionary}. One approach in tackling \eqref{eq:sparsereconstruction} is to incrementally generate non-zero elements in $\bfx$ by iteratively selecting and updating the element in $\bfx$ minimizing the remaining residual until a prescribed accuracy threshold is reached. Methods following this approach include \emph{matching pursuit} and \emph{orthogonal matching pursuit}, \cite{mallat1993matching,tropp2004greed}.  A main disadvantage of matching pursuit approaches is that we cannot ensure optimally. Further, each iteration requires the computation of inner products with the remaining dictionary adding to the computational complexity of the algorithm. 

A common alternative approach to provide a non-negative and sparse $\bfx$ is based on convex relaxation. We approximate \eqref{eq:sparsereconstruction} by 
\begin{equation}\label{eq:lasso}
    \min_{\bfx\geq {\bf0}} \quad \thf\norm[2]{\bfA\bfx -\bfb}^2 + \lambda \norm[1]{\bfx}, 
\end{equation}
where $\bfA = \bfC\bfG$ and $\lambda>0$ is an appropriate regularization parameter~\cite{candes2006stable}. When $\bfx$ is unrestricted, the optimization problem \eqref{eq:lasso} is the so-called least absolute shrinkage and selection operator (LASSO) problem. Methods such as \emph{basis pursuit}, \emph{Alternating Direction Method of Multipliers} (ADMM), and \emph{Variable Projected Augmented Lagrangian} (VPAL) have been shown to efficiently solve an unconstrained version  of \eqref{eq:lasso}; see \cite{chen2001atomic,boyd2011distributed, chung2022variable}.  Additionally, Krylov subspace methods are a major tool for large-scale least squares problems~\cite{saad2003iterative}.  However, none of these approaches are directly applicable in this setting because they do not inherently enforce non-negativity or promote sparse solutions. Consequently, alternative approaches have been developed to address these constraints. Here, we follow another approach, we propose to use a method designed for non-negative least squares and extend these approaches to also provide sparse solutions. 

Our work is structured as follows. In \Cref{sec:background} we provide details on the modified residual norm steepest descent methods which have been developed to solve least squares problems with non-negativity constraints and we further provide a background on patch dictionary learning. In \Cref{sec:mainconstributions} we discuss our proposed methods to solve non-negative least squares problems with sparsity priors.  We demonstrate the effectiveness of our approaches with various numerical experiments in \Cref{sec:numericalexperiments} including deblurring (\Cref{sec:deblur}), image completion (\Cref{sec:indicator}), and computer tomography  (\Cref{sec:tomography}). We close up our investigations with some conclusions and further outlook in \Cref{sec:conclusion}.

\section{Background}\label{sec:background}

Because images are inherently non-negative, we are interested in the constrained optimization problem 
\begin{equation}\label{eq:nnleastsqures}
    \min_{\bfx\geq {\bf0}} \quad f(\bfx) \equiv \thf\norm[2]{\bfA\bfx -\bfb}^2.
\end{equation}
While standard constrained optimization techniques may be employed for \eqref{eq:nnleastsqures}, they often suffer from poor computational efficiency \cite{nocedal1999numerical}. To address this issue, specialized methods like the expectation-maximization (EM) algorithm or modified residual norm steepest descent (\emph{MRNSD}) have been proposed. These methods have shown great success by directly handling the linear least squares fitting term with non-negativity constraints \cite{Kaufman1993:MRNSD, NagyZdenek2000:MRNSD}.

\subsection{Modified Residual Norm Steepest Descent (MRNSD)}
The main idea behind MRNSD is to enforce non-negativity by introducing an implicit, element-wise exponential variable transformation and performing a gradient descent update directly on $\bfx$ with appropriate step size selection.  Let us introduce the element-wise variable transformation $x_i = \e^{z_i}$, for $i = 1,\ldots, K$, short we write $\bfx = \e^\bfz$. This mapping ensure $\bfx\geq {\bf0}$ and the optimization problem \eqref{eq:nnleastsqures} reduces to 
\begin{equation}
    \min_{\bfz} \quad \hat{f}(\bfz) \equiv \thf \norm[2]{\bfA \e^\bfz - \bfb}^2.
\end{equation}
Note the gradient of $\hat{f}$, the objective function over $\bfz$-space, is given as
\begin{equation}
    \nabla_\bfz \hat{f}(\bfz) = \diag{\e^\bfz} \bfA\t (\bfA \e^\bfz - \bfb) =\diag{\bfx} \nabla_\bfx f(\bfx).
\end{equation}
MRNSD performs a steepest descent direction update with respect $\bfx$ using the gradient of $\widehat{f}$, i.e.,
\begin{align}
    \bfs = - \nabla_\bfz \hat{f}(\bfz)=  - \diag{\bfx} \bfg,
\end{align}
where $\bfg = \nabla_\bfx f(\bfx) = \bfA\t(\bfA\bfx -\bfb)$ is the gradient of the $f$, the objective function in $\bfx$-space.
We note that $\bfs$ is a descent direction in $\bfx$-space, i.e., $\bfs^\top \bfg < 0$, because of the assumed positivity of $\bfx$.  This ensures $\diag{\bfx}$ is a positive definite matrix, which in turn ensures we move in a descent direction. Further note that an optimal step size can be computed, i.e.,
\begin{align}
\alpha = \frac{-\bfs\t\bfg}{\bfs\t\bfA\t\bfA\bfs}.
\end{align}
However, since MRNSD updates $\bfx$ directly (and not the surrogate variable $\bfz$), there may exist a selected step length $\alpha>0$ for which the update $\bfx + \alpha \bfs$ can have negative elements. Hence, thresholding must be performed to maintain non-negativity in $\bfx$ which is provided by
\begin{align}
\alpha = \min \left\{-\bfs^\top \bfg / (\bfs^\top \bfA^\top \bfA \bfs), \min_{s_i < 0} \{-x_i/s_i)\}_{i = 1}^K\right\}.
\end{align}

The MRNSD method is summarized in \Cref{alg:mrnsd}.
\begin{algorithm}
\caption{Modified Residual Norm Steepest Descent (MRNSD)~\cite{NagyZdenek2000:MRNSD, Kaufman1993:MRNSD}}\label{alg:mrnsd}
\begin{algorithmic}[1]
\State Solve $\min_{\bfx \ge 0} \frac{1}{2}\|\bfA\bfx - \bfb\|_2^2$
\While{not converged}
    \State $\bfg = \bfA^\top (\bfA \bfx - \bfb)$
    \State $\bfs = -\diag{\bfx}\bfg$ \Comment{modified steepest descent direction}
    \State $\alpha = \min \left\{-\bfs^\top \bfg / (\bfs^\top \bfA^\top \bfA \bfs), \min_{s_i < 0} \{-x_i/s_i)\}_{i = 1}^K\right\}$ 	 
    					\Comment{optimal step size with bounded line search}\label{alg:mrnsd_step_size}    
    \State $\bfx \gets \bfx + \alpha \bfs$
\EndWhile
\end{algorithmic}
\end{algorithm}

Note that MRNSD only performs a surrogate variable transformation implicitly but never constructs nor updates $\bfz$. While MRNSD has been shown to effectively address non-negativity constraints, the utilized exponential mapping prevents values of $\bfx$ from becoming zero and hence does not naturally promote sparsity in the solution. Achieving sparsity requires ad hoc threshholding in postprocessing, and motivates our exploration of efficient methods that enforce both non-negativity and sparsity during optimization. 

\subsection{Dictionary Learning}
\begin{figure}
    \centering
    \begin{tikzpicture}
        \node (n0) at (0,0) {\includegraphics[width=0.15\linewidth]{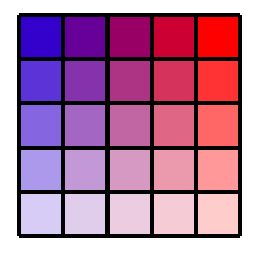}};  
        \node[below=0.0cm of n0.south, anchor=north] (y){\begin{tabular}{c}  
        $\bfY$ \\ $M\times N$\end{tabular}};

        \node[right=1.0cm of n0.east, anchor=west] (n1) {\includegraphics[width=0.3\linewidth]{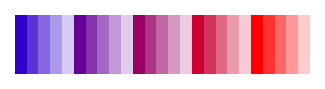}};

        \node at (y -| n1) {\begin{tabular}{c}  
        $\texttt{patchify}(\bfY)$ \\ $pq\times r$\end{tabular}};
        
    \end{tikzpicture}

\caption{Illustration of image patchification.  For simplicity, we illustrate non-overlapping patches that exactly partition the image with $r = (M/p)(N/q)$ total patches.  In {\sc Matlab}, this operation is given by \texttt{im2col(Y,[p,q],'distinct')} and is reversible.}
\label{fig:patchify}
\end{figure}

Before we turn our attention to sparsity-promoting alternatives to MRNSD, we next provide some details on patch dictionary learning and the resulting sparse representations~\cite{SoltaniAH15, SoltaniKilmerHansen2016}.  Suppose we have a set of $t$ training $M\times N$ images $\{\bfY_1,\cdots, \bfY_t\}$ with $\vec{\bfY_j} \in \Ycal \subset \Rbb_+^{MN}$, where $\myVec(\cdot)$ vectorizes the image column-wise. We collect (overlapping) patches of size $p\times q$ with $p< M$ and $q < N$,  ``patchify'' the images in our training set, and vectorize those patches and store them as columns of a matrix $\bfY^{\rm patch} \in \Rbb_+^{pq \times tr}$ where $r$ is the number of patches per image; see~\Cref{fig:patchify} for an illustration. 

A goal in dictionary learning is to construct an appropriate non-negative patch dictionary $\bfD\in \Rbb_+^{pq\times s}$ with a corresponding set of non-negative, column-wise sparse coefficients $\bfX^{\rm patch} \in \Rbb_+^{s\times r}$ such that $\bfY^{\rm patch} \approx \bfD \bfX^{\rm patch}$.  The desired properties of the decomposition ensure that each patch of our training set will be constructed using only a few dictionary atoms, producing a ``parts-of-the-whole'' type of basis~\cite{LeeSeung1999:nmf}.  

\begin{figure}
    \centering
    \subfloat[Flower dictionary patches from a subset of images resized to $64\times 64$ from~\cite{tfflowers}.  The dictionary has $400$ atoms, meaning $\bfD_{\rm flower}\in \Rbb_+^{256\times 400}$. \label{fig:flower_dictionary}]{\begin{tikzpicture}
        \def\w{1}
        \def\s{0.05}
        \foreach \i in {1,...,5}{
            \node at (\i*\w+\i*\s,0) {\includegraphics[width=\w cm]{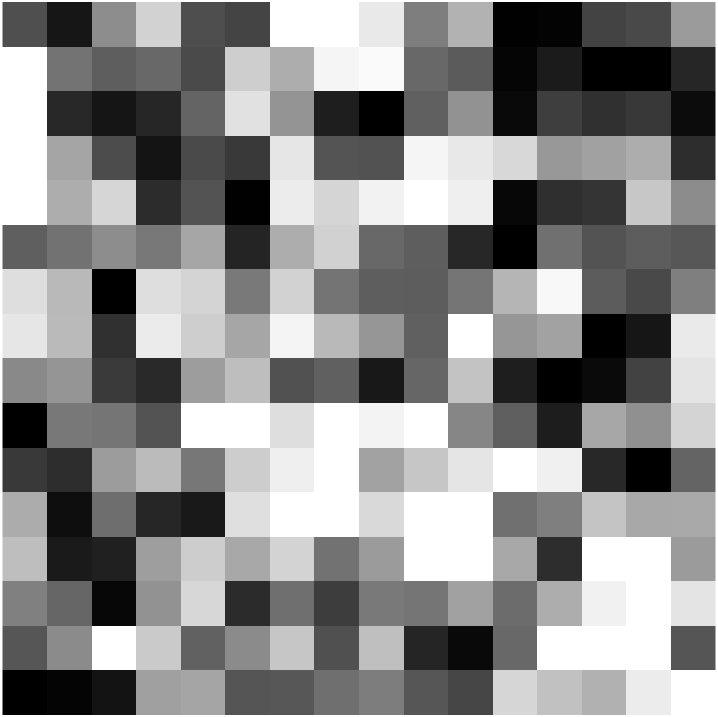}};

            \pgfmathsetmacro{\j}{int(\i + 5)}
            \node at (\i*\w+\i*\s,-\w-\s) {\includegraphics[width=\w cm]{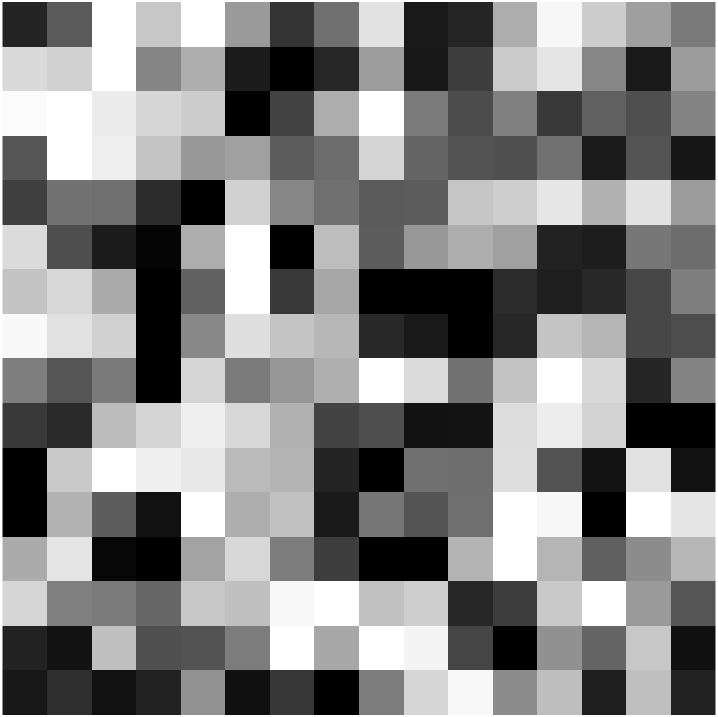}};
        }
        \end{tikzpicture}
        }
\hspace{0.5cm}
    \subfloat[Earth dictionary patches from a single high-resolution image of the earth of size $675 \times  1,\!200$~\cite{earthImg}. 
    The dictionary has $400$ atoms, meaning $\bfD_{\rm earth}\in \Rbb_+^{256\times 400}$. \label{fig:earth_dictionary}]{\begin{tikzpicture}
         \def\w{1}
        \def\s{0.05}
        \def\y{0}
        \foreach \i in {1,...,5}{
            \node at (\i*\w+\i*\s,-\y) {\includegraphics[width=\w cm]{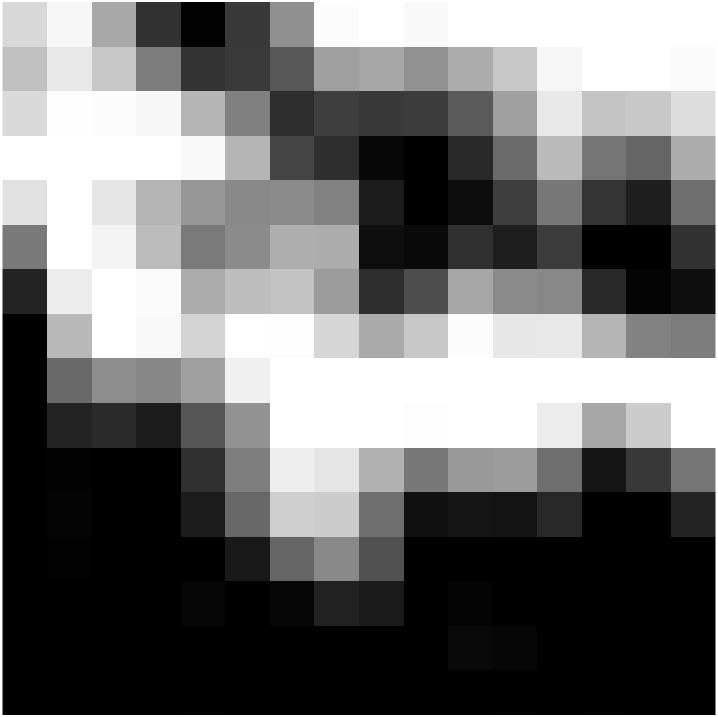}};

            \pgfmathsetmacro{\j}{int(\i + 10)}
            \node at (\i*\w+\i*\s,-\w-\s-\y) {\includegraphics[width=\w cm]{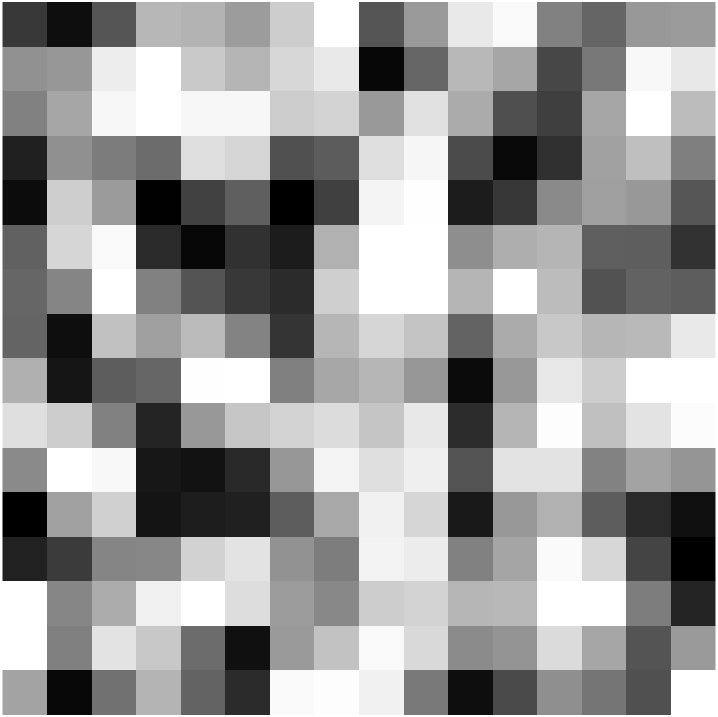}};
        }
    \end{tikzpicture}
    }
    \caption{Sample of $10$ patches from dictionaries generated from two different sources.  For training, we collect overlapping patches of size $16\times 16$.  The patches share similar binary-looking patterns because of the box constraints~\eqref{eq:boxConstraints}.}
    \label{fig:dictionary_patches}
\end{figure}

Following the presentation in~\cite{SoltaniAH15}, we form this non-negative matrix factorization by solving the optimization problem
    \begin{align}\label{eq:DLOpt}
        \min_{\bfD \in \Dcal, \bfX^{\rm patch} \in \Rbb_+^{s\times r}} \tfrac{1}{2}\norm[\fro]{\bfY^{\rm patch} - \bfD \bfX^{\rm patch}}^2 + \beta \sum_{j=1}^r \|\bfx_j^{\rm patch}\|_1,
    \end{align}
where $\bfx_j^{\rm patch}\in \Rbb^{s}$ is the $j$-th column of $\bfX^{\rm patch}$, $\beta>0$ is the regularization parameter for the sparsity-promoting $\ell_1$-regularization, and $\Dcal$ restricts the values of the dictionary to lie between $0$ and $1$. In our experiments, we use box constraints
    \begin{align}\label{eq:boxConstraints}
        \Dcal := \left\{\bfD \in \Rbb^{pq\times s} \mid 0 \le d_{ij}  \le 1 \text{ for }i=1,\dots,pq \text{ and }j=1,\dots,s\right\}.
    \end{align}
The dictionary learning problem~\eqref{eq:DLOpt} is non-convex with respect to $(\bfD, \bfX^{\rm patch})$ jointly, but is convex with respect to each variable separately. This motivates the use of an alternating optimization strategy. To this end, we consider the alternating direction method of multipliers (ADMM) approach proposed by~\cite{SoltaniAH15}, which reformulates~\eqref{eq:DLOpt} as
    \begin{subequations}\label{eq:DLOptADMM}
    \begin{alignat}{3}
    &\min_{\bfD, \bfX^{\rm patch}, \bfU, \bfV}&&\quad\tfrac{1}{2} \norm[\fro]{\bfY^{\rm patch} - \bfU \bfV}^2 + \beta \sum_{i=1}^t \norm[1]{\bfx_j^{\rm patch}} + \chi_{\Dcal}(\bfD) + \chi_{\Rbb_+^{s\times r}}(\bfX^{\rm patch}), \\
    &\text{subject to} &&\quad  \bfD = \bfU \quad \text{ and } \quad \bfX^{\rm patch} = \bfV,
    \end{alignat}
    \end{subequations}
where $\chi_{\Zcal}$ is the characteristic function over the set $\Zcal$ such that $\chi_{\Zcal}(z) = 0$ if $z \in \Zcal$ and $\chi_{\Zcal}(z) = +\infty$ if $z\not\in \Zcal$.   Introducing the auxiliary $\bfU$ and $\bfV$ ensures that the function is separable over the optimization variables, a necessary condition for ADMM. We solve~\eqref{eq:DLOptADMM} using the augmented Lagrangian approach~\cite{boyd2011distributed}; for further details, see~\cite{SoltaniAH15, SoltaniKilmerHansen2016, NewmanKilmer2020}.  Examples of generated dictionary patches are shown in~\Cref{fig:dictionary_patches}. 

The patch dictionary presentation is related to the original global dictionary formulation in~\Cref{sec:introduction}. 
Let $\bfX_i\in \Rbb_+^{pq \times r}$ be the coefficients of the patches drawn from training image $\bfY_i$ and assume the patches are non-overlapping. Then,
    \begin{align}
        \myVec(\bfY_i) = \bfP\myVec(\bfD \bfX_i) = \underbrace{\bfP(\bfI \otimes \bfD)}_{\bfG} \myVec(\bfX_i),
    \end{align}
where $\otimes$ denotes the Kronecker product, $\bfI$ is the $r\times r$ identity matrix, and $\bfP$ is a permutation matrix mapping patches to linear indices~\cite{matrixCookbook}.  We require the permutation matrix because $\myVec(\bfX_i)$ stacks vectorized patches column-wise, but $\myVec(\bfY_i)$ concatenates the columns of the image vertically, which order the pixels differently.

\section{Non-negativity and Sparsity-Promoting Mappings} 
\label{sec:mainconstributions}

We consider two alternative approaches based on MRNSD to promote sparse solutions.  In the first approach, we introduce an $\ell_1$-regularization term following the presentation~\cite{NewmanKilmer2020}. In the second approach, we modify the implicit non-negativity mapping to allow for entries to be set equal to zero.  In both cases, choose an optimal step size based on the new iterations. 

\paragraph{Sparsity-Promoting MRNSD}

We consider a sparsity-promoting MRNSD that minimizes the non-negativity-constrained lasso problem as introduced in~\eqref{eq:lasso}. The $\ell_1$-norm promotes sparsity when the solution switches from positive to negative entries. However, because unregularized MRNSD imposes the implicit mapping that enforces strict positivity of the solution, the $\ell_1$-regularization will not promote sparsity.  In response, we slightly modify the implicit mapping as $\bfx = \e^\bfz - \bfvarepsilon$ where $\bfvarepsilon > {\bf0}$ to allow for small negative values. Following proximal operator theory~\cite{ParikhBoyd2014:proximal}, the resulting sparsity-promoting MRNSD iteration is
    \begin{align}
    \bfx \gets S_{\alpha \lambda}(\bfx + \alpha \bfs),
    \end{align}
where the soft-thresholding operator $S_{\alpha \lambda}$ is applied entry-wise and given by the piecewise function\footnote{In {\sc Matlab}, we write \texttt{softThreshold(y,c) = sign(y) .* max(abs(y) - c, 0)}.} 
    \begin{align}\label{eq:softThreshold}
           (S_{\gamma}(\bfx))_i = \begin{cases}
            x_i + \gamma, &\mbox{for } x_i < -\gamma,\\
            0, & \mbox{for }  |x_i| \le \gamma, \\
            x_i - \gamma, & \mbox{for } x_i > \gamma,
        \end{cases} \qquad \mbox{for } i = 1,\ldots, K.
    \end{align}    
In this work, we choose the optimal step size for the soft-threshold-updated solution by solving
    \begin{align}\label{eq:optalpha}
        \alpha^* \in \min_{\alpha \in [0,u]}\tfrac{1}{2}\|\bfA S_{\alpha\lambda}(\bfx + \alpha \bfs) - \bfb\|_2^2, \qquad \text{where} \quad u = \min_{s_i < -\lambda} \{-x_i/(s_i + \lambda)\}_{i = 1}^K
    \end{align}
at each iteration. Solving the scalar optimization problem \eqref{eq:optalpha} can be performed efficiently utilizing nonlinear solvers for single variables such as golden section search with parabolic interpolation \cite{ascher2011first}. We summarize sparsity-promoting MRNSD (\texttt{spMRNSD}) in~\Cref{alg:sparsemrnsd}. We note that in practice we can pre-apply the operator $\bfA$ before solving~\eqref{eq:optalpha}, significantly reducing the computational cost. 

Note that the upper bound $u$ in~\eqref{eq:optalpha} is slightly different from the boundary for the original MRNSD step size (\Cref{alg:mrnsd}, \Cref{alg:mrnsd_step_size}). Specifically, when we step in $\bfx$-space in MRNSD, we adjust the step size $\alpha$ to ensure that the step $x_i + \alpha s_i$ does not violate the non-negativity constraint. Similarly, for \texttt{spMRNSD}, we must ensure that the smallest value in $S_{\alpha \lambda}$ remains non-negative; that is, whenever $x_i + \alpha s_i < -\alpha \lambda$ or rewritten $x_i + \alpha (s_i + \lambda) < 0$ is satisfied. Because $x_i$, $\lambda$ and $\alpha$ are assumed to be non-negative, this condition arises when $s_i$ is negative and sufficiently smaller than $-\lambda$.  Thus, $\alpha$ cannot be larger than $-x_i / (s_i + \lambda)$, which completes the upper bound of~\eqref{eq:optalpha}.

The original sparsity-promoting MRNSD from~\cite{NewmanKilmer2020} used the optimal MRNSD step size in~\Cref{alg:mrnsd}, \Cref{alg:mrnsd_step_size}. Instead of the optimal step size for the soft-thresholding operator. We note that the MRNSD step size is acceptable for \texttt{spMRNSD} in the sense that it will not introduce negative values in $\bfx$, but may produce inferior updates. However, an optimal step size based on the soft-thresholding operator could potentially be larger because it allows for negative search directions up to the level of $-\lambda$.  Additionally, a larger step size will make the middle region of~\eqref{eq:softThreshold} wider, potentially yielding sparser solutions.  This being said, the MRNSD step size may be larger or yield sparser solutions for certain problems, but these cases are difficult to predict and further analysis will be considered for future directions. 

\begin{algorithm}
\caption{Sparsity-Promoting MRNSD (\texttt{spMRNSD})~\cite{NewmanKilmer2020}}\label{alg:sparsemrnsd}
\begin{algorithmic}[1]
\State Solve $\min_{\bfx \ge {\bf0}}\ \frac{1}{2}\|\bfA\bfx - \bfb\|_2^2 + \lambda \|\bfx\|_1$
\While{not converged}
    \State $\bfg = \bfA^\top (\bfA \bfx - \bfb)$
    \State $\bfs = -\diag\bfx\bfg$
    \State compute  $\alpha^*$ solving \eqref{eq:optalpha}
    \State $\bfx \gets S_{\alpha^* \lambda}(\bfx + \alpha^* \bfs)$
\EndWhile
\end{algorithmic}
\end{algorithm}

\paragraph{Sparsity-Promoting Non-negative Mapping}

\begin{figure}
    \centering
    \begin{tikzpicture}
        \def\n{8}
        \draw[white] (-\n,0) -- (\n,0);
        
        \node at (0,0) {\includegraphics[width= 0.75\textwidth]{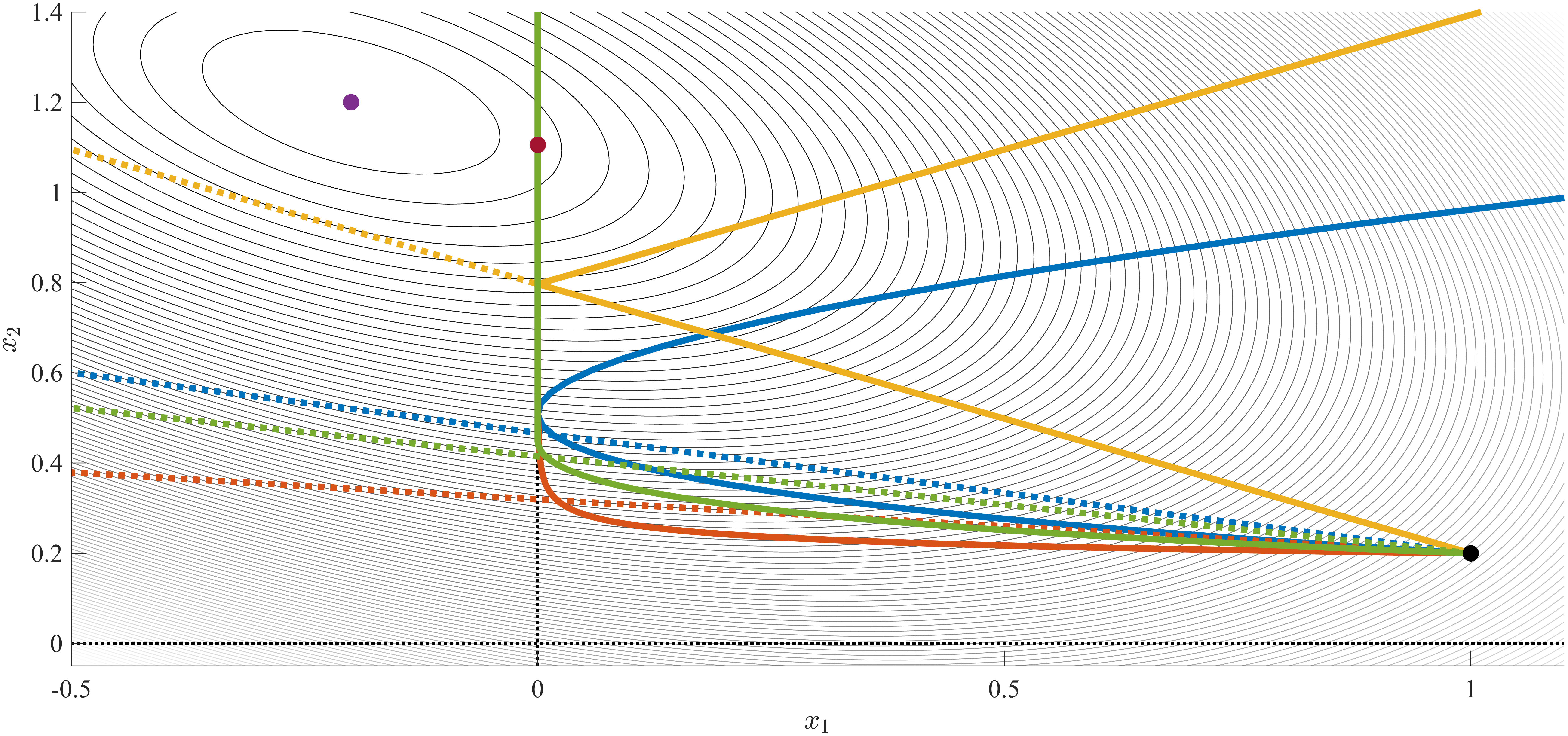}};

        \node at (-1.7,1.8) {$\widehat{\bfx}$};
        \node at (-3.05,2.15) {$\bfx^*$};
        \node at (5.75,-1.45) {$\bfx$};
        
        \node at (-1.95,3.25) {\color{mycolor5} $w_{1,0}(\bfz)$};
        \node at (5.5,3.25) {\color{mycolor3} $|\bfz|$};
        \node at (6.5,1.5) {\color{mycolor1} $\bfz^2$};
        \node at (0,-1.6) {\color{mycolor2} $e^{\bfz}$};
    \end{tikzpicture}
    
    \caption{Consider the toy problem $\min_{\bfx\geq {\bf0}} \ \norm[2]{\bfA\bfx-\bfb}^2$, with $\bfA = \begin{bmatrix}20 & 5\\ 5 & 20 \end{bmatrix}$ and $\bfb = \begin{bmatrix}2 \\  23 \end{bmatrix}$. The unconstrained minimizer is at $\bfx^* = 1/10\begin{bmatrix}-2 & 12 \end{bmatrix}\t$ (purple), while the minimizer of the constrained optimization problem is obtained at $\widehat \bfx = \begin{bmatrix}0 &  1.1059 \end{bmatrix}\t$ (dark red). We illustrate the different steepest descent paths in the $\bfx$-space. Let the current iteration start from $\bfx = \begin{bmatrix}0.15 & 0.13 \end{bmatrix}\t$. We depict various choices for MRNSD (line search in $\bfx$, dashed lines) and spNNGD (line search in $\bfz$, solid lines), with four different functions for $w$, i.e., $w(\bfz) = \bfz^2$ (blue), $w(\bfz) = \e^\bfz $ (red), $w(\bfz) = \abs{\bfz}$ (yellow), $w(\bfz) = w_{1,0}(\bfz)$ from \eqref{eq:nnmapping} (green). Note that with optimal step size selection, spNNGD with the $w_{1,0}(\bfz)$ would find the minimizer in just one iteration.}
    \label{fig:mappings}
\end{figure}

As a competing approach to \texttt{spMRNSD}, we focus on adapting the non-negativity MRNSD mapping. 
The core of MRNSD is to represent the solution via the \emph{implicit} non-negative mapping $\bfx = \e^\bfz$.  Here, we consider a broader class of point-wise, non-negativity mappings $w: \Rbb^K\to \Rbb_+^K$ and solve
\begin{align}
    \min_{\bfz} \quad \tfrac{1}{2}\norm[2]{\bfA w(\bfz) - \bfb}^2.
\end{align}
There are two key differences between our approach and traditional MRNSD. The first is that MRNSD performs all of the updates in $\bfx$-space, potentially yielding small step sizes near the boundary of the non-negative domain. In comparison, our approach will take steps in $\bfz$-space, which is unconstrained and can allow for larger steps. The second difference is the broader choice of mapping. MRNSD restricts itself to the exponential mapping. However, any other non-negativity mapping may be viable as well, e.g., $\bfx = w(\bfz) = \bfz^2$, or $\bfx = w(\bfz) = \abs{\bfz}$. Here, we propose the following,
    \begin{align} \label{eq:nnmapping}
        w_{a,c}(z) = \begin{cases}
            \e^{a(z - c)} - a(z-c) - 1, & \mbox{for } z > c,\\
            0,  & \mbox{for } z \le c,
        \end{cases}
    \end{align}
where $a > 0$ controls the steepness of the mapping and $c$ controls where the mapping switches from exponential to zero. The mapping is continuously differentiable and allows for the solution $\bfx = w(\bfz)$ to have entries equal to zero. An illustrative toy example of how different mappings perform is provided in \Cref{fig:mappings}.

This new method provides sparsity while ensuring non-negativity in $\bfx$, hence we refer to this method as Gradient Descent with a Sparsity-Promoting Non-negative Mapping (spNNGD). We summarize the proposed approach in~\Cref{alg:gdnn}. 

\begin{algorithm}
\caption{Gradient Descent with Sparsity-Promoting Non-negative Mapping (\texttt{spNNGD})}\label{alg:gdnn}
\begin{algorithmic}[1]
\State Solve $\min_{\bfz} \frac{1}{2}\|\bfA w(\bfz) - \bfb\|_2^2$
\While{not converged}
    \State $\bfg = \bfA^\top (\bfA w(\bfz) - \bfb)$
    \State $\bfs = -\diag{w'(\bfz)}\bfg$
    \State choose an optimal step size
        $$\alpha^* \in \min_{\alpha > 0} \ \tfrac{1}{2}\|\bfA w(\bfz + \alpha \bfs) - \bfb\|_2^2\hspace*{199pt}$$
    \State $\bfz \gets \bfz + \alpha^* \bfs$
\EndWhile
\end{algorithmic}
\end{algorithm}

This approach shares some similarities with sparsity-promoting MRNSD. First, the sparsity pattern will not change during the iterations; that is, once a zero appears in the search direction, we do not change the value of $\bfz$, so the zero remains throughout the iterations.   This is because $w'(\bfz)$ will have the same pattern of zeros as $w(\bfz)$. Second, the implicit vertical shift defined by $\bfvarepsilon$ plays a similar role to $a$ and $c$ in that it enables zero values to exist in the mapping.


\section{Numerical Experiments}
\label{sec:numericalexperiments}

\subsection{Experiment Setup}

In the following, we discuss four experiments comparing the effectiveness of our two proposed methods.  Specifically, we investigate an image deblurring, an image completion, a medical tomography, and a super-resolution problem. For reproducible science, we made the experiments and code available at \url{https://github.com/elizabethnewman/DL4IP}.

\paragraph{Metrics and Presented Results}
We consider three metrics that measure (relative) accuracy and sparsity, respectively, given by
    \begin{align}
        \text{rel. residual} &= \frac{\|\bfA \bfx - \bfb\|_2}{\|\bfb\|_2}, &
        \text{rel. error} &= \frac{\|\bfy_{\rm true} - \bfG \bfx\|_2}{\|\bfy_{\rm true}\|_2}, & 
        &\mbox{and} &
        \text{rel. sparsity} &= \frac{\texttt{nnz}(\bfx)}{\texttt{nnz}(\bfy_{\rm true})},
    \end{align}
where $\texttt{nnz}(\cdot)$ computes the number of nonzero entries for the given array. For all metrics, lower values indicate better performance.  While the relative residual and the relative error are common measures of quality, the relative sparsity metric, storing the cost of the coefficients $\bfx$ relative to storing the original image, $\bfy_{\rm true}$, is a critical measure in dictionary learning. The image's storage cost is reduced when the relative sparsity is below $1$. For every experiment and each algorithm, we present the reconstruction results using the above metrics and the convergence results, including the relative residual, a proxy for the relative sparsity, and the optimal step sizes for each iteration.

\begin{figure}
    \centering
    \subfloat[Model operators $\bfC$ where $\texttt{numel}(\cdot)$ returns the total number of elements in the given array. \label{fig:operators_C}]{\begin{tikzpicture}
        \def\s{0.2}
        \def\b{2}
        \def\sep{0.75}

        \node[anchor=east] at (-1.25,-\b) {$\texttt{size}(\bfC) =$};
        \node[anchor=east] at (-1.25,-\b-0.6) {$\frac{\texttt{nnz}(\bfC)}{\texttt{numel}(\bfC)} \approx $};

        \node[draw] (n0) at (0,0) {\includegraphics[width=0.15\linewidth]{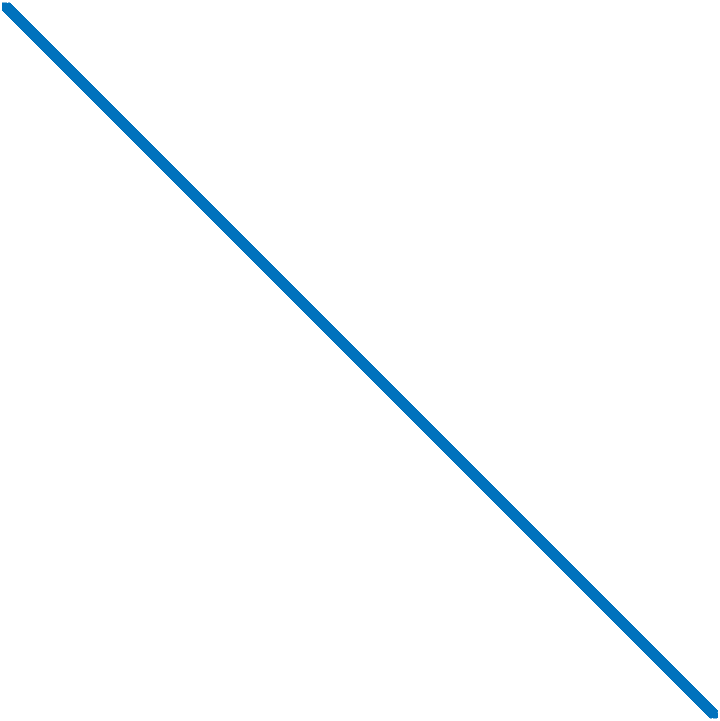}};

        \path let \p1 = (n0) in node  at (\x1,-\b) {$65,\!536 \times 65,\!536$};
        \path let \p1 = (n0) in node  at (\x1,-\b-0.6) {$3.8\times 10^{-4}$ \vphantom{$\frac{\texttt{nnz}(\bfC)}{\texttt{numel}(\bfC)} \approx$}};

        \node[above=0.0cm of n0.north, anchor=south]{Deblurring};

        \node[right=\sep cm of n0.north east, anchor=north west, draw] (n1) {\includegraphics[width=0.075\linewidth]{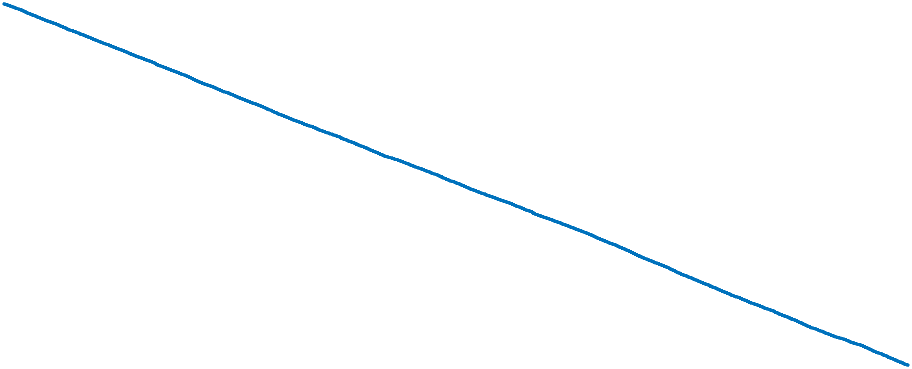}};

        \path let \p1 = (n1) in node  at (\x1,-\b) {$6,\!554 \times 16,\!384$};
        \path let \p1 = (n1) in node  at (\x1,-\b-0.6) {$6.1\times 10^{-5}$ \vphantom{$\frac{\texttt{nnz}(\bfC)}{\texttt{numel}(\bfC)} \approx$}};

        \node[above=0.0cm of n1.north, anchor=south]{Dead Pixel};

        \node[right=\sep cm of n1.north east, anchor=north west, draw] (n2) {\includegraphics[width=0.15\linewidth]{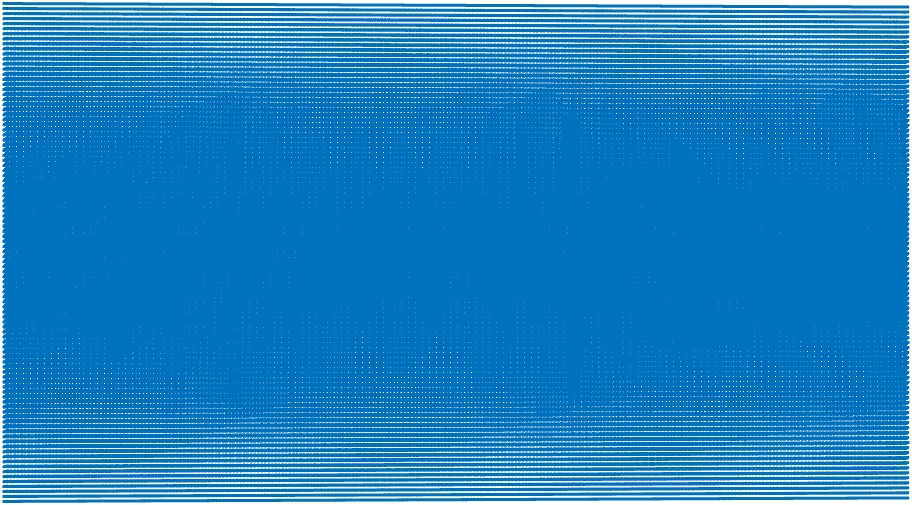}};

        \path let \p1 = (n2) in node  at (\x1,-\b) {$36,\!200 \times 65,\!536$};
        \path let \p1 = (n2) in node  at (\x1,-\b-0.6) {$3.5\times 10^{-3}$ \vphantom{$\frac{\texttt{nnz}(\bfC)}{\texttt{numel}(\bfC)} \approx$}};

        \node[above=0.0cm of n2.north, anchor=south]{Tomography};

        \node[right=\sep cm of n2.north east, anchor=north west, draw] (n3) {\includegraphics[width=0.3\linewidth]{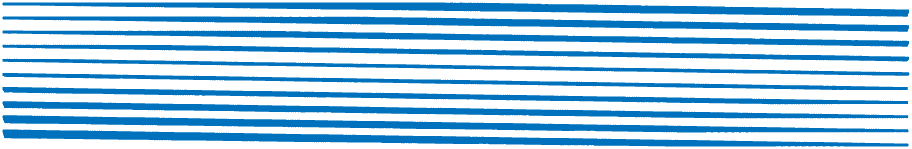}};

        \path let \p1 = (n3) in node  at (\x1,-\b) {$40,\!960 \times 262,\!144$};
        \path let \p1 = (n3) in node  at (\x1,-\b-0.6) {$2.5\times 10^{-4}$ \vphantom{$\frac{\texttt{nnz}(\bfC)}{\texttt{numel}(\bfC)} \approx$}};

        \node[above=0.0cm of n3.north, anchor=south]{Superresolution};

    \end{tikzpicture}}

    \bigskip
    
    \subfloat[Finite differencing regularization operators $\bfL$ \label{fig:operators_L}]{\begin{tikzpicture}
        \node[draw] (n0) at (0,0) {\includegraphics[width=0.1\linewidth]{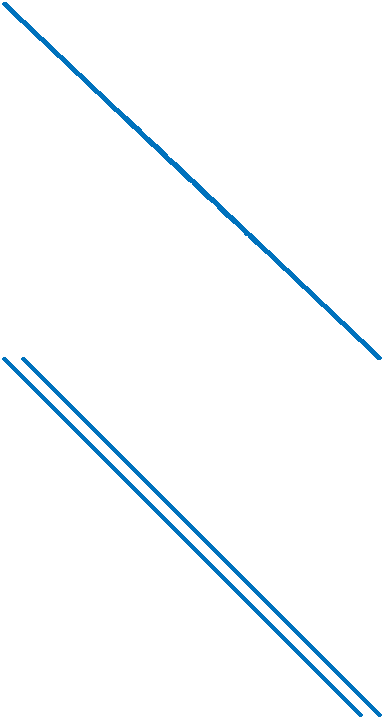}};
        \node[above=0.0cm of n0.north, anchor=south] (fd) {Pixel Smoother};

        \node[below=0.0cm of n0.south, anchor=north] {\scriptsize $[(N-1)M + M(N-1)] \times MN$};

        \node[draw, right=2.0cm of n0.north east, anchor=north west] (n2) {\includegraphics[width=0.1\linewidth]{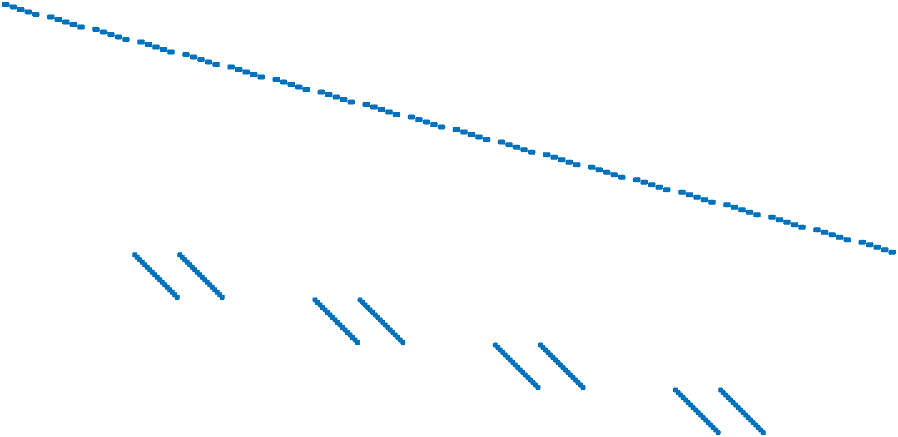}}; 
        
        \node (ps2) at (n2 |- fd) {Patch Border Smoother};

        \node[below=0.0cm of n2.south, anchor=north] {\scriptsize $[(\frac{M}{p}(\frac{N}{q} - 1) + \frac{N}{q}(\frac{M}{p} - 1)] \times MN$};

        \node[draw, right=2.0cm of n2.north east, anchor=north west] (n1) {\includegraphics[width=0.1\linewidth]{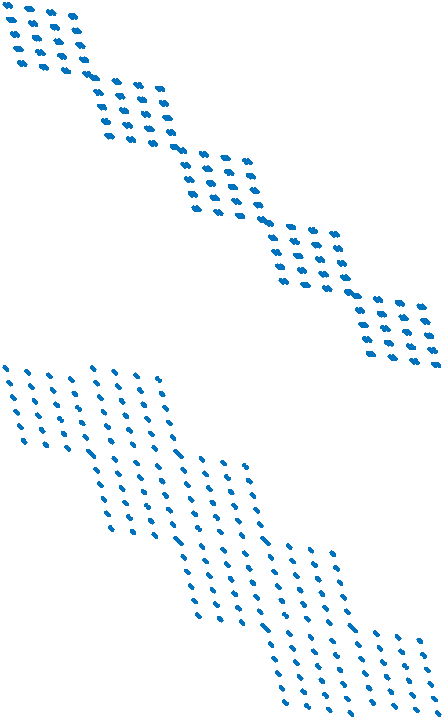}};
        
        \node (ps) at (n1 |- fd) {Patch Smoother};

        \node[below=0.0cm of n1.south, anchor=north] {\scriptsize $[pq((\frac{M}{p}(\frac{N}{q} - 1) + \frac{N}{q}(\frac{M}{p} - 1))] \times MN$};

    \end{tikzpicture}}
    \caption{Sparsity patterns of operators. The operators are illustrated with the correct proportions based on their row and column dimensions. The width of the operators is correctly scaled (e.g., the superresolution operator has twice as many columns as the deblurring operator).}
    \label{fig:operators}
\end{figure}
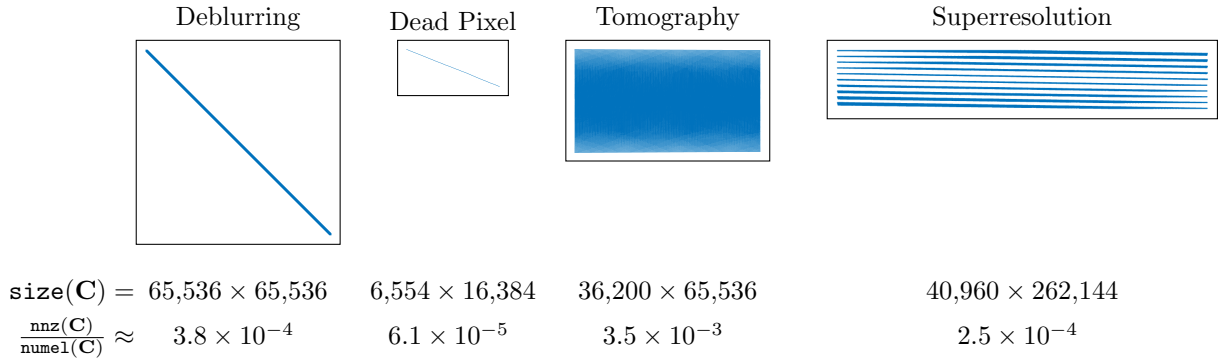
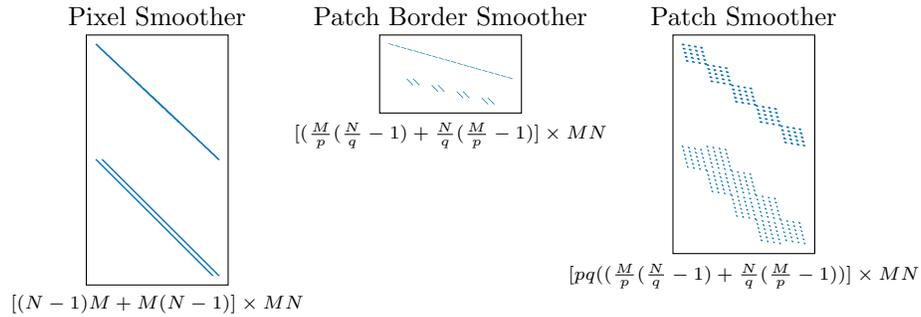

\paragraph{Operators.} 
Each experiment in this section relies on a specific linear operator $\bfC$. The shapes and sparsity patterns vary significantly for each application. For completeness, we illustrate the sparsity patterns of the operators in~\Cref{fig:operators_C}.  In every problem, we incorporate a generalized Tikhonov regularization on the approximation $\mu\calR(\bfy) = \mu\norm[2]{\bfL\bfy}^2$.  We utilize different regularization operators $\bfL$ depending on the problem.  All of the operators we consider promote smoother approximations, and hence are variations of two-dimensional finite difference operators.  By a ``finite difference'' operator, we mean that each row of $\bfL$ contains a single $1$ and a single $-1$ to compute the difference between pixels.  The three finite differencing operators we consider are the following: the difference between (1) all adjacent pixels, (2) adjacent pixels at the boundary of patches, and (3) adjacent patches.  The sparsity patterns of the three operators are depicted in~\Cref{fig:operators_L}. 

\subsection{Deblurring}
\label{sec:deblur}
\begin{figure}
    \centering
    
    \begin{tikzpicture}
    \def\w{0.2}
    
    \node (orig) at (0,0) {\includegraphics[width=\w\linewidth]{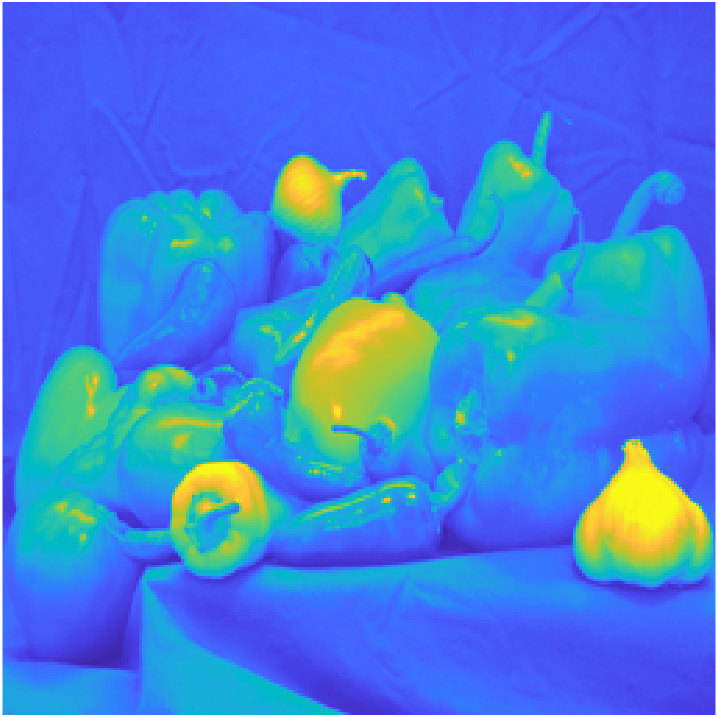}};
    \node[above=0.0cm of orig.north, anchor=south] {\begin{tabular}{c}Original, $\bfy_{\rm true}$\\ $256 \times 256$\end{tabular}};

    \node[right=0.0cm of orig.east, anchor=west] (orig) {\includegraphics[width=\w\linewidth]{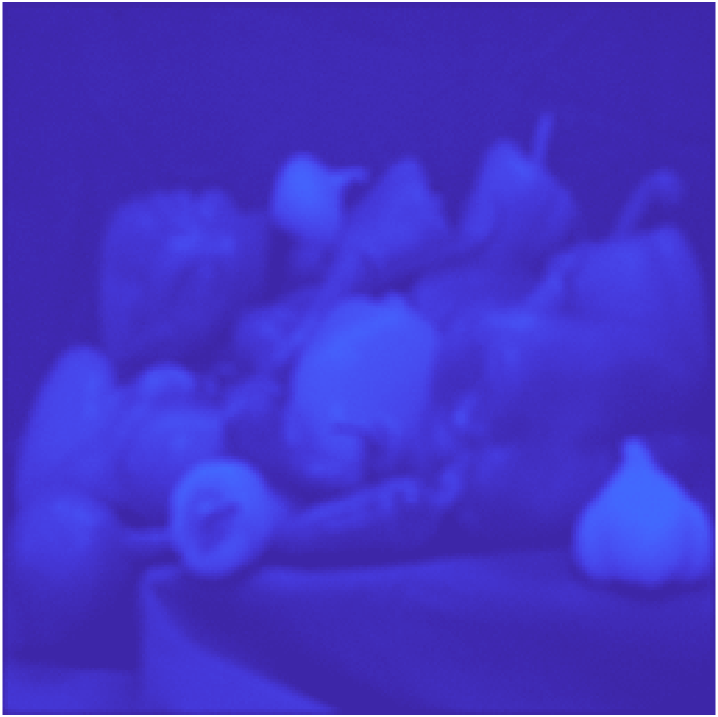}};
    \node[above=0.0cm of orig.north, anchor=south] {\begin{tabular}{c}Blurred, Noisy, $\bfb$ \\ $256\times 256$ \end{tabular}};

    \node[right=0.0cm of orig.east, anchor=west] (recon_mrnsd){\includegraphics[width=\w\linewidth]{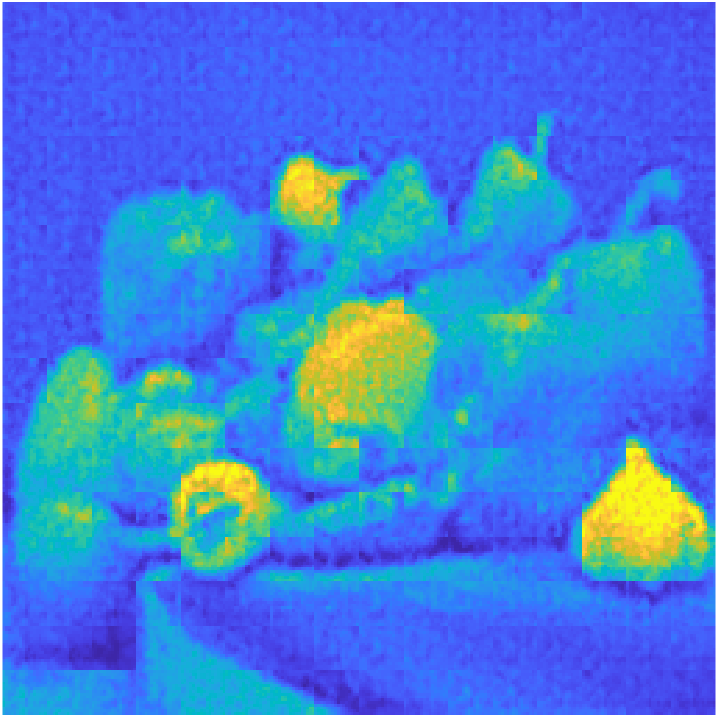}};
    \node[above=0.0cm of recon_mrnsd.north, anchor=south] {\begin{tabular}{c} \texttt{spMRNSD}, $\bfG\bfx$ \\ $\lambda = 10^{-8}$ \end{tabular}};

    \node[below=0.0cm of recon_mrnsd.south, anchor=north] (diff_mrnsd){\includegraphics[width=\w\linewidth]{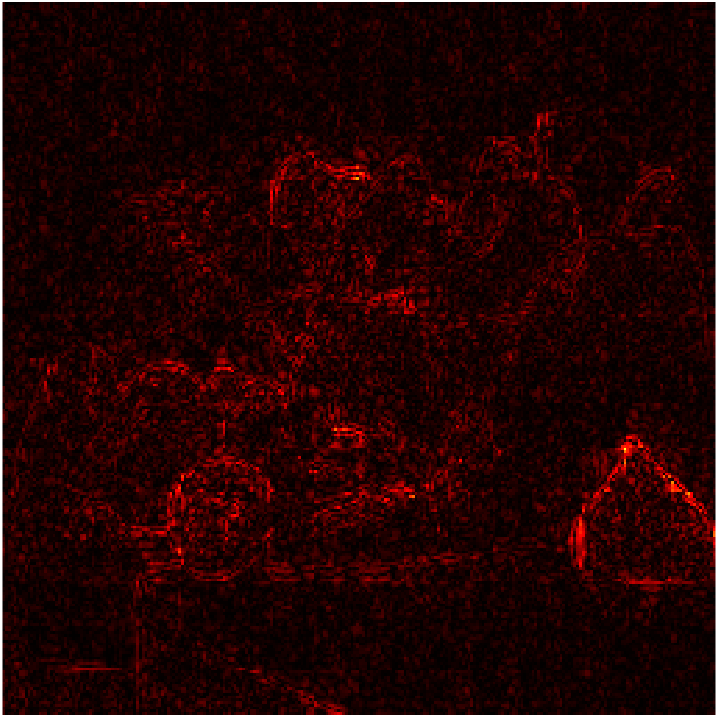}};

    \node[right=0.0cm of recon_mrnsd.east, anchor=west] (recon_gdnn){\includegraphics[width=\w\linewidth]{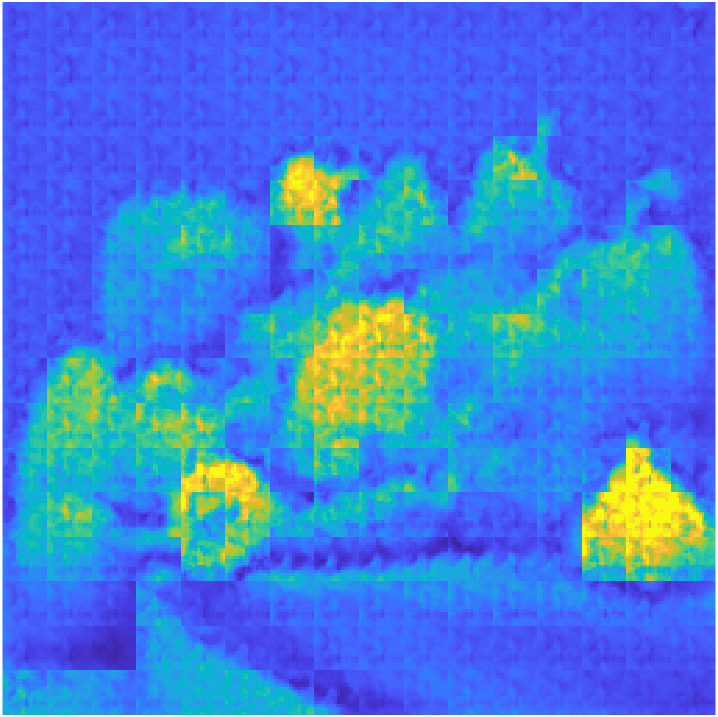}};
    \node[above=0.0cm of recon_gdnn.north, anchor=south] {\begin{tabular}{c}\texttt{spNNGD},  $\bfG\bfx$ \\ $a = 0.1$, $c = -0.75$ \end{tabular}};

    \node[below=0.0cm of recon_gdnn.south, anchor=north] (diff_gdnn){\includegraphics[width=\w\linewidth]{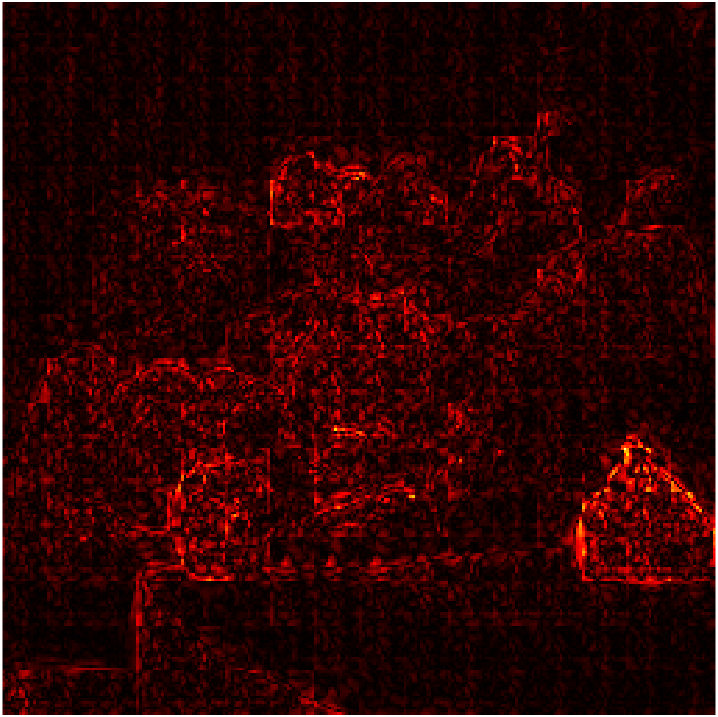}};

    \node[right=0.0cm of recon_gdnn.east, anchor=west] (parula) {\includegraphics[height=\w\linewidth]{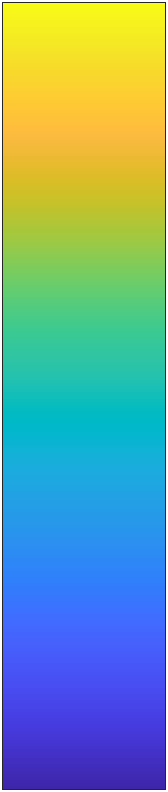}};
    \node[below right=0.0cm of parula.north east, anchor=north west] {$0.011$};
    \node[right=0.0cm of parula.south east, anchor=south west] {$0$};
    
    \node[above right=0.0cm of diff_gdnn.east, anchor=west] (hot) {\includegraphics[height=\w\linewidth]{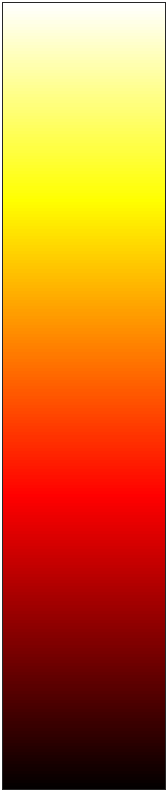}};
    \node[below right=0.0cm of hot.north east, anchor=north west] {$0.007$};
    \node[above right=0.0cm of hot.south east, anchor=south west] {$0$};

    \node[below=0.0cm of diff_mrnsd.south, anchor=north] (err_mrnsd) {\vphantom{$\frac{\|\bfy_{\rm true} - {\bfG\bfx}\|_2}{\|\bfy_{\rm true}\|_2} \approx$}
    $\bf 1.15\times 10^{-1}$};
    \node[below=0.1cm of err_mrnsd.south, anchor=north] (sparse_mrnsd) {\vphantom{$\frac{\texttt{nnz(coefficients)}}{\texttt{nnz(image)}} \approx$}$\bf 0.15$};

    \node[below=0.0cm of diff_gdnn.south, anchor=north] (err_gdnn) {\vphantom{$\frac{\|\bfy_{\rm true} - {\bfG\bfx}\|_2}{\|\bfy_{\rm true}\|_2} \approx$ }$1.57\times 10^{-1}$};
    \node[below=0.1cm of err_gdnn.south, anchor=north] (sparse_gdnn) {\vphantom{$\frac{\texttt{nnz(coefficients)}}{\texttt{nnz(image)}} \approx$}$0.23$};

    \node[above=0.0cm of diff_mrnsd.west, anchor=south, rotate=90] {\begin{tabular}{c} absolute difference \\ $|\bfy_{\rm true} - \bfG\bfx|$ \end{tabular}};

    \node[left=1.0cm of err_mrnsd.west, anchor=east]  (err) {relative error, $\frac{\|\bfy_{\rm true} - {\bfG\bfx}\|_2}{\|\bfy_{\rm true}\|_2} \approx$};

    \node[below=0.1cm of err.south east, anchor=north east] {relative sparsity, $\frac{\texttt{nnz}(\bfx)}{\texttt{nnz}(\bfy_{\rm true})} \approx$};

    \node[draw, below=0.0cm of sparse_mrnsd.south, anchor=north] (sol_mrnsd) {\includegraphics[scale=0.2]{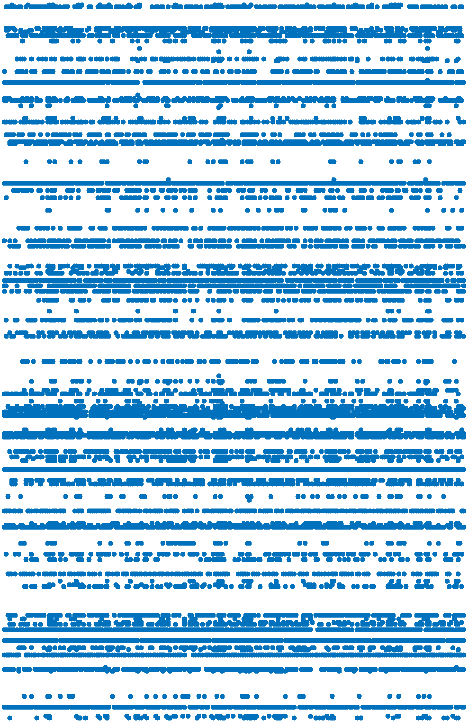}};

    \node[draw, below=0.0cm of sparse_gdnn.south, anchor=north] {\includegraphics[scale=0.2]{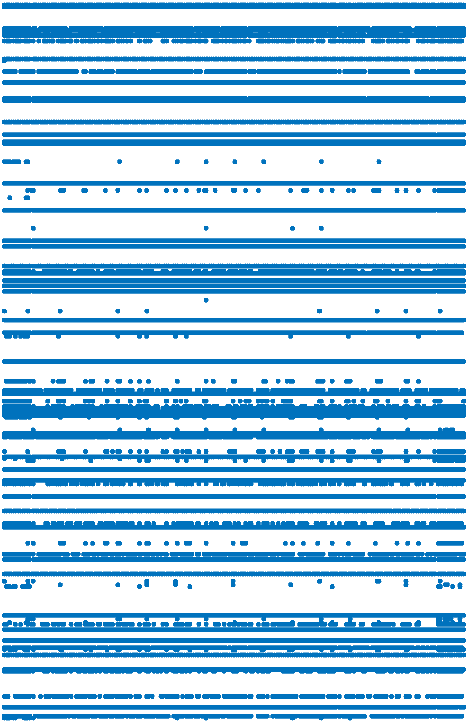}};

    \node[anchor=east] at (err.east |- sol_mrnsd) {matricized coefficients, $\bfx \equiv $};

    \end{tikzpicture}

\caption{Reconstructions for the image deblurring task. The top row (left-to-right) contains the original image, the blurred, noisy image, the reconstruction using \texttt{spMRNSD}, and the reconstruction using \texttt{spNNGD}. The middle row displays the absolute differences between the original and the reconstructions with the relative error and sparsity metrics below. Here, \texttt{spMRNSD} produces the best approximation in terms of both metrics of interest. The bottom row depicts the sparsity patterns of matricized solutions $\bfx$ of size $400 \times 256$ where each column corresponds to the coefficients for one vectorized dictionary patch.}
\label{fig:deblur_recon}
\end{figure}
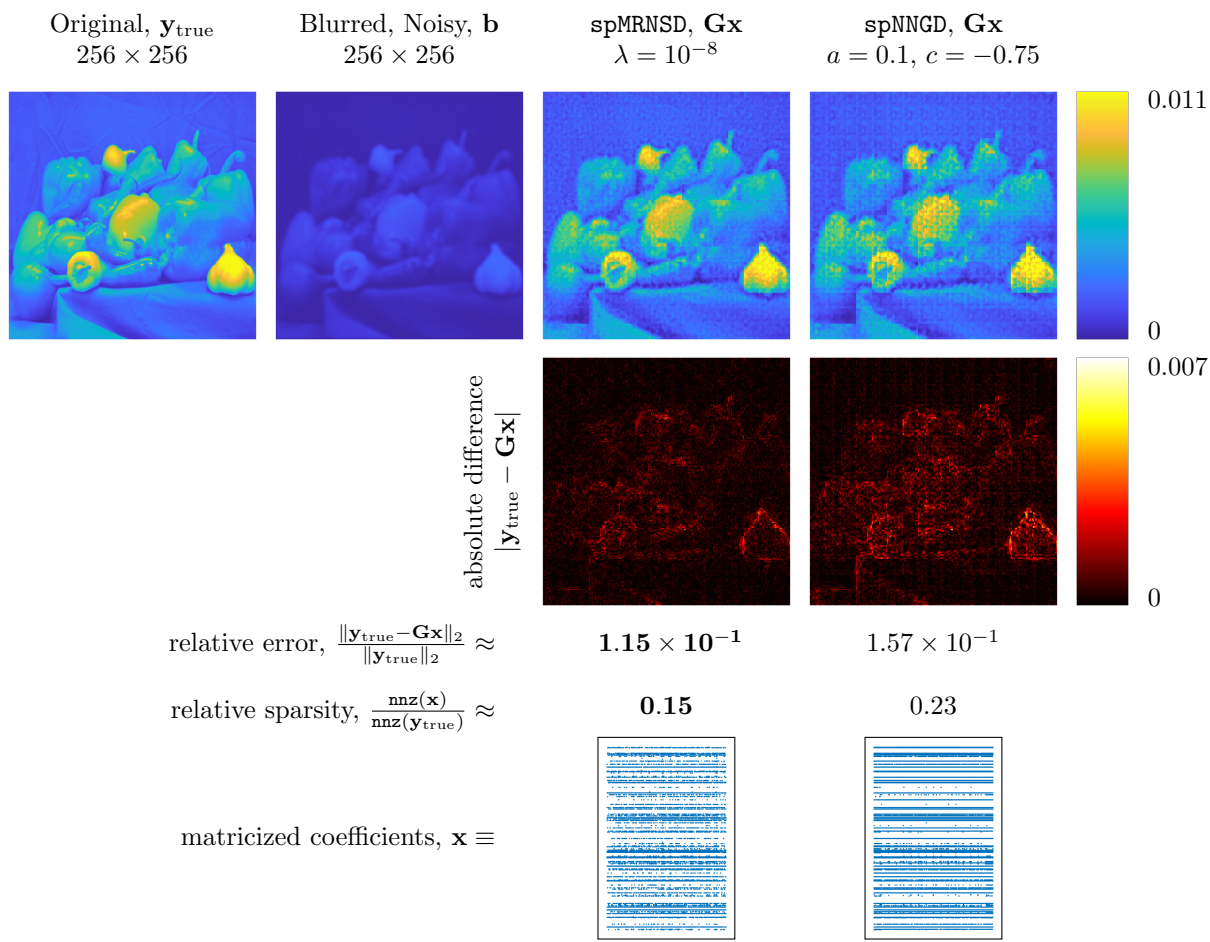

\begin{figure}
    \centering
    \subfloat[Relative Residual \label{fig:deblur_relerr}]{\includegraphics[width=0.32\linewidth]{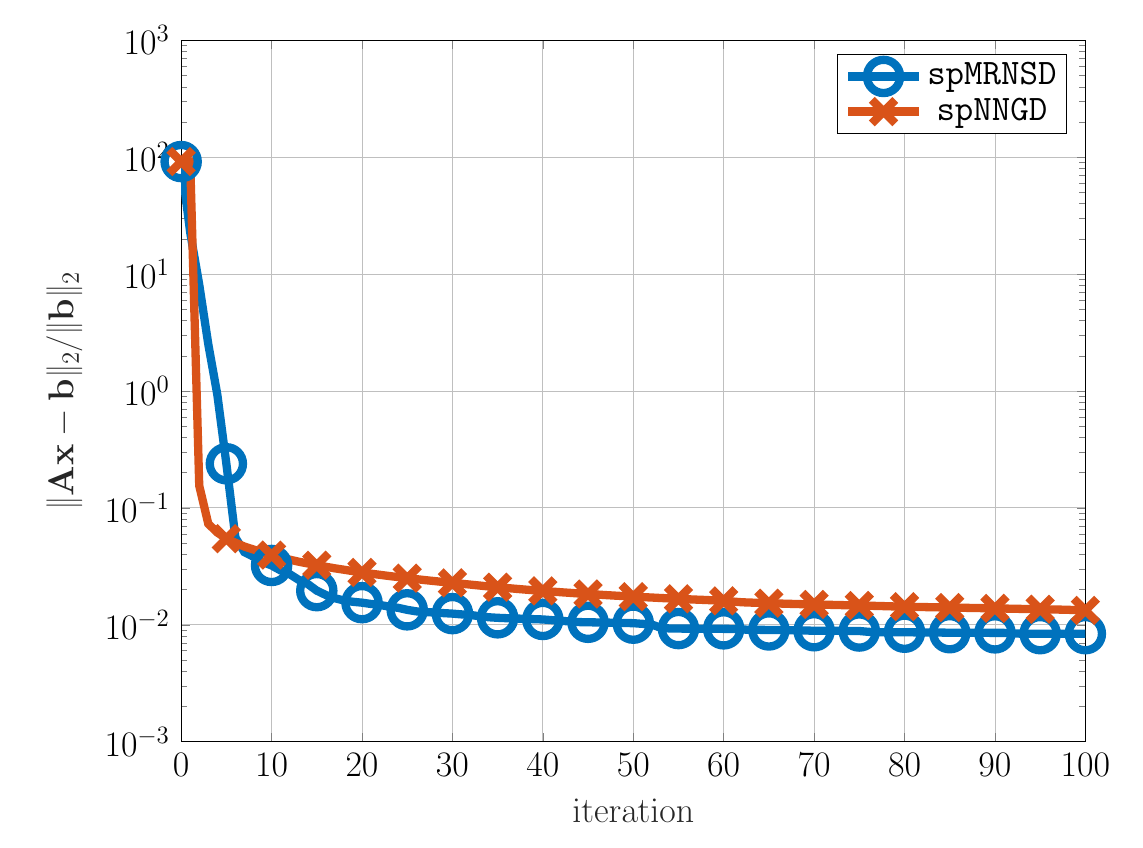}}
    \subfloat[Proxy for Relative Sparsity\label{fig:deblur_sparsity}]{\includegraphics[width=0.32\linewidth]{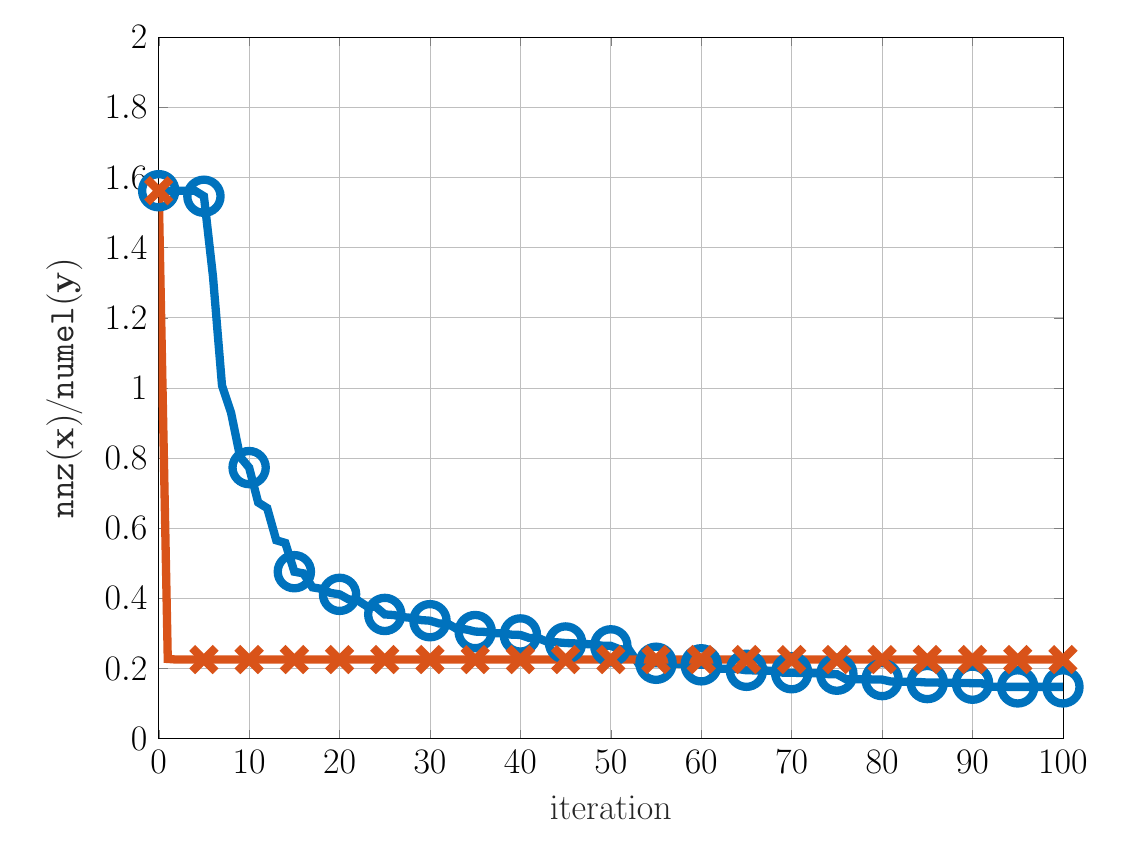}}
    \subfloat[Optimal Step Sizes\label{fig:deblur_alpha}]{\includegraphics[width=0.32\linewidth]{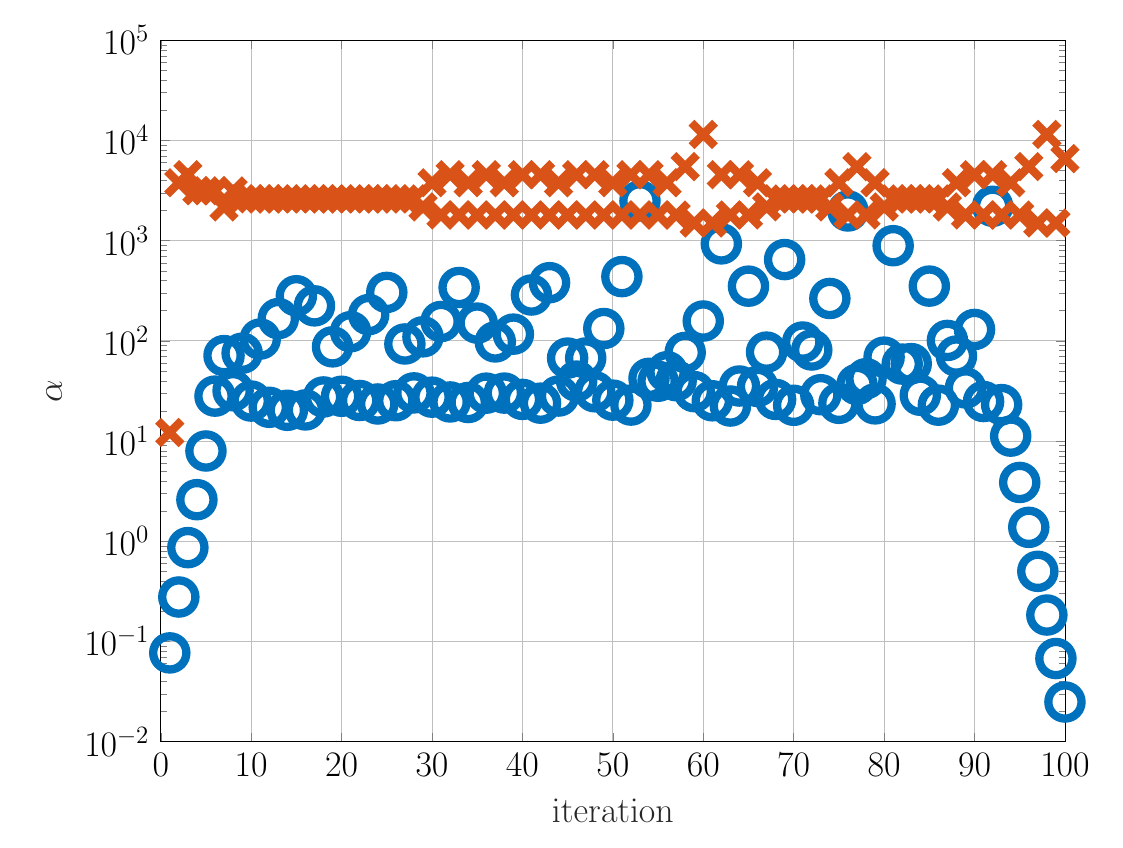}}

    \caption{Convergence behavior for the deblurring problem. The relative residual norm (left) converges at similar rates for both  algorithms.   The proxy for relative sparsity (middle) converges faster for \texttt{spNNGD}.  On average, the step sizes (right) for \texttt{spNNGD} are larger and more consistent.}
    \label{fig:deblur_convergence}
\end{figure}

Capturing digital images can often introduce noise, blur, and other distortions that hide the true representation. 
Image deblurring has been a long-studied technique to remove unwanted artifacts in and enhance digital images~\cite{HansenNagyOLeary2006:deblur}. 
Here, we consider applying the Gaussian blurring operator to an $N\times N$ image. 
To be precise, we construct the blurring operator $\bfC \in \Rbb^{N^2\times N^2}$ as a banded, Toeplitz matrix with exponentially decaying bands\footnote{In {\sc Matlab}, we write \texttt{ C1 = toeplitz(exp(-(0:b-1) / (2 * sigma * sigma)), zeros(1, n - b)); C = kron(C1,C1); } to construct this blurring operator.}. 
We use a bandwidth of $b = 3$ and an exponential decay rate of $\sigma = 4$. 
Given a vectorized true image $\bfy_{\rm true} \in \Rbb_+^{N^2}$, we construct a blurred, noisy image $\bfb =\bfC \bfy_{\rm true} + \bfeta$ where each entry of the noise vector $\bfeta$ is drawn from a Gaussian distribution $\Ncal({\bf0}, \beta\|\bfC\bfy_{\rm true} \|_2^2\bfI_{N^2})$ with noise level $\beta = 10^{-4}$.  
We include a pixel smoothing Tikhonov regularizer term with regularization parameter $\mu = 10^{-4}$. 
We train for $100$ iterations for both \texttt{spNNGD} and \texttt{spMRNSD} and report the results in~\Cref{fig:deblur_recon} and the convergence behavior in~\Cref{fig:deblur_convergence}.

In~\Cref{fig:deblur_recon}, we observe that both \texttt{spNNGD} and \texttt{spMRNSD} produce qualitatively good approximations to the original image for similar sparsity levels. Both sets of coefficients require less than $25\%$ of the storage cost of the original image. Both approximations produce some blocky artifacts due to the dictionary patch representation; this is most noticeable in the \texttt{spNNGD} reconstruction of small vegetables, such as the garlic in the bottom right. In the bottom row, we observe similar sparsity patterns in the matricized coefficients for \texttt{spMRNSD} and \texttt{spNNGD}.  In particular, the two have some rows with many nonzeros and some sparse rows, and the rows with high density are similar for both solutions. It can be observed, that the solution obtained from \texttt{spMRNSD} is slightly sparser than that obtained from \texttt{spNNGD}.

The convergence behavior in~\Cref{fig:deblur_relerr} shows that \texttt{spNNGD} and \texttt{spMRNSD} converge to similar levels of accuracy at similar rates.  In comparison, the sparsity convergence plot~\Cref{fig:deblur_sparsity}, we see that \texttt{spNNGD} converges quickly to a sparse solution and then stagnates, whereas \texttt{spMRNSD} converges to a sparse solution more gradually.  This behavior indicates that \texttt{spNNGD} fixes a sparsity pattern early and learns the best coefficients for that pattern while \texttt{spMRNSD} adjusts the sparsity pattern during iterations.  We note that both methods are constrained by the current sparsity pattern because the search direction $\bfs = \diag{w'(\bfz)}\bfg$ has the same zero entries as the approximate solution $\bfx = w(\bfz)$.   The step size plot in~\Cref{fig:deblur_alpha} supports our original hypothesis that we take larger step sizes in $\bfz$-space (\texttt{spNNGD}) than in $\bfx$-space (\texttt{spMRNSD}).

\begin{figure}
    \centering

    \subfloat[\texttt{spNNGD}: the effect of steepness $a$ and shift $c$ in $w_{a,c}(z)$ on relative error and sparsity. \label{fig:deblur_param_gdnn}]{\begin{tikzpicture}

        \node at (0,0) (n0) {\includegraphics[height=0.225\linewidth]{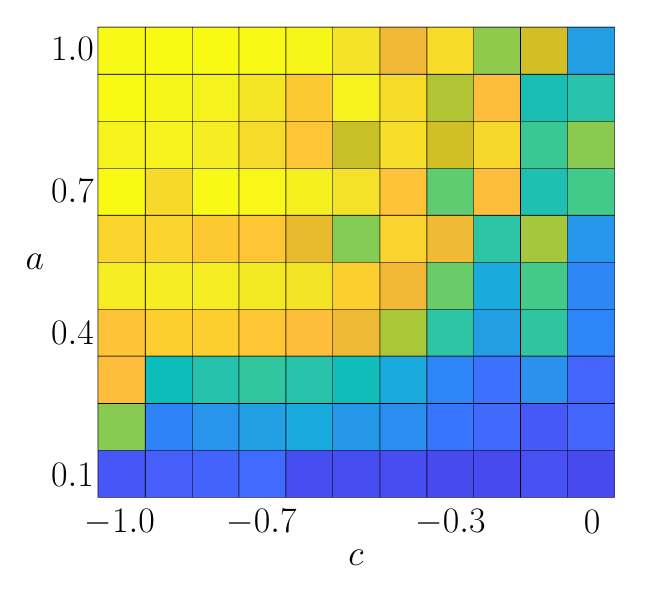}};
        \node[above=0.0cm  of n0.north, anchor=south] (tmp1) {Relative Error};

        \node[right=0.0cm of n0.east, anchor=west] (n1) {\includegraphics[height=0.225\linewidth]{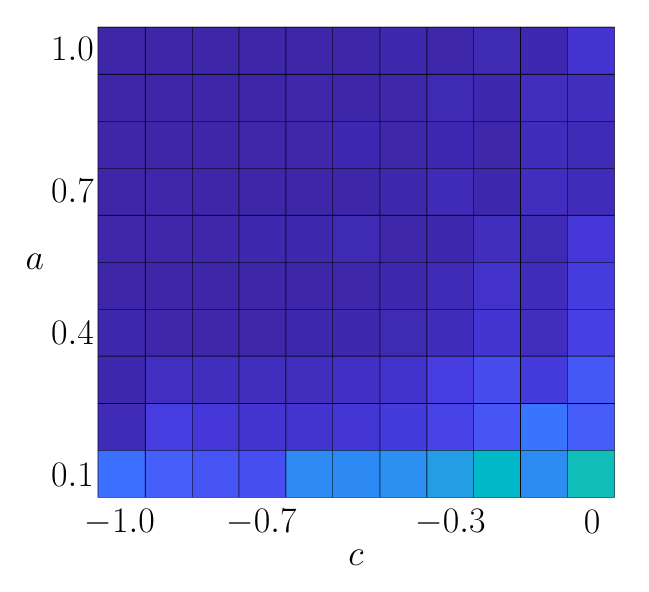}};
        \node[above=0.0cm  of n1.north, anchor=south] (tmp2) {Relative Sparsity};

        \node[right=0.0cm of n1.east, anchor=west] (parula) {\includegraphics[height=0.15\linewidth]{parula.png}};

        \node[below right=0.0cm of parula.north east, anchor=north west] {$1$};
        \node[above right=0.0cm of parula.south east, anchor=south west] {$0$};
        
    \end{tikzpicture}
    }
    \hspace{0.25cm}
    \subfloat[\texttt{spMRNSD}: the effect of the regularization parameter $\lambda$ on relative error and sparsity. \label{fig:deblur_param_mrnsd}]{
    \begin{tikzpicture}
        \node at (0,0) (n0) {\includegraphics[height=0.25\linewidth]{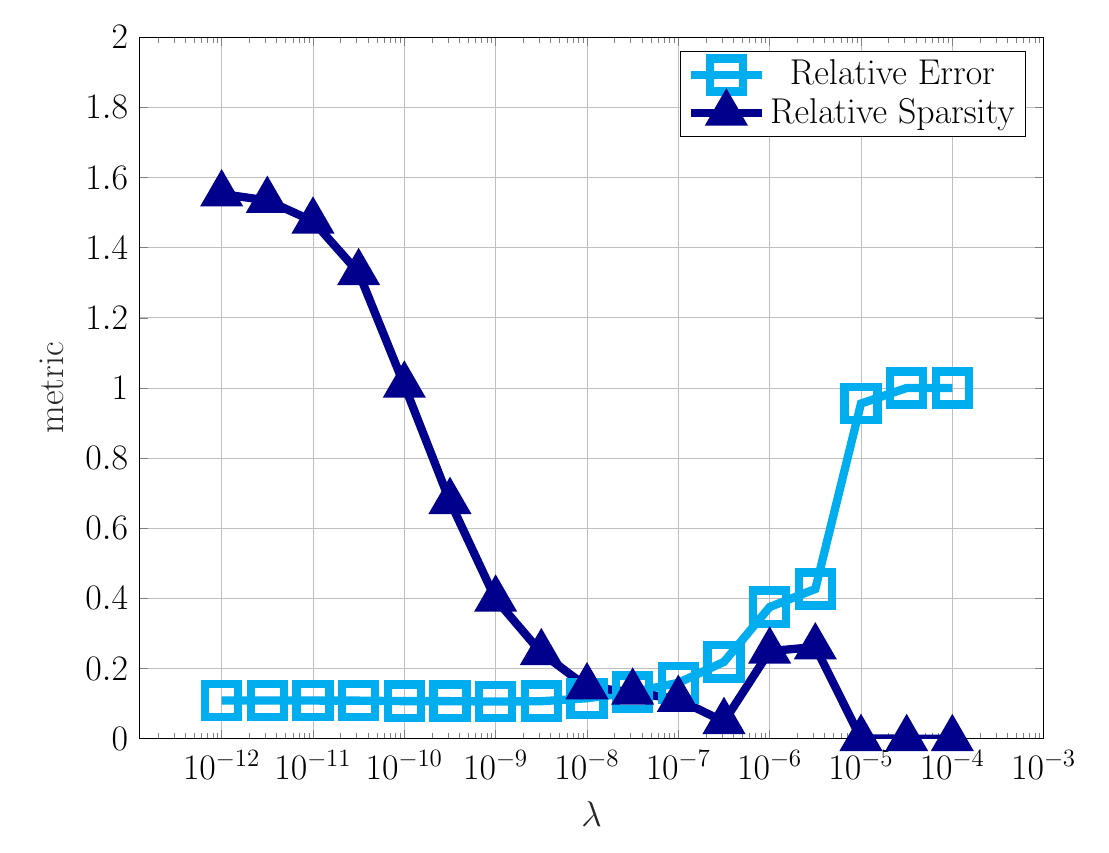}};
    \end{tikzpicture}
    }
    \caption{Parameter sensitivity exploration for the deblurring problem.  
    }
    \label{fig:deblur_paramSensitivity}
\end{figure}

Examining the relative error plot, it becomes apparent that the relative error of \texttt{spNNGD} is more susceptible to changes in the steepness parameter $a$ rather than the shift parameter $c$. This observation is substantiated by the relatively constant rows in the relative error plot. We also observe a positive correlation between $a$ and $c$, implying that increasing both parameters simultaneously results in small changes to the relative error. However, while the relative error exhibits a degree of sensitivity to the choice of parameters,  the sparsity plot appears to be less sensitive to its parameter selection.

For \texttt{spMRNSD}, we observe that the relative error remains mainly constant across a wide range of regularization parameter $\lambda$ values, while the sparsity undergoes drastic changes. This behavior is consistent with expectations, that larger values for $\lambda$ will lead to poorer approximations and sparser solutions.

There are several major takeaways from our exploration of the image deblurring problem. First, the \texttt{spMRNSD} produced better results in terms of approximation quality and sparsity, but \texttt{spNNGD} is competitive. In terms of convergence, \texttt{spNNGD} obtains a sparser solution in fewer iterations than \texttt{spMRNSD} with larger step sizes, but as a result, prescribes a sparsity pattern of the solution early. In comparison, \texttt{spMRNSD} gradually increases sparsity which can bring benefits in the long run.  In terms of sensitivity to parameters, the \texttt{spNNGD} parameters can drastically change the relative error, but the sparsity stays more consistent.   Conversely, the regularization parameter for \texttt{spMRNSD} can drastically change the sparsity but does not impact the relative error as dramatically.

\subsection{Image Completion}
\label{sec:indicator}

Matrix completion arises in numerous fields~\cite{candes2008exact, li2019survey}, such as collaborative filtering to predict users preferences with incomplete information~\cite{candes2009matrix}, pixels~\cite{Bertalmio2000:inpainting},  recovery of incomplete signals~\cite{Du2013:traffic}, and digital image inpainting to restore damaged or missing information. There are similarly many different approaches to completing matrices, including forming a low-rank approximation and solving an inverse problem. In this section, we explore the inverse problem approach for digital inpainting applications. 

\begin{figure}
    \centering
    
    \begin{tikzpicture}
    \def\w{0.2}
    
    \node (orig) at (0,0) {\includegraphics[width=\w\linewidth]{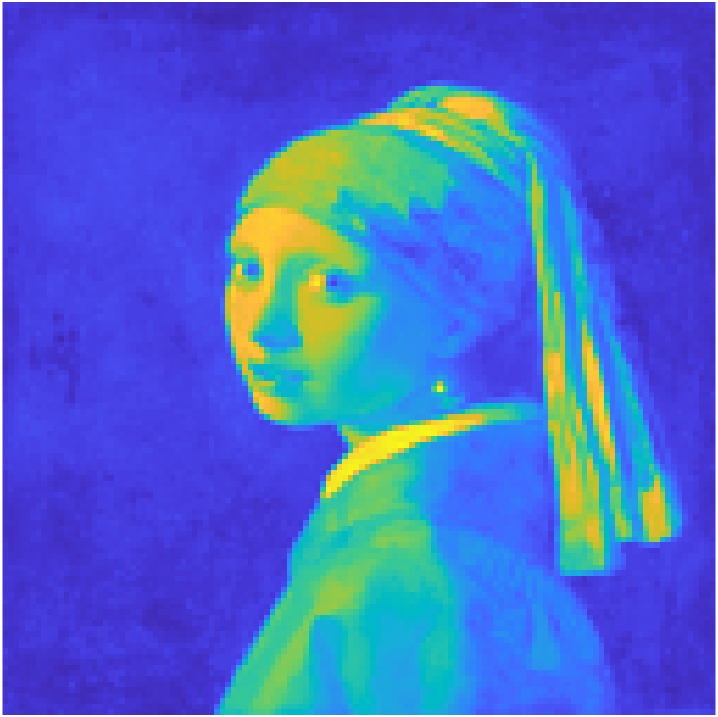}};
    \node[above=0.0cm of orig.north, anchor=south] {\begin{tabular}{c} Original, $\bfy_{\rm true}$ \\ $128 \times 128$ \end{tabular}};

    \node[right=0.0cm of orig.east, anchor=west] (orig){\includegraphics[width=\w\linewidth]{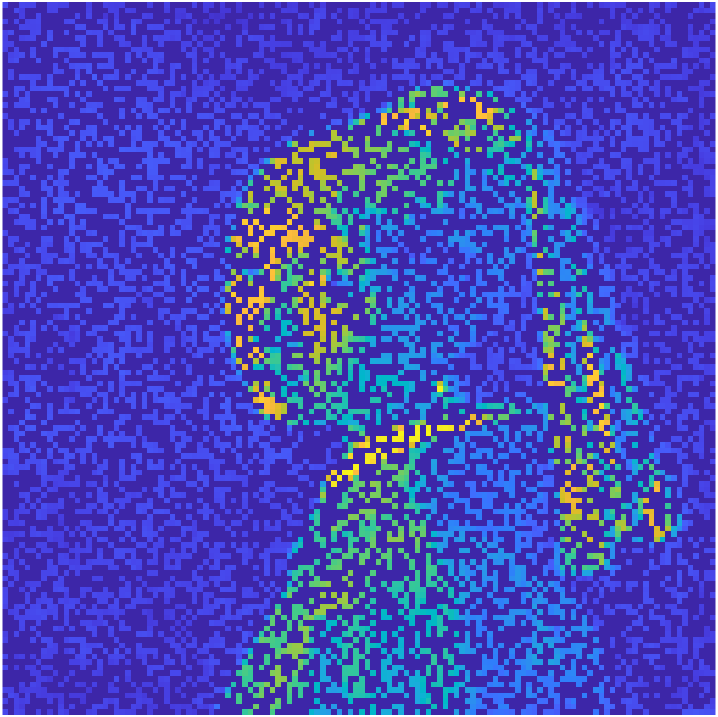}};
    \node[above=0.0cm of orig.north, anchor=south] {\begin{tabular}{c} Missing Pixels, $\bfb$ \\ $60\%$ removed \end{tabular}};

    \node[right=0.0cm of orig.east, anchor=west] (recon_mrnsd){\includegraphics[width=\w\linewidth]{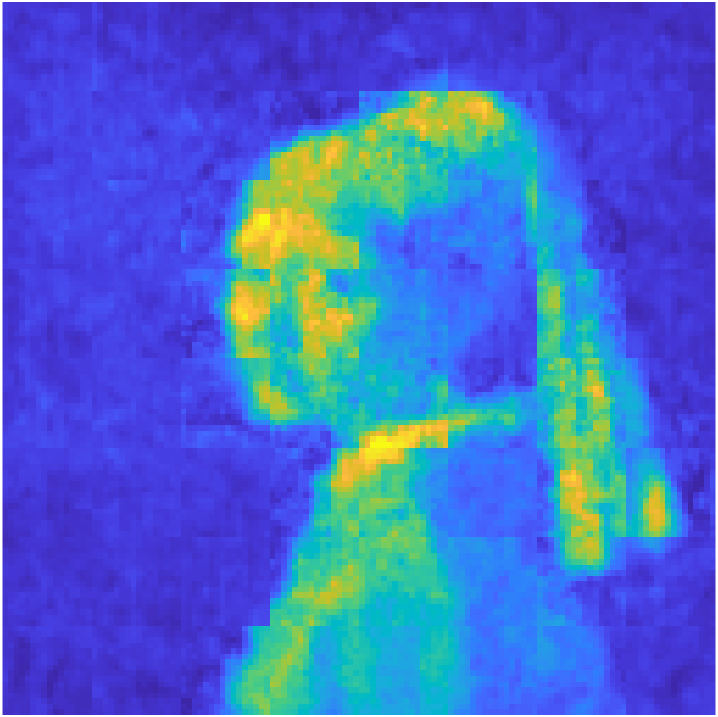}};
    \node[above=0.0cm of recon_mrnsd.north, anchor=south] {\begin{tabular}{c} \texttt{spMRNSD}, $\bfG\bfx$ \\ $\lambda = 10^{-3}$\end{tabular}};

    \node[below=0.0cm of recon_mrnsd.south, anchor=north] (diff_mrnsd){\includegraphics[width=\w\linewidth]{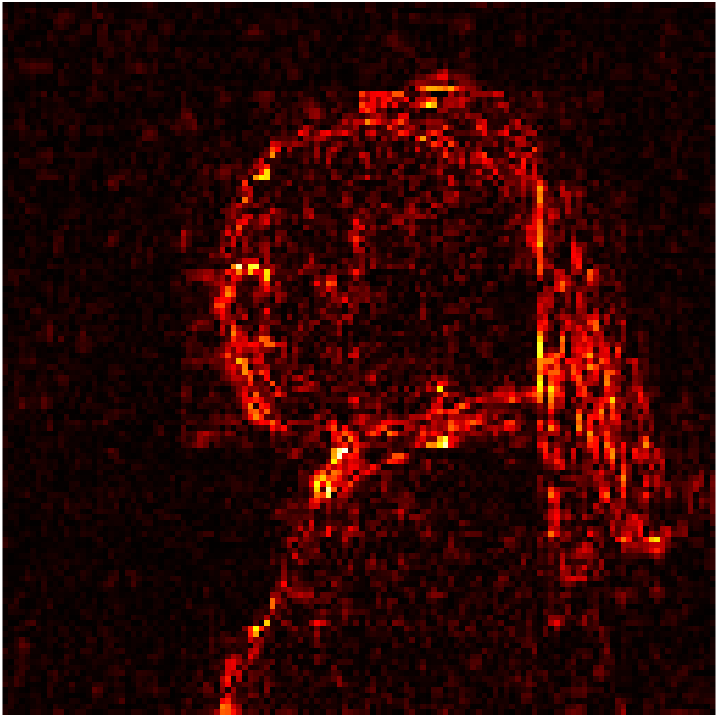}};

    \node[right=0.0cm of recon_mrnsd.east, anchor=west] (recon_gdnn){\includegraphics[width=\w\linewidth]{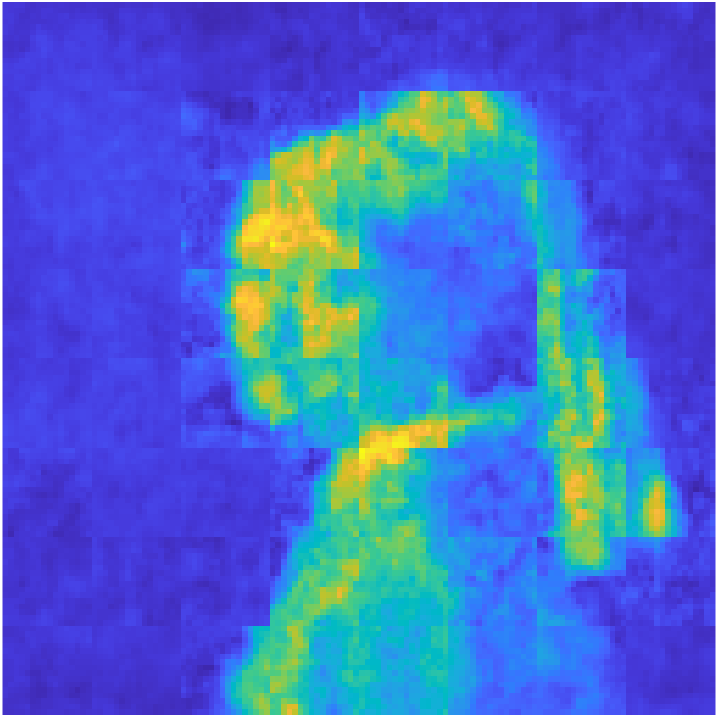}};
    \node[above=0.0cm of recon_gdnn.north, anchor=south] {\begin{tabular}{c}\texttt{spNNGD}, $\bfG\bfx$ \\ $a = 1$, $c = -0.5$ \end{tabular}};

    \node[below=0.0cm of recon_gdnn.south, anchor=north] (diff_gdnn){\includegraphics[width=\w\linewidth]{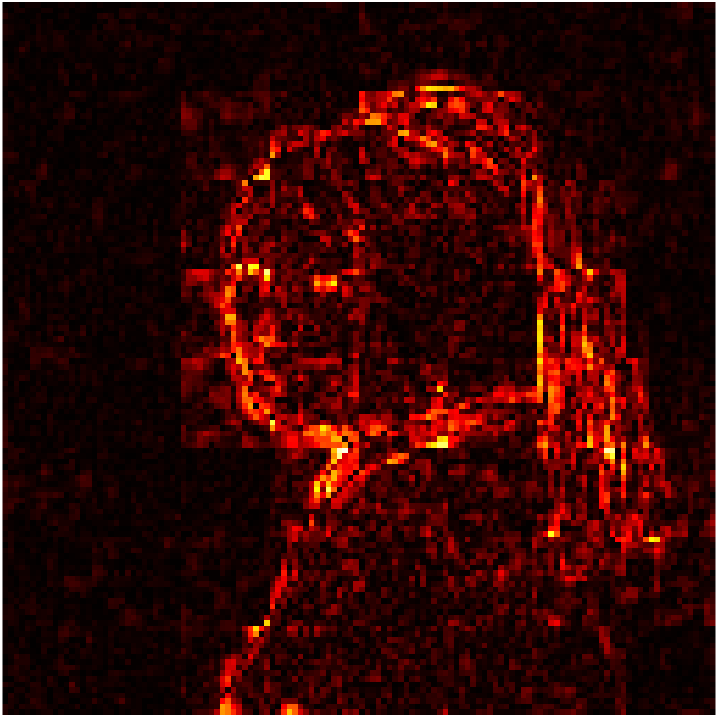}};
    
    \node[above=0.0cm of diff_mrnsd.west, anchor=south, rotate=90] {\begin{tabular}{c} absolute difference \\ $|\bfy_{\rm true} - \bfG\bfx|$ \end{tabular}};

    \node[right=0.0cm of recon_gdnn.east, anchor=west] (parula) {\includegraphics[height=\w\linewidth]{parula.png}};
    \node[below right=0.0cm of parula.north east, anchor=north west] {$1$};
    \node[right=0.0cm of parula.south east, anchor=south west] {$0.047$};
    \node[above right=0.0cm of diff_gdnn.east, anchor=west] (hot) {\includegraphics[height=\w\linewidth]{hot.png}};
    
    \node[below right=0.0cm of hot.north east, anchor=north west] {$0.5$};
    \node[above right=0.0cm of hot.south east, anchor=south west] {$0$};


    \node[below=0.0cm of diff_mrnsd.south, anchor=north] (err_mrnsd) {\vphantom{$\frac{\|\bfy_{\rm true} - \bfG \bfx\|_2}{\|\bfy_{\rm true}\|_2} \approx$} $\bf 1.65\times 10^{-1}$};
    \node[below=0.1cm of err_mrnsd.south, anchor=north] (sparse_mrnsd) {\vphantom{$\frac{\texttt{nnz($\bfx$)}}{\texttt{nnz($\bfy_{\rm true}$)}} \approx$} $\bf 0.09$};

    \node[below=0.0cm of diff_gdnn.south, anchor=north] (err_gdnn) { \vphantom{$\frac{\|\bfy_{\rm true} - \bfG \bfx\|_2}{\|\bfy_{\rm true}\|_2} \approx$} $1.83\times 10^{-1}$};
    \node[below=0.1cm of err_gdnn.south, anchor=north] (sparse_gdnn) {\vphantom{$\frac{\texttt{nnz($\bfx$)}}{\texttt{nnz($\bfy_{\rm true}$)}} \approx$} $0.28$};

    \node[left=1.0cm of err_mrnsd.west, anchor=east]  (err) {relative error, $\frac{\|\bfy_{\rm true} - \bfG\bfx\|_2}{\|\bfy_{\rm true}\|_2} \approx$ };

    \node[below=0.1cm of err.south east, anchor=north east] {relative sparsity, $\frac{\texttt{nnz($\bfx$)}}{\texttt{nnz($\bfy_{\rm true}$)}} \approx$};

     \node[draw, below=0.0cm of sparse_mrnsd.south, anchor=north] (sol_mrnsd) {\includegraphics[scale=0.2]{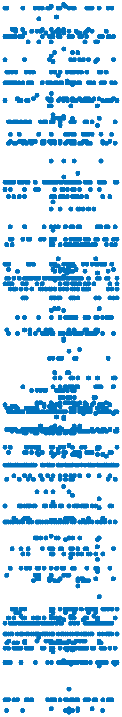}};

    \node[draw, below=0.0cm of sparse_gdnn.south, anchor=north] {\includegraphics[scale=0.2]{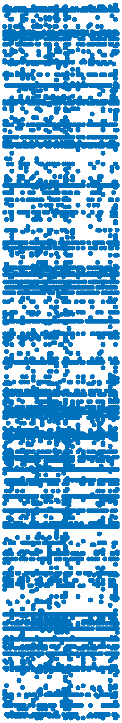}};

    \node[anchor=east] at (err.east |- sol_mrnsd) {matricized coefficients, $\bfx \equiv $};
    
    \end{tikzpicture}

\caption{Reconstructions for the image completion task. The top row (left-to-right) contains the original image, the blurred, noisy image, the reconstruction using \texttt{spMRNSD}, and the reconstruction using \texttt{spNNGD}. The middle row displays the absolute differences between the original and the reconstructions with the relative error and sparsity metrics below. Here, \texttt{spMRNSD} produces the best approximation in terms of both metrics of interest. The bottom row depicts the sparsity patterns of matricized solutions $\bfx$ of size $400 \times 64$ where each column corresponds to the coefficients for one vectorized dictionary patch.}
\label{fig:indicator_recon}
\end{figure}

\begin{figure}
    \centering
    \subfloat[Relative Residual \label{fig:indicator_relerr}]{\includegraphics[width=0.32\linewidth]{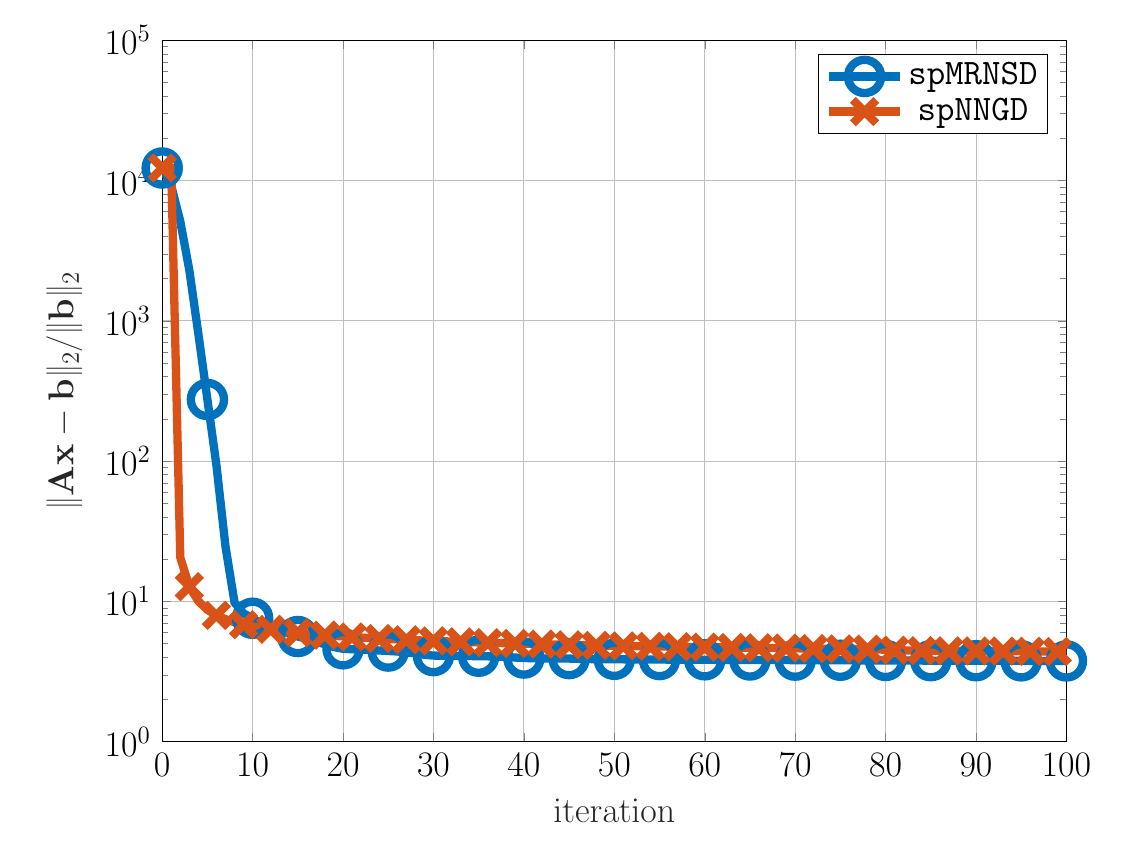}}
    \subfloat[Proxy for Relative Sparsity\label{fig:indicator_sparsity}]{\includegraphics[width=0.32\linewidth]{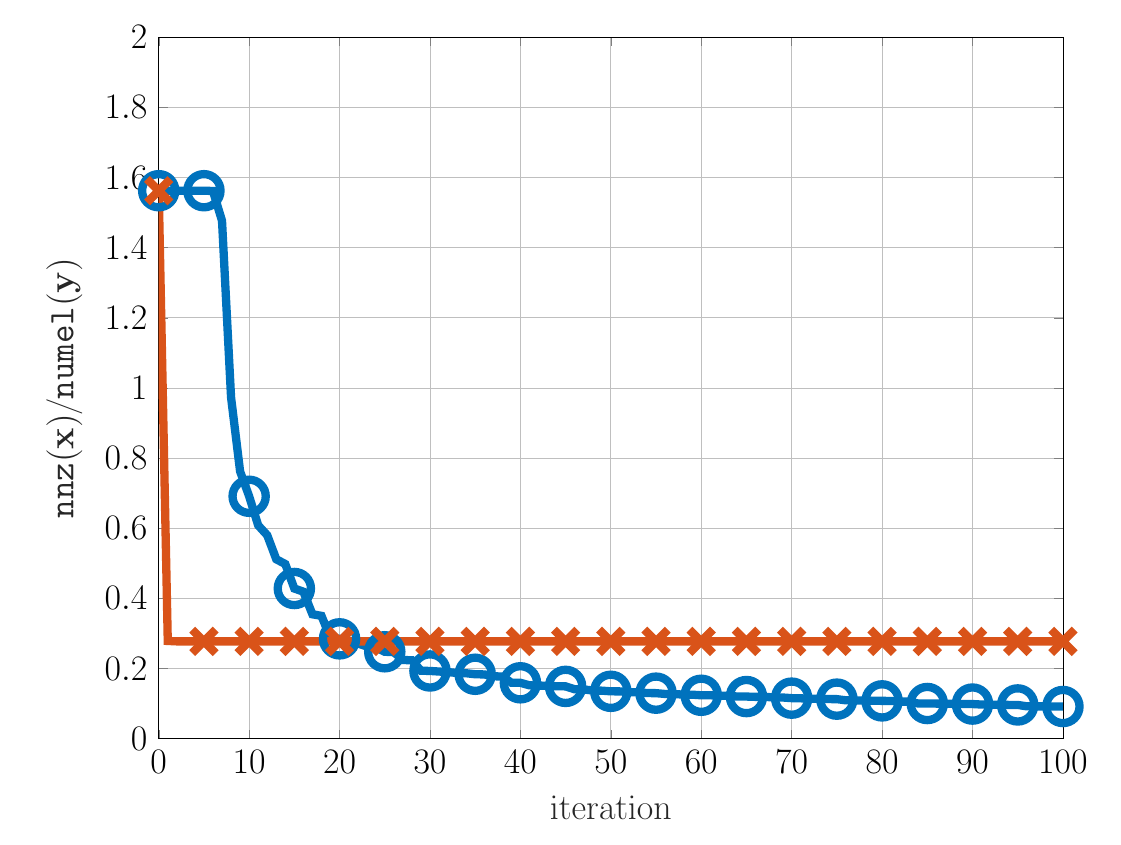}}
    \subfloat[Optimal Step Sizes \label{fig:indicator_alpha}]{\includegraphics[width=0.32\linewidth]{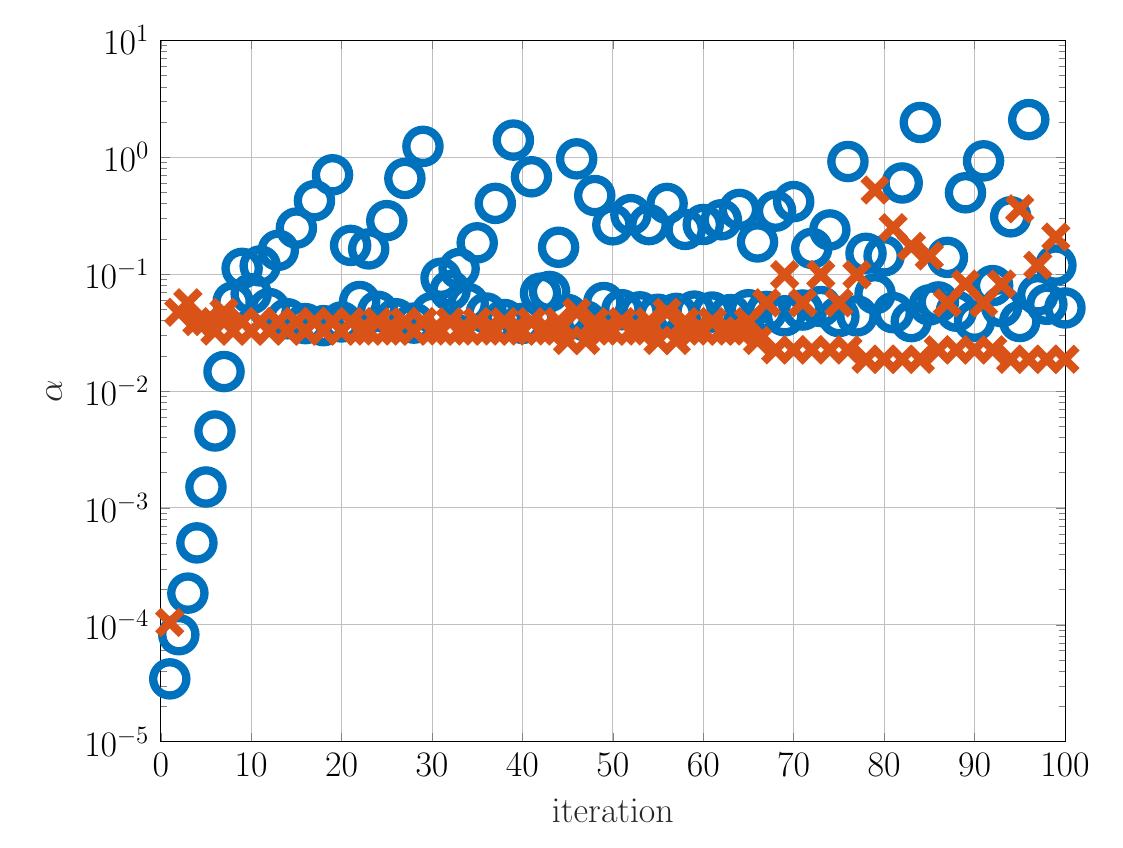}}

    \caption{Convergence behavior for the image completion problem. The relative residual norm (left) converges at similar rates for both algorithms.   The proxy for relative sparsity (middle) converges faster for \texttt{spNNGD}. On average, the step sizes (right) for \texttt{spMRNSD} are larger and the step sizes for \texttt{spNNGD} are more consistent.}
    \label{fig:indicator_convergence}
\end{figure}

We construct the missing pixel operator $\bfC \in \Rbb^{m\times N^2}$ with $m \ll N^2$ that samples uniformly random pixels from an $N\times N$ image. In practice, the rows of $\bfC$ are a subset of rows from the $N^2\times N^2$ identity matrix. As a result, $\bfb = \bfC \bfy$ contains the same as the original image with a percentage of pixels randomly set to zero. 
In our experiment, we use Johannes Vermeer's \emph{Girl with the Pearl Earring} painting,\footnote{\url{https://en.wikipedia.org/wiki/Girl_with_a_Pearl_Earring}} resized to $N\times N = 128\times 128$, and removed $60\%$ of the pixels in our experiments.  We report the representation results in~\Cref{fig:indicator_recon} and the convergence behavior in~\Cref{fig:indicator_convergence}.

In~\Cref{fig:indicator_recon}, we observe that both \texttt{spNNGD} and \texttt{spMRNSD} produce qualitatively good approximations of the original image. The approximation obtained from \texttt{spMRNSD} is a sharper approximation to the original, particularly for the facial features in the image. Additionally, the \texttt{spMRNSD}-generated coefficients are notably sparser than those generated by \texttt{spNNGD}. The sparsity patterns of the un-vectorized coefficients reflect these metrics.  

The convergence behavior in~\Cref{fig:indicator_relerr} and~\Cref{fig:indicator_sparsity} matches that of~\Cref{sec:deblur}. 
In particular, it shows that both algorithms converge to similar levels of accuracy at similar rates and \texttt{spNNGD} quickly reaches a sparse solution while \texttt{spMRNSD} converges more gradually.  Notably and surprisingly, the step size plot in~\Cref{fig:indicator_alpha} reveals that larger steps are taken in $\bfx$-space (\texttt{spMRNSD}) than in $\bfz$-space (\texttt{spNNGD}).  The difference in behavior in terms of both relative sparsity and step size can be attributed to the image completion operator and data normalization.  Specifically, $\bfC$ has rows from an identity matrix, the dictionary $\bfG$ contains entries between $0$ and $1$, and the right-hand side is normalized so its maximum value is $1$.  Thus, the approximation $\bfC\bfG w_{a,c}(\bfz)$ is more directly sensitive to values in the coefficients. For \texttt{spNNGD}, a large update in $\bfz$ could introduce unwanted zeros or could yield large values for $w_{a,c}(\bfz)$. For this reason, we take smaller step sizes. In comparison, for \texttt{spMRNSD}, the soft thresholding operator reduces large values and zeroes out small values, enabling slightly larger step sizes. 

To summarize, the main takeaways from the image completion experiment are that \texttt{spNNGD} and \texttt{spMRNSD} produce approximations with similar relative errors, and \texttt{spMRNSD} produces a sparser solution. As before, \texttt{spNNGD} converges to a sparser solution more quickly, but the gradual progress of \texttt{spMRNSD} produces better metrics overall.  The operator plays a role in determining the effectiveness of the chosen method and may inform the preferred approach.

\subsection{Tomography}\label{sec:tomography}
Tomography is an imaging technique that reconstructs the characteristics of an object's interior by taking cross-sectional measurements from the exterior. Tomography is widely used in a variety of applications ranging from but not limited to geophysics, atmospheric science, physics, materials science, and medical imaging. In the medical field, computer tomography (CT) is a common non-invasive technique used for diagnostic purposes and surgical planning. CT employs penetrating X-rays to pass through an object, with the amount of energy absorbed dependent on the object's properties along the X-ray path. Detectors on the opposite side measure the intensity of the X-rays emitted by the source. By rotating both the X-ray source and receiver around the object, measurements are collected in a sinogram. CT is a classic example of an inverse problem, where the goal is to deduce the object's energy absorbency from the sinogram. For illustration, we consider the 2D computer tomography
\cite{natterer2001mathematics,kak2002principles}. The ill-posedness of this problem requires regularization and a common regularization functional is total variation. 

\begin{figure}
\centering
\begin{tikzpicture}[scale=1.7]
    \pgftext{\includegraphics[width=90pt]{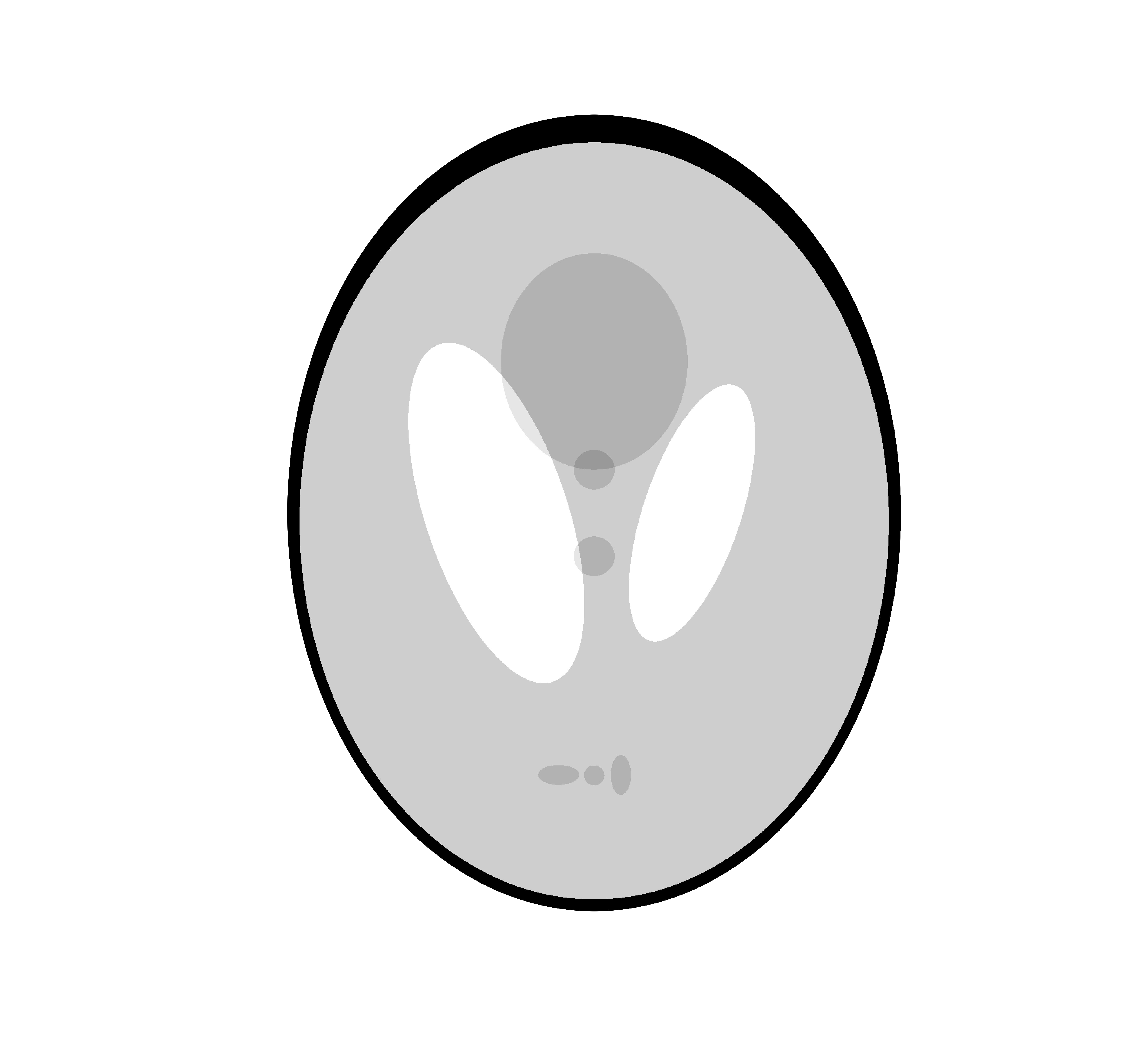}} at (0pt,0pt);

    \begin{scope}[rotate=45,shift={(-0.4,-1.1)}]
      \coordinate [] (A) at (0.7,0);
      \coordinate [] (B) at (0,0);
      \coordinate [] (C) at (0.5,-0.5);
      \path[clip] (A) -- (B) -- (C) -- cycle;
      \draw [matlab1!20,fill=matlab1!20] (B) circle (0.5);
      \draw [matlab1, thick] (A) -- (B);
      \draw [matlab1, thick] (B) -- (C);
      \node[matlab1] at (0.28,-0.1) {};
    \end{scope}
    %
    %
    \begin{scope}[rotate=45]
        \foreach \x in {-1,-0.8,-0.6,-0.4,-0.2,0,0.2,0.4,0.6,0.8,1}
            {
            \draw[ultra thick,-latex,text opacity=1,matlab2]  (-1.5,\x) node{} -- (1.5,\x);
            }
    \end{scope}
    %
    \begin{scope}[rotate=45]
        \draw[gray,fill=gray!20,text opacity=1] (-1.5,-1.1) node[below=0.5,left = 1,black]{x-ray source} rectangle +(-0.2,2.2);
        \draw[black] (-1.6,0) node{};
    \end{scope}
    %
    \begin{scope}[rotate=45]
        \draw[gray,fill=gray!20,text opacity=1] (1.5,-1.1) -- ++(0,2.2) -- ++(0.3,0) -- ++(0.3,-0.3) node[right=0.3,black]{detector} -- ++(0,-1.6) -- ++(-0.3,-0.3) -- cycle ; 
        \draw[black] (1.8,0) node{};
    \end{scope}
\end{tikzpicture}
\caption{Illustration of parallel beam tomography, \cite{ruthotto2018optimal}.}
\label{fig:parallel_tomography}
\end{figure}

\begin{figure}
    \centering
    
    \begin{tikzpicture}
    \def\w{0.2}
    
    \node (orig) at (0,0) {\includegraphics[width=\w\linewidth]{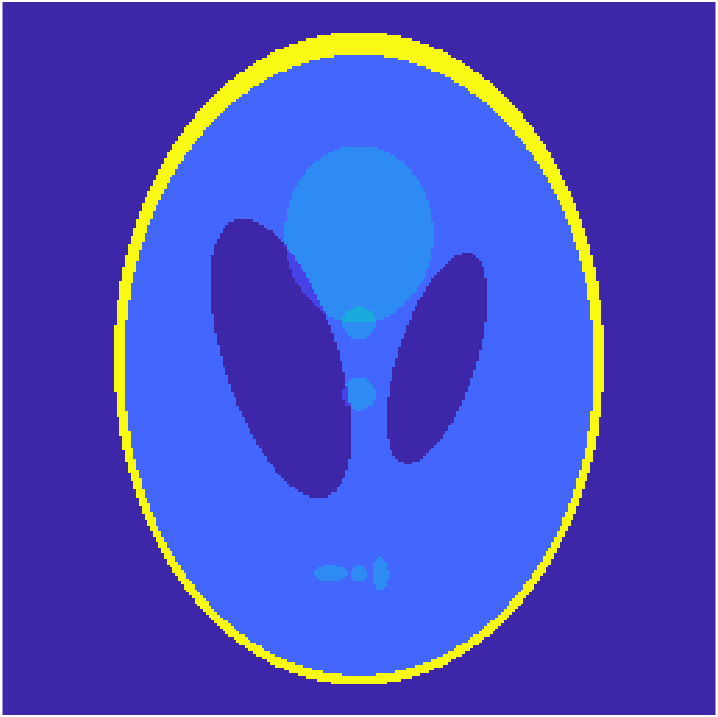}};
    \node[above=0.0cm of orig.north, anchor=south] {\begin{tabular}{c} Original, $\bfy_{\rm true}$ \\ $256 \times 256$ \end{tabular}};

    \node[right=0.0cm of orig.east, anchor=west] (orig) {\includegraphics[width=\w\linewidth]{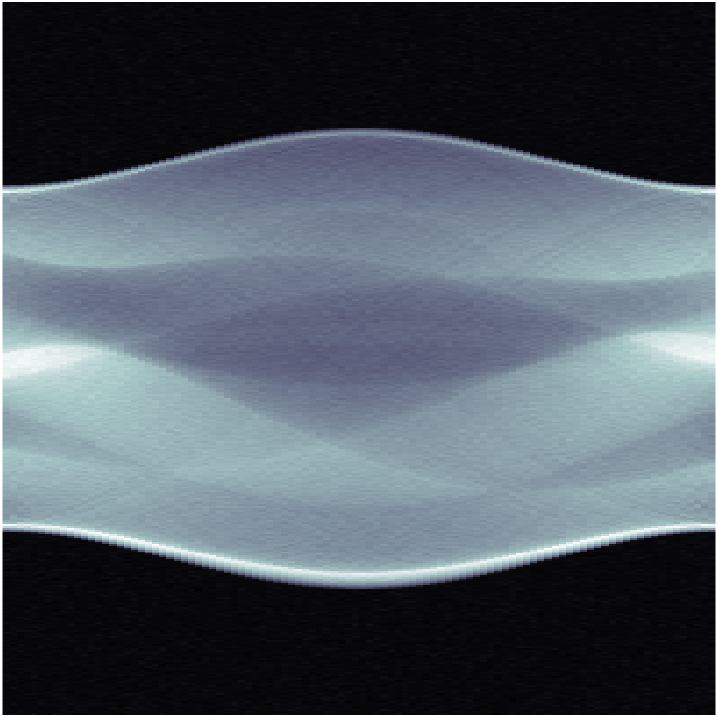}};
    \node[above=0.0cm of orig.north, anchor=south] {\begin{tabular}{c} Sinogram, $\bfb$ \\ $n_{\rm rays} \times n_{\rm proj} = 362\times 100$ \end{tabular}};

    \node[right=0.0cm of orig.east, anchor=west] (recon_mrnsd){\includegraphics[width=\w\linewidth]{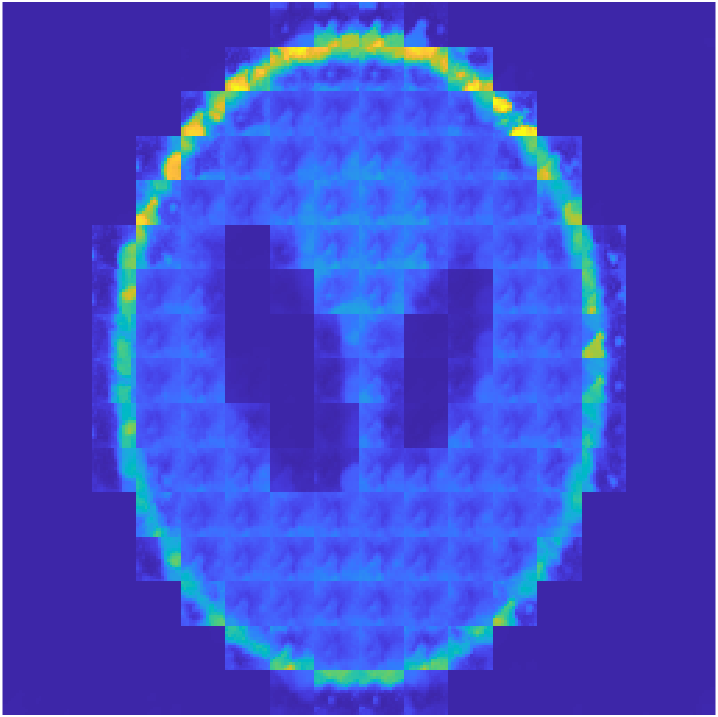}};
    \node[above=0.0cm of recon_mrnsd.north, anchor=south] {\begin{tabular}{c} \texttt{spMRNSD}, $\bfG\bfx$ \\ $\lambda = 10^{0}$\end{tabular}};

    \node[below=0.0cm of recon_mrnsd.south, anchor=north] (diff_mrnsd){\includegraphics[width=\w\linewidth]{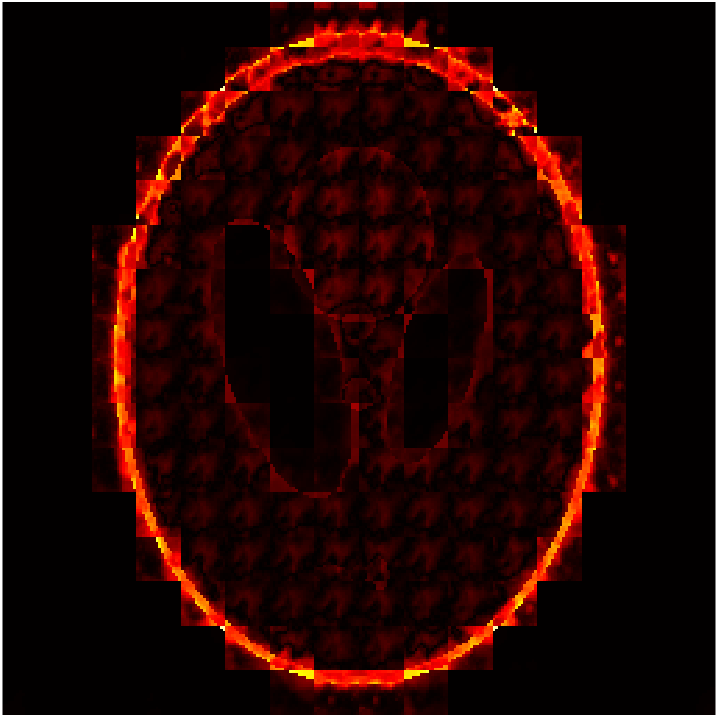}};
    
    \node[right=0.0cm of recon_mrnsd.east, anchor=west] (recon_gdnn){\includegraphics[width=\w\linewidth]{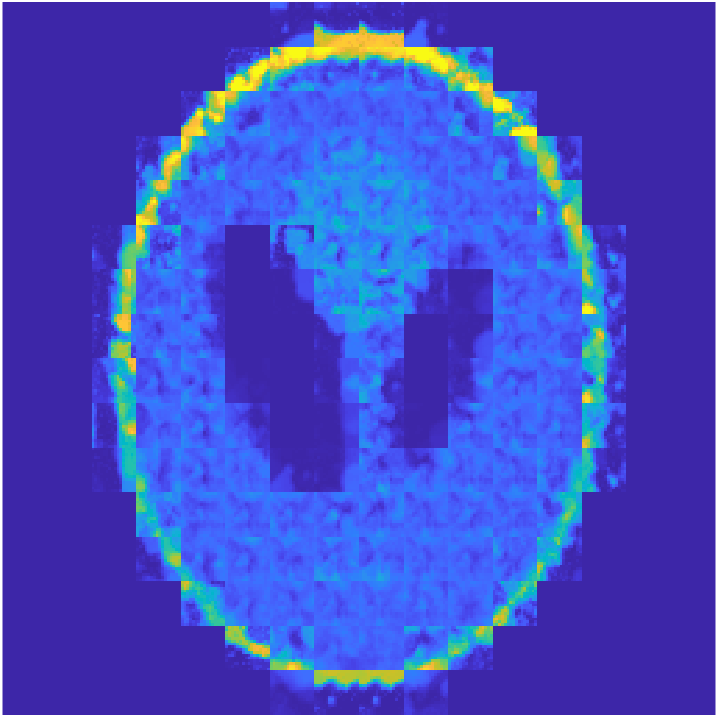}};
    \node[above=0.0cm of recon_gdnn.north, anchor=south] {\begin{tabular}{c} \texttt{spNNGD}, $\bfG\bfx$ \\ $a = 0.5$, $c = -1$\end{tabular}};

    \node[below=0.0cm of recon_gdnn.south, anchor=north] (diff_gdnn){\includegraphics[width=\w\linewidth]{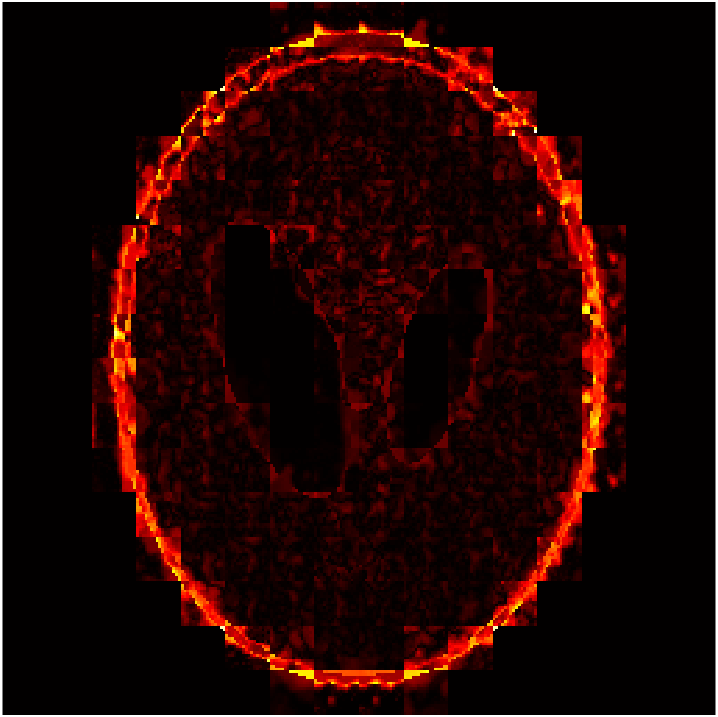}};

    \node[above=0.0cm of diff_mrnsd.west, anchor=south, rotate=90] {\begin{tabular}{c} absolute difference \\ $|\bfy_{\rm true} - \bfG\bfx|$ \end{tabular}};

    \node[right=0.0cm of recon_gdnn.east, anchor=west] (parula) {\includegraphics[height=\w\linewidth]{parula.png}};
    \node[below right=0.0cm of parula.north east, anchor=north west] {$1$};
    \node[right=0.0cm of parula.south east, anchor=south west] {$0$};
    \node[above right=0.0cm of diff_gdnn.east, anchor=west] (hot) {\includegraphics[height=\w\linewidth]{hot.png}};
    
    \node[below right=0.0cm of hot.north east, anchor=north west] {$1.08$};
    \node[above right=0.0cm of hot.south east, anchor=south west] {$0$};


    \node[below=0.0cm of diff_mrnsd.south, anchor=north] (err_mrnsd) { \vphantom{$\frac{\|\bfy_{\rm true} - \bfA\bfx\|_2}{\|\bfy_{\rm true}\|_2} \approx$} $4.99\times 10^{-1}$};
    \node[below=0.1cm of err_mrnsd.south, anchor=north] (sparse_mrnsd) { \vphantom{$\frac{\texttt{nnz(coefficients)}}{\texttt{nnz(image)}} \approx$} $0.45$};

    \node[below=0.0cm of diff_gdnn.south, anchor=north] (err_gdnn) {\vphantom{$\frac{\|\bfy_{\rm true} - \bfG\bfx\|_2}{\|\bfy_{\rm true}\|_2} \approx$} $\bf 4.84\times 10^{-1}$};
    \node[below=0.1cm of err_gdnn.south, anchor=north] (sparse_gdnn) {\vphantom{$\frac{\texttt{nnz(coefficients)}}{\texttt{nnz(image)}} \approx$} $\bf 0.17$};

    \node[left=1.0cm of err_mrnsd.west, anchor=east]  (err) {relative error, $\frac{\|\bfy_{\rm true} - \bfG\bfx\|_2}{\|\bfy_{\rm true}\|_2} \approx$ };

    \node[below=0.1cm of err.south east, anchor=north east] {relative sparsity, $\frac{\texttt{nnz($\bfx$)}}{\texttt{nnz($\bfy_{\rm true}$)}} \approx$};

    \node[draw, below=0.0cm of sparse_mrnsd.south, anchor=north] (sol_mrnsd) {\includegraphics[scale=0.2]{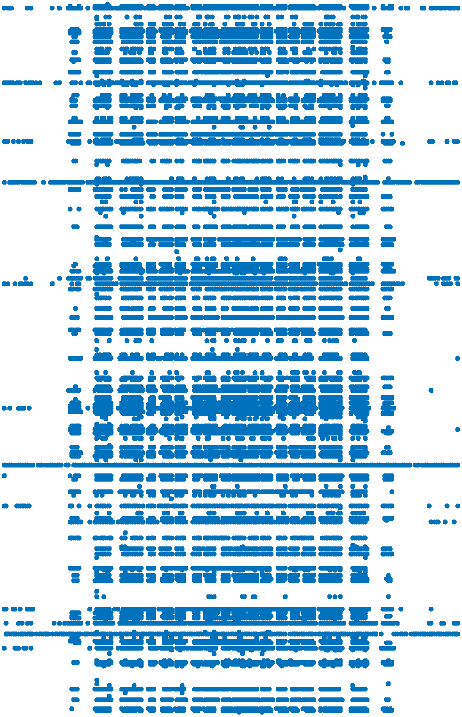}};

    \node[draw, below=0.0cm of sparse_gdnn.south, anchor=north] {\includegraphics[scale=0.2]{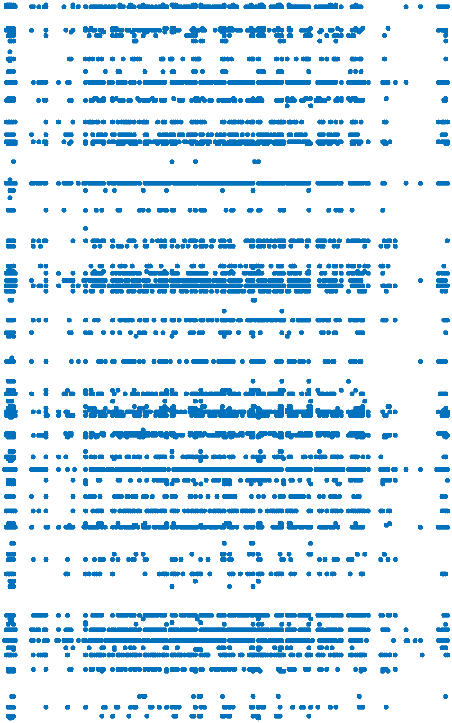}};

    \node[anchor=east] at (err.east |- sol_mrnsd) {matricized coefficients, $\bfx \equiv $};
    
    \end{tikzpicture}

\caption{Reconstructions for the computing tomography task.  The top row (left-to-right) contains the original image, the blurred, noisy image, the reconstruction using \texttt{spMRNSD}, and the reconstruction using \texttt{spNNGD}. The middle row displays the absolute differences between the original and the reconstructions with the relative error and sparsity metrics below. Here, \texttt{spMRNSD} produces the best approximation in terms of both metrics of interest. The bottom row depicts the sparsity patterns of matricized solutions $\bfx$ of size $400 \times 256$ where each column corresponds to the coefficients for one vectorized dictionary patch.}
\label{fig:tomography_recon}
\end{figure}

\begin{figure}
    \centering
    \subfloat[Relative Residual \label{fig:tomography_relerr}]{\includegraphics[width=0.32\linewidth]{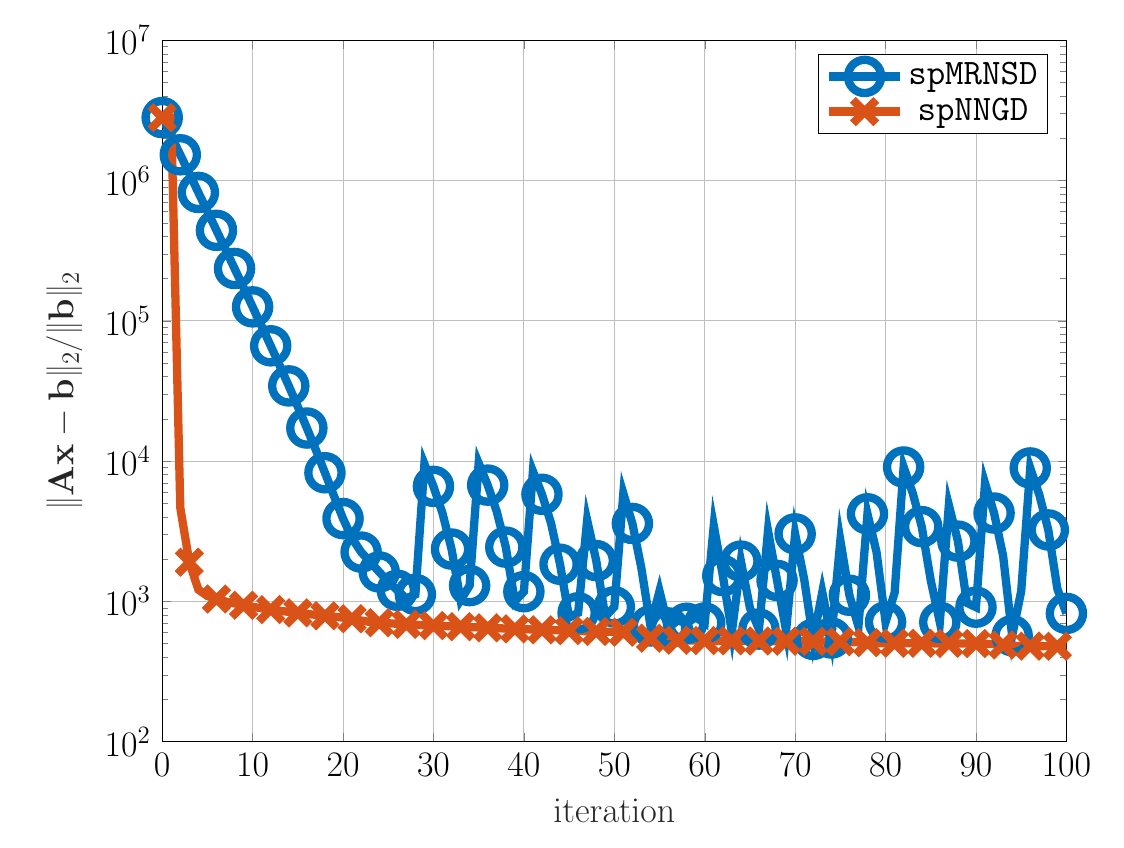}}
    \subfloat[Proxy for Relative Sparsity \label{fig:tomography_sparsity}]{\includegraphics[width=0.32\linewidth]{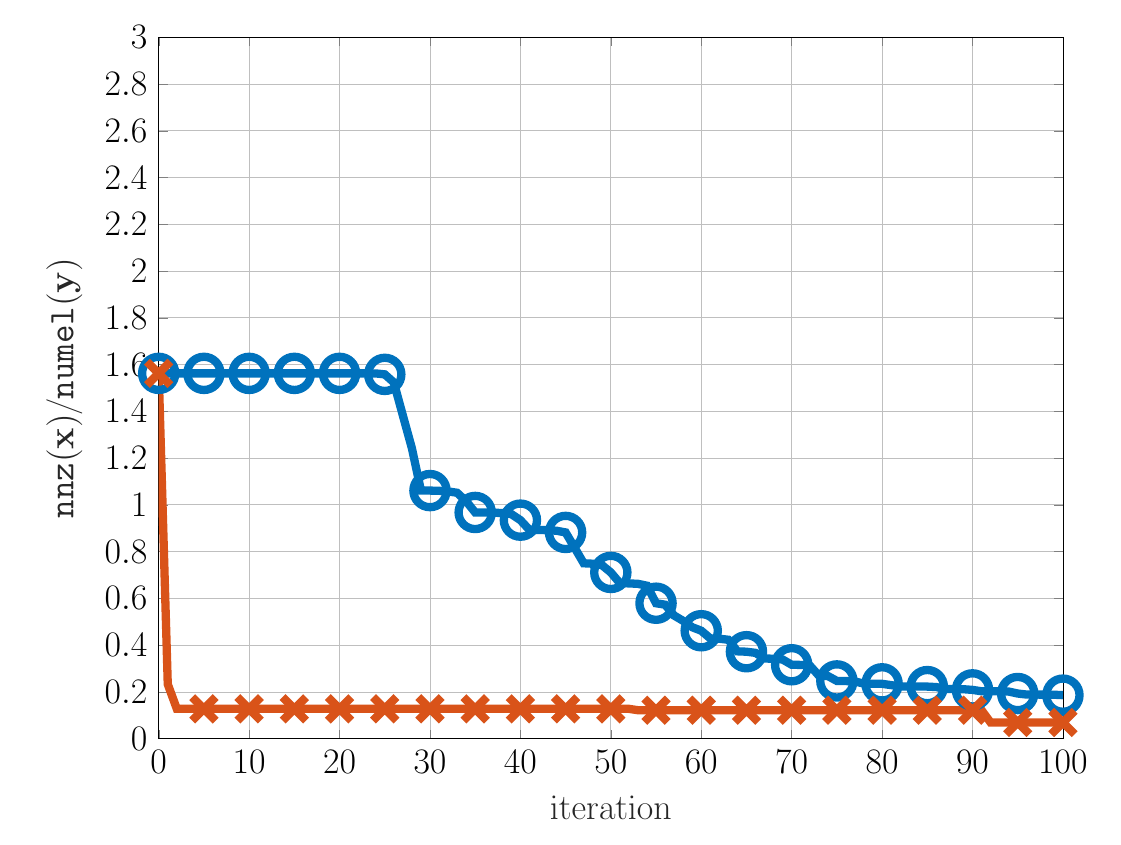}}
    \subfloat[Optimal Step Sizes \label{fig:tomography_alpha}]{\includegraphics[width=0.32\linewidth]{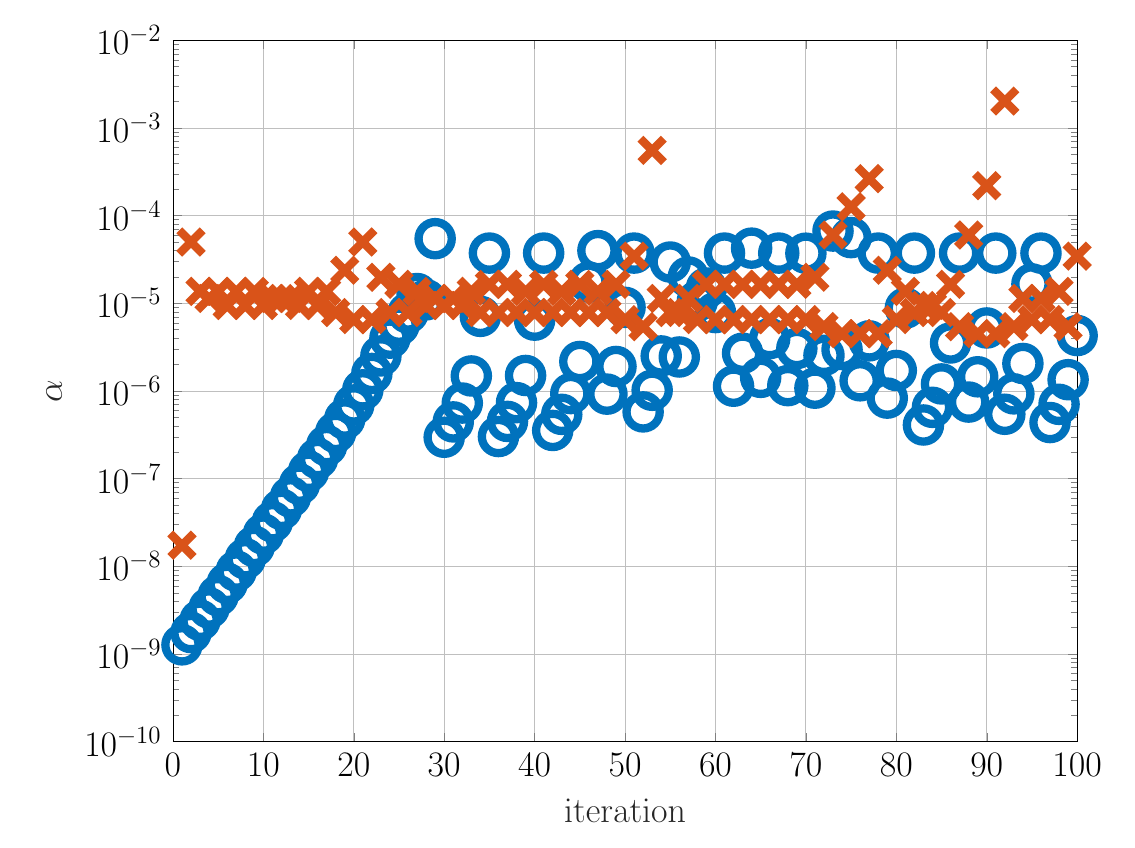}}

    \caption{Convergence behavior for the computed tomography problem.  The relative residual norm (left) converges faster for \texttt{spNNGD} and exhibits more erratic behavior for \texttt{spMNRSD}.   The proxy for relative sparsity (middle) converges faster for \texttt{spNNGD}. On average, the step sizes (right) for \texttt{spNNGD} are larger and more consistent.}
    \label{fig:tomography_convergence}
\end{figure}

In this experiment, we leverage the setup from  IRTools~\cite{IRTools} and AIRTools II~\cite{AIRToolsII} to construct a parallel beam tomography problem; see~\Cref{fig:parallel_tomography} for illustration.   We reconstruct the Shepp-Logan phantom of size $N \times N = 256 \times 256$.  The model matrix $\bfC$ is constructed using $\lfloor \sqrt{2}N\rfloor = 362$ rays and  $100$ rotation angles equispaced between $0$ and $179$ degrees.  For the generalized Tikhonov regularization,  we use a patch smoothing operator with $\mu = 1$. We train for $100$ iterations, use the flower dictionary for representation, and report our results in~\Cref{fig:tomography_recon} and~\Cref{fig:tomography_convergence}.

In~\Cref{fig:tomography_convergence}, we observe \texttt{spNNGD} produces the best approximation results combined with the sparsest solution.  The approximation is particularly sharp in the center of the Shepp-Logan phantom, while the reconstruction obtained from \texttt{spMRNSD} retains many artifacts from the dictionary patches. The sparsity patterns of the solutions are slightly different for the two algorithms as well.  The coefficients from \texttt{spMRNSD} are dense in the center while the \texttt{spNNGD} contain more nonzeros at the edges.   We note that for both algorithms, the largest errors occur on the edges of the Shepp-Logan phantom.  This is perhaps unsurprising because the dictionary elements from the flower dictionary do not contain similar curved features.  A different dictionary prior could yield improved performance on the edge features, but learning the dictionaries is outside the scope of this work. 

 We also mention that the structure of the data motivations for using a patch smoother rather than a pixel smoother or patch border smoother.  We train a dictionary on natural images, but the Shepp-Logan phantom is constructed synthetically.  The Shepp-Logan patches are closer to piecewise constant than natural image patches. Using a patch smoother, we encourage patches of the images to be more similar, which leads to larger regions of almost constant pixel values.

As in the previous experiments, \texttt{spNNGD} converges more quickly to a smaller relative error and a sparser solution and takes larger step sizes on average. However, we observe starkly different convergence behavior in~\Cref{fig:tomography_relerr} than in the previous experiments.  The convergence of \texttt{spMRNSD} is erratic after approximately $30$ iterations.  This behavior is due to the soft-thresholding operator and the mismatch of a search direction in $\bfz$-space and an update in $\bfx$-space.  These two factors make it possible for the algorithm to step in a non-descent direction $\bfx$-space.  We see this occur when the optimal $\alpha$ is large and as the solution becomes sparser.  In contrast, \texttt{spNNGD}, which computes search directions and steps in $\bfz$-space, consistently decays.

\subsection{Superresolution}
\label{sec:superresolution}
Superresolution is a technique to reconstruct a high-resolution image by a sequence of low-resolution images to reveal image details that are unobservable in lower resolution. This task arises in various applications, from enlarging and restoring blurry or damaged photos to improving medical imaging and enhancing surveillance footage~\cite{yang2010image}. For simplicity, let us assume we desire to reconstruct a high-resolution image $\bfx_{\rm true}\in\bbR^{N^2}$ from a sequence of low-resolution images $\bfb_1,\ldots,\bfb_M \in \bbR^{\nu^2}$ with $\nu^2 \ll N^2$.  We assume that the (ridged) image transformation, such as shifts, rotations, and scaling of the high-resolution image to each low-resolution image is known and denoted by $\bfS_1, \ldots, \bfS_M \in \bbR^{N^2 \times N^2}$. With a downsampling/restriction operator $\bfR \in \bbR^{\nu^2 \times N^2}$ we state corresponding superresolution problem as
\begin{equation}
  \min_{\bfx\geq 0 } \quad \norm[2]{
    \begin{bmatrix}
        \bfR\bfS_1 \\
        \vdots \\
        \bfR\bfS_M \\
    \end{bmatrix}
    \bfy_{\rm true} -
    \begin{bmatrix}
        \bfb_1 \\
        \vdots \\
        \bfb_M \\
    \end{bmatrix}
    }^2 + \mu \norm[2]{\bfL\bfy_{\rm true}}^2
  \end{equation}
where we selected a $\ell_2$ regularization term, \cite{slagel2019sampled}. 

\begin{figure}
    \centering
    
    \begin{tikzpicture}
    \def\w{0.2}
    \pgfmathsetmacro{\ww}{\w/3}
    
    \node (orig) at (0,0) {\includegraphics[width=\w\linewidth]{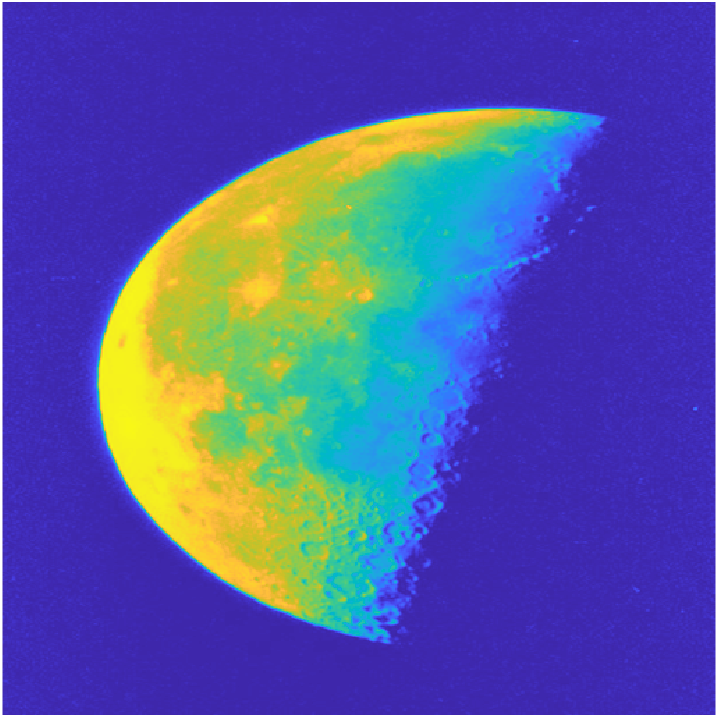}};
    \node[above=0.0cm of orig.north, anchor=south] (orig_header) {\begin{tabular}{c} Original, $\bfy_{\rm true}$ \\ $512\times 512$ \end{tabular}};

    \node[right=-0.05cm of orig.north east, anchor=north west] (tmp0)
    {\includegraphics[width=\ww\linewidth]{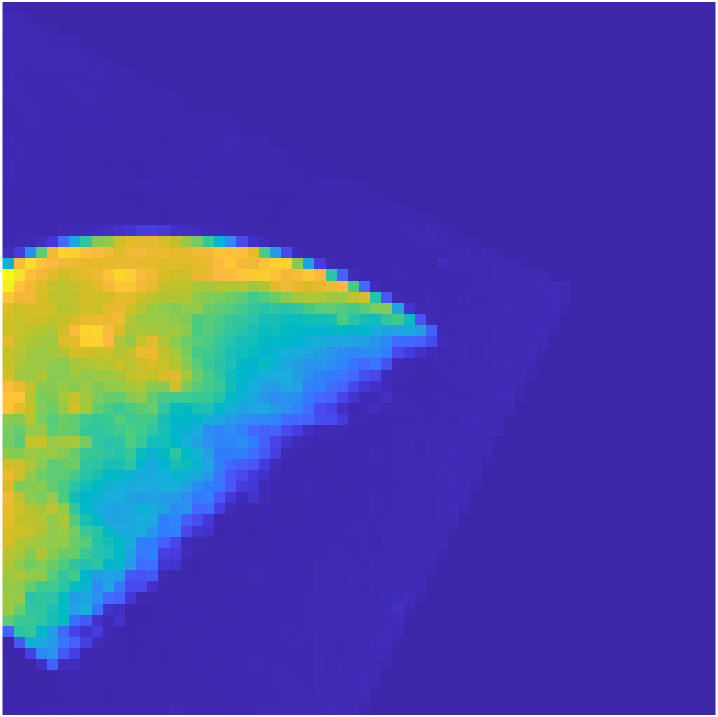}};

    \node[right=-0.2cm of tmp0.east, anchor=west] (tmp)
    {\includegraphics[width=\ww\linewidth]{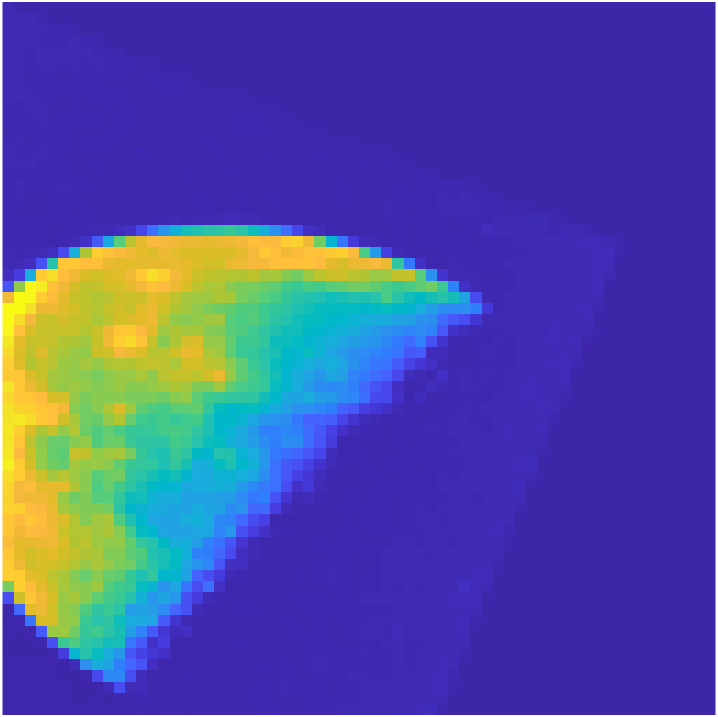}};
    
    \node[right=-0.2cm of tmp.east, anchor=west] (tmp)
    {\includegraphics[width=\ww\linewidth]{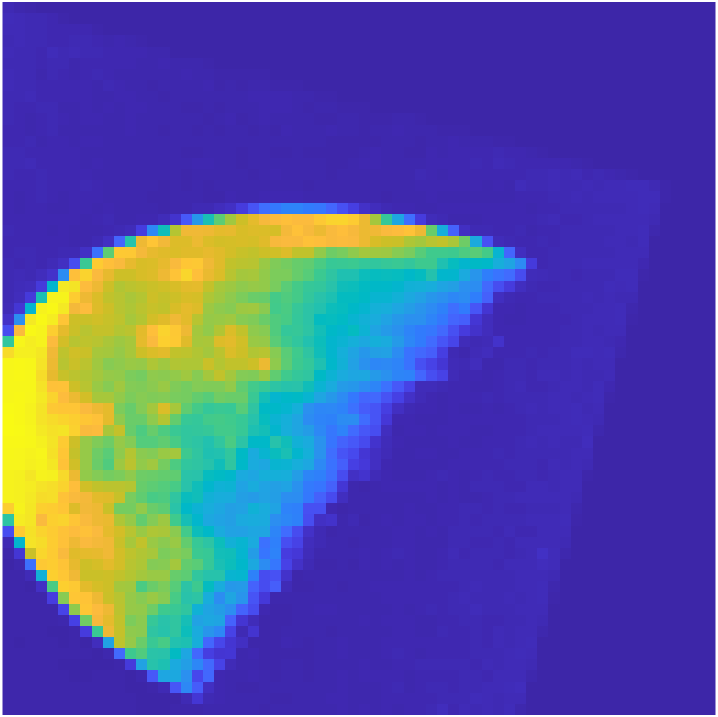}};

    \node[below=-0.2cm of tmp0.south, anchor=north] (tmp1)
    {\includegraphics[width=\ww\linewidth]{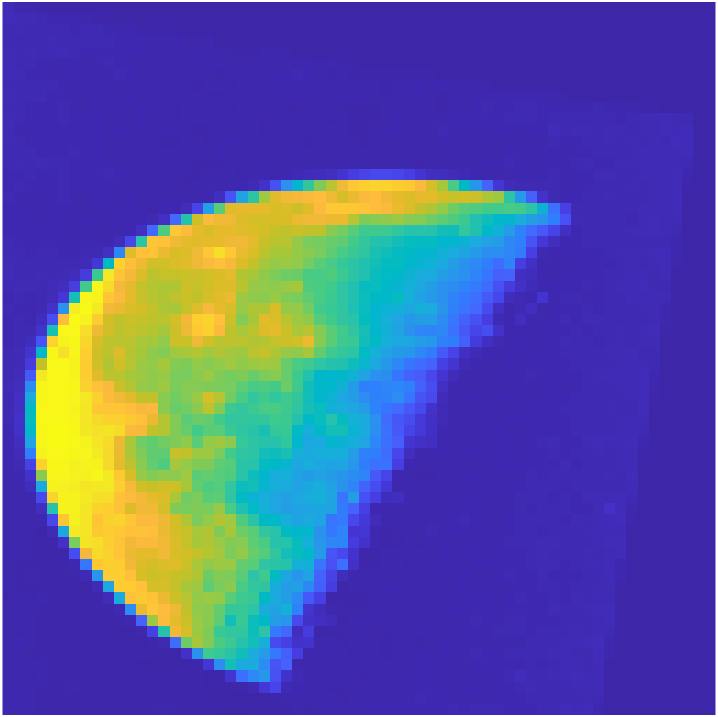}};

    \node[right=-0.2cm of tmp1.east, anchor=west] (tmp)
    {\includegraphics[width=\ww\linewidth]{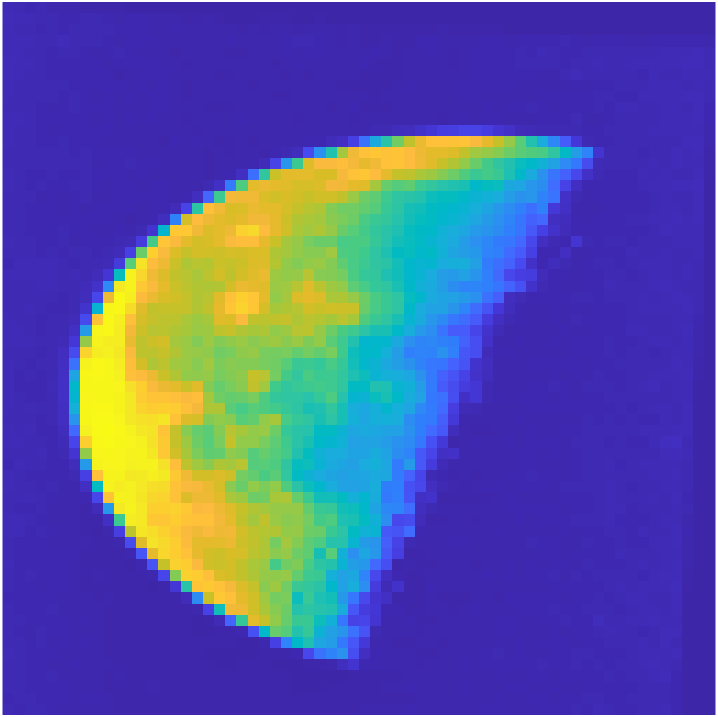}};
    
    \node[right=-0.2cm of tmp.east, anchor=west] (tmp)
    {\includegraphics[width=\ww\linewidth]
    {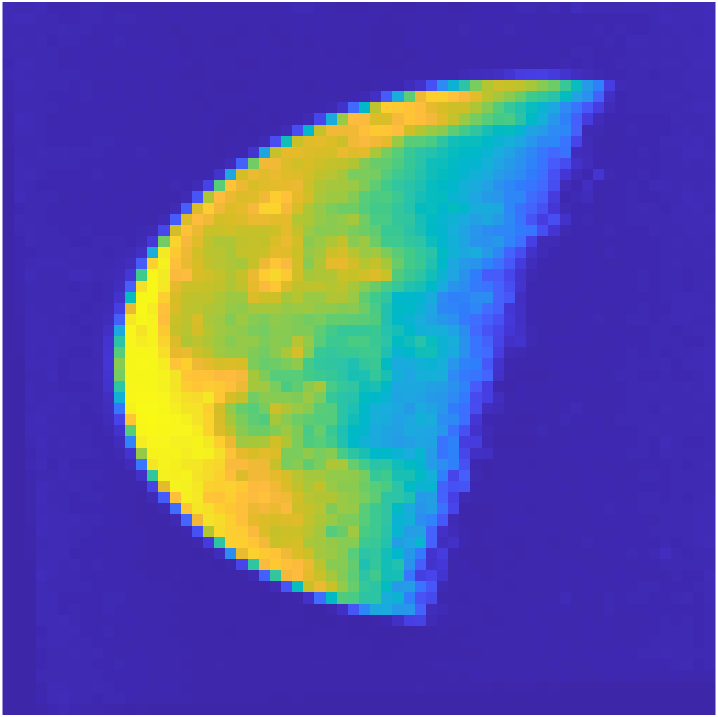}};

    \node[below=-0.2cm of tmp1.south, anchor=north] (tmp2)
    {\includegraphics[width=\ww\linewidth]{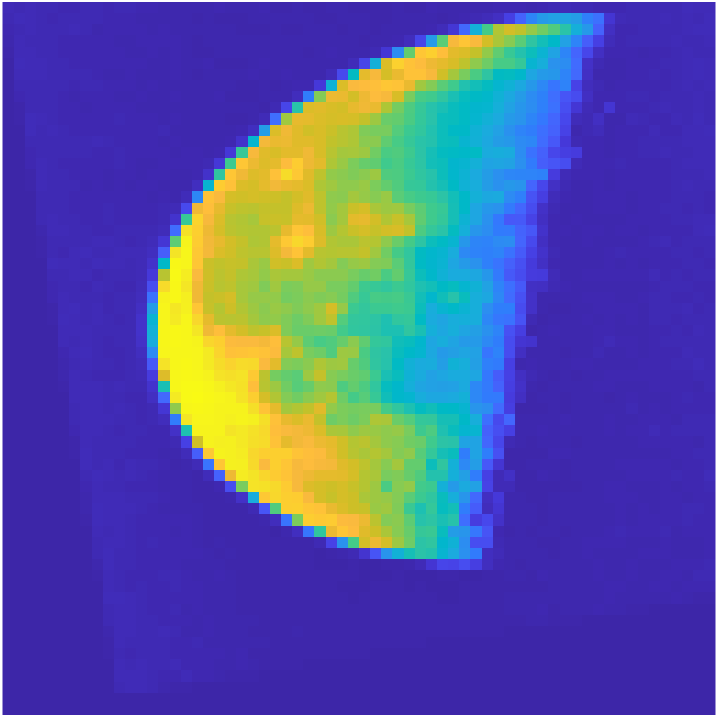}};

    \node[right=-0.2cm of tmp2.east, anchor=west] (tmp)
    {\includegraphics[width=\ww\linewidth]{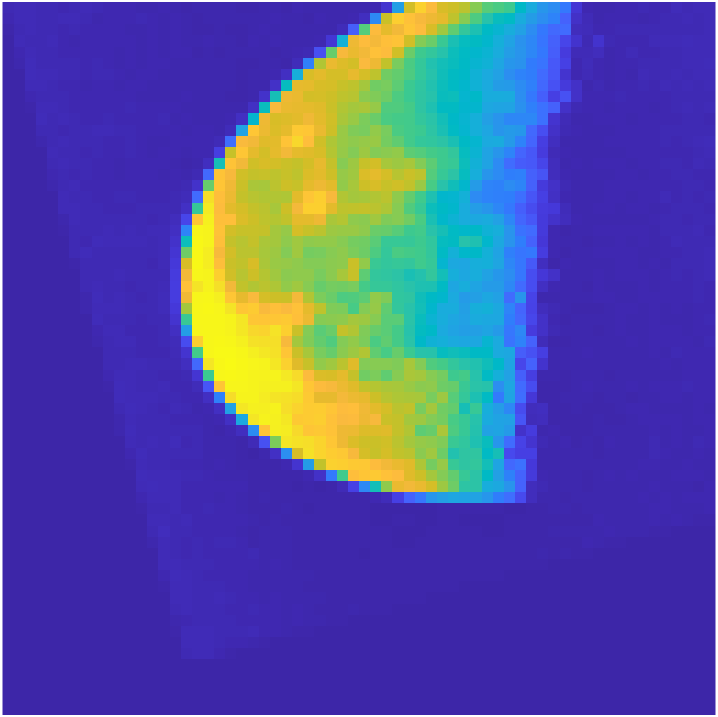}};
    
    \node[right=-0.2cm of tmp.east, anchor=west] (tmp)
    {\includegraphics[width=\ww\linewidth]{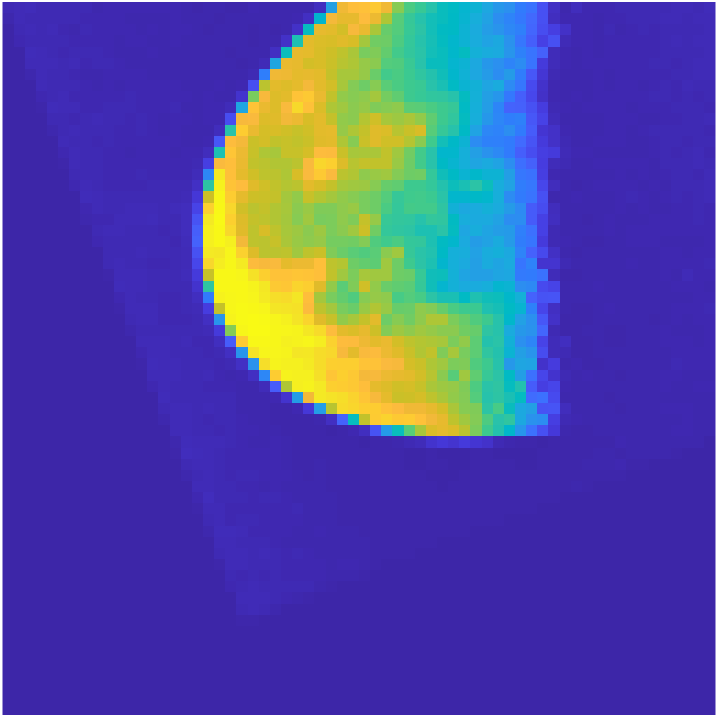}};

    \node[below=-0.2cm of tmp2.south, anchor=north] (tmp3)
    {\includegraphics[width=\ww\linewidth]{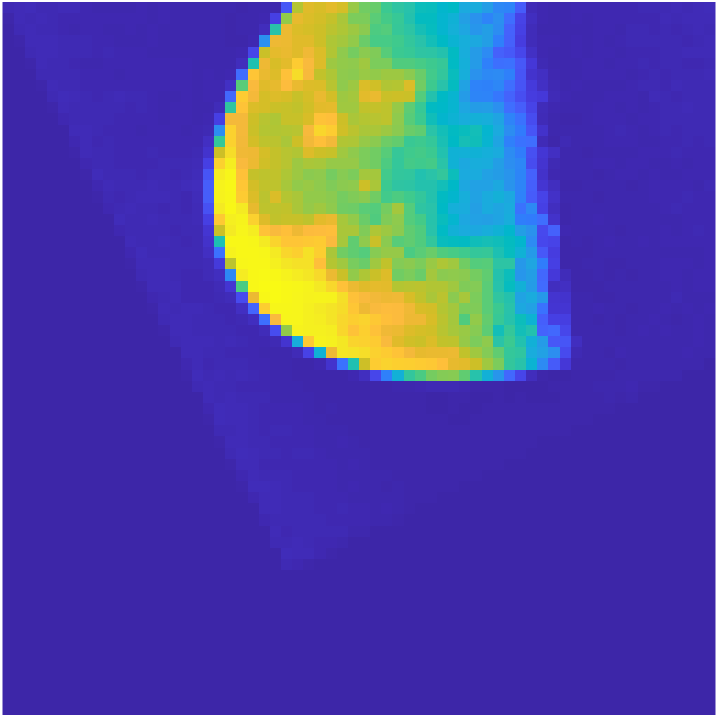}};

    \node[right=0.0cm of orig.east, anchor=west]  (orig) {\phantom{\includegraphics[width=\w\linewidth]{superresolution_gdnn_orig.png}}};

    \node[above=0.0cm of orig.north, anchor=south] {\begin{tabular}{c} Low Resolution, $\bfb_i$ \\ $10$, $64\times 64$ images \end{tabular}};

    \node[right=0.0cm of orig.east, anchor=west] (recon_mrnsd){\includegraphics[width=\w\linewidth]{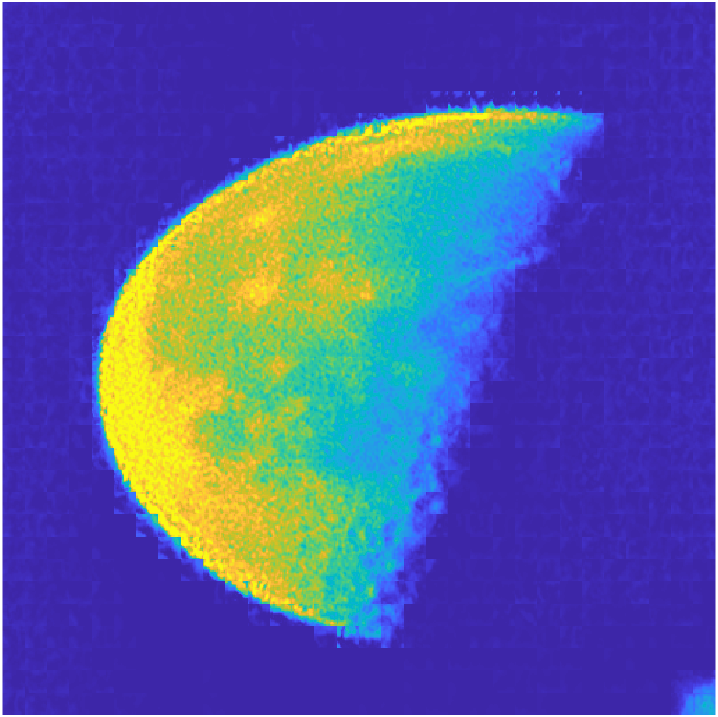}};
    \node[above=0.0cm of recon_mrnsd.north, anchor=south] {\begin{tabular}{c} \texttt{spMRNSD}, $\bfG\bfx$ \\ $\lambda = 10^{-4}$ \end{tabular}};

    \node[below=0.0cm of recon_mrnsd.south, anchor=north] (diff_mrnsd){\includegraphics[width=\w\linewidth]{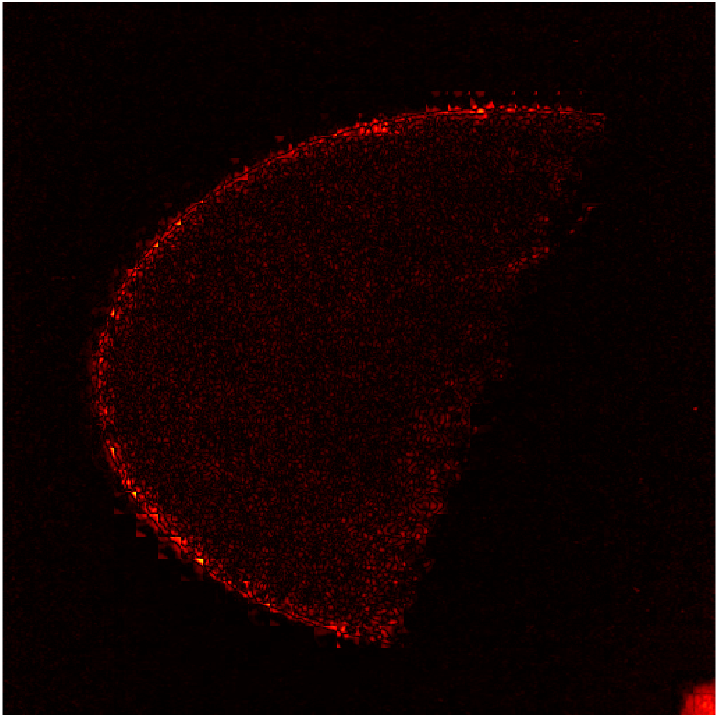}};

    \node[right=0.0cm of recon_mrnsd.east, anchor=west] (recon_gdnn){\includegraphics[width=\w\linewidth]{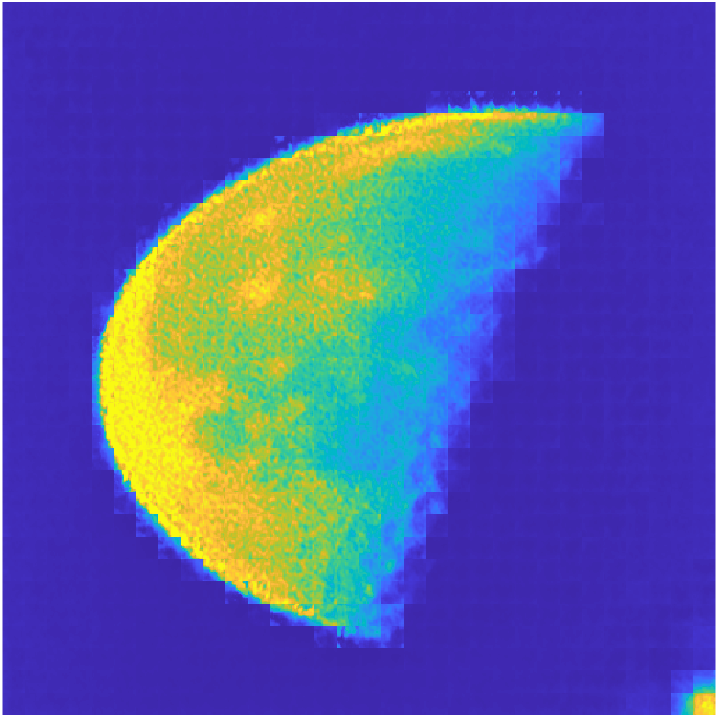}};
    \node[above=0.0cm of recon_gdnn.north, anchor=south] (spNNGD) {\begin{tabular}{c} \texttt{spNNGD}, $\bfG\bfx $ \\ $a=0.2975$, $c = -0.05$\end{tabular}};

    \node[below=0.0cm of recon_gdnn.south, anchor=north] (diff_gdnn){\includegraphics[width=\w\linewidth]{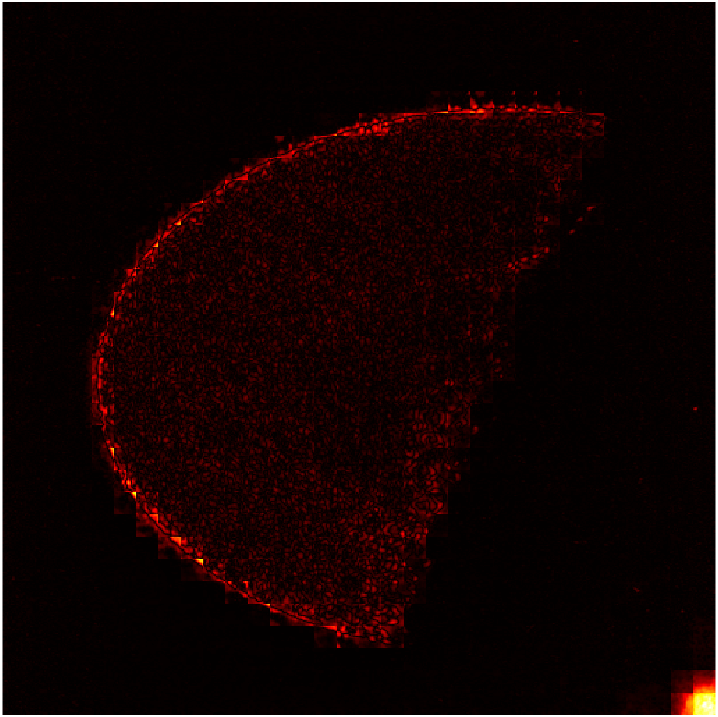}};

        \node[above=0.0cm of diff_mrnsd.west, anchor=south, rotate=90] {\begin{tabular}{c} absolute difference \\ $|\bfy_{\rm true} - \bfG\bfx|$ \end{tabular}};

    \node[right=0.0cm of recon_gdnn.east, anchor=west] (parula) {\includegraphics[height=\w\linewidth]{parula.png}};
    \node[below right=0.0cm of parula.north east, anchor=north west] {$1$};
    \node[right=0.0cm of parula.south east, anchor=south west] {$0$};
    \node[above right=0.0cm of diff_gdnn.east, anchor=west] (hot) {\includegraphics[height=\w\linewidth]{hot.png}};
    
    \node[below right=0.0cm of hot.north east, anchor=north west] {$1$};
    \node[above right=0.0cm of hot.south east, anchor=south west] {$0$};


    \node[below=0.0cm of diff_mrnsd.south, anchor=north] (err_mrnsd) {\vphantom{$\frac{\|\bfx - \hat{\bfx}\|_2}{\|\bfx\|_2} \approx$} $\bf 1.10\times 10^{-1}$};
    \node[below=0.1cm of err_mrnsd.south, anchor=north] (sparse_mrnsd) {\vphantom{$\frac{\texttt{nnz(coefficients)}}{\texttt{nnz(image)}} \approx$} $\bf 0.19$};

    \node[below=0.0cm of diff_gdnn.south, anchor=north] (err_gdnn) {\vphantom{$\frac{\|\bfx - \hat{\bfx}\|_2}{\|\bfx\|_2} \approx$} $1.35\times 10^{-1}$};
    \node[below=0.1cm of err_gdnn.south, anchor=north] (sparse_gdnn) {\vphantom{$\frac{\texttt{nnz(coefficients)}}{\texttt{nnz(image)}} \approx$} $0.98$};

    \node[left=1.0cm of err_mrnsd.west, anchor=east]  (err) {relative error, $\frac{\|\bfy_{\rm true} - \bfG\bfx\|_2}{\|\bfy_{\rm true}\|_2} \approx$ };

    \node[below=0.1cm of err.south east, anchor=north east] {relative sparsity, $\frac{\texttt{nnz($\bfx$)}}{\texttt{nnz($\bfy_{\rm true}$)}} \approx$};

    \node[draw, below=0.0cm of sparse_mrnsd.south, anchor=north] (sol_mrnsd) {\includegraphics[scale=0.2]{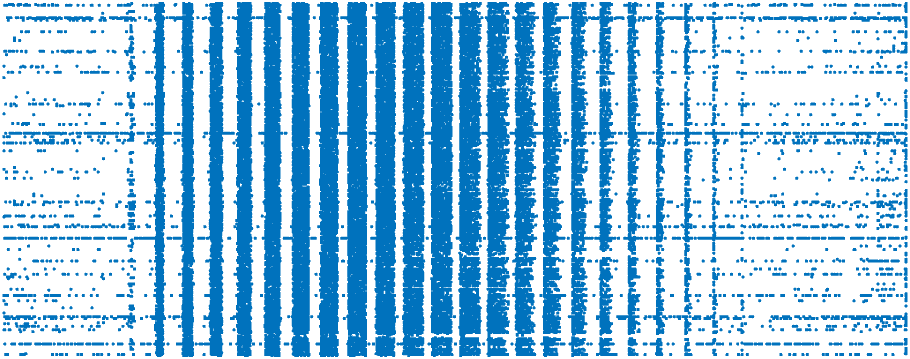}};

    \node[below=0.0cm of sparse_gdnn.south, anchor=north] {\includegraphics[scale=0.2]{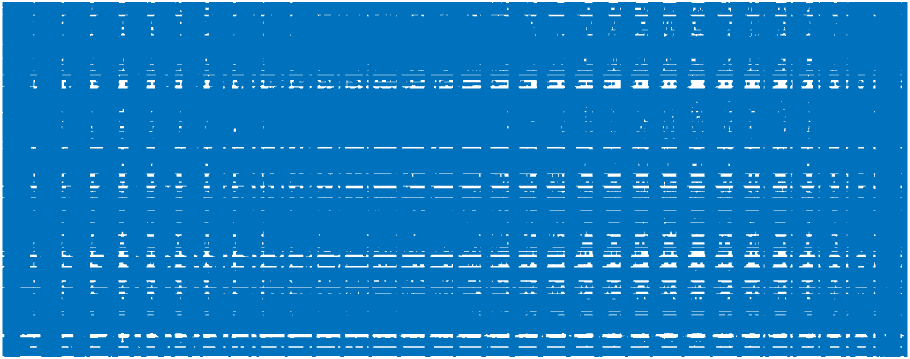}};

    \node[draw, below=0.0cm of sparse_gdnn.south, anchor=north] {\phantom{\includegraphics[scale=0.2]{superresolution_gdnn_solution.png}}};

    \node[anchor=east] at (err.east |- sol_mrnsd) {matricized coefficients, $\bfx \equiv $};

    \end{tikzpicture}

\caption{Reconstructions for the superresolution task. The top row (left-to-right) contains the original image, the blurred, noisy image, the reconstruction using \texttt{spMRNSD}, and the reconstruction using \texttt{spNNGD}.  The middle row displays the absolute differences between the original and the reconstructions with the relative error and sparsity metrics below. Here, \texttt{spMRNSD} produces the best approximation in terms of both metrics of interest.  The bottom row depicts the sparsity patterns of matricized solutions $\bfx$ of size $400 \times 1024$ where each column corresponds to the coefficients for one vectorized dictionary patch.  Note that the lower right patch for both approximations has significantly more error than the other patches due to the boundary conditions of the operator.}
\label{fig:superresolution_recon}
\end{figure}

\begin{figure}
    \centering
    \subfloat[Relative Residual \label{fig:superresolution_relerr}]{\includegraphics[width=0.32\linewidth]{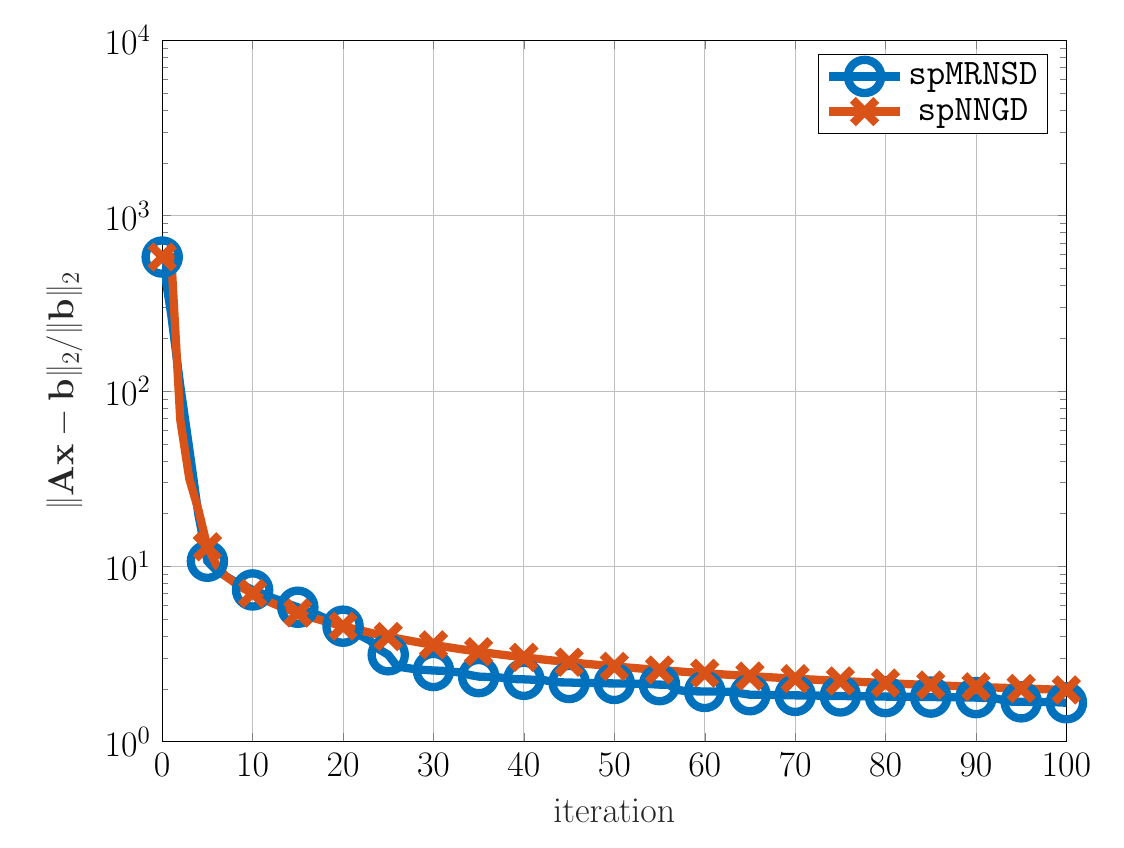}}
    \subfloat[Proxy for Relative Sparsity \label{fig:superresolution_sparsity}]{\includegraphics[width=0.32\linewidth]{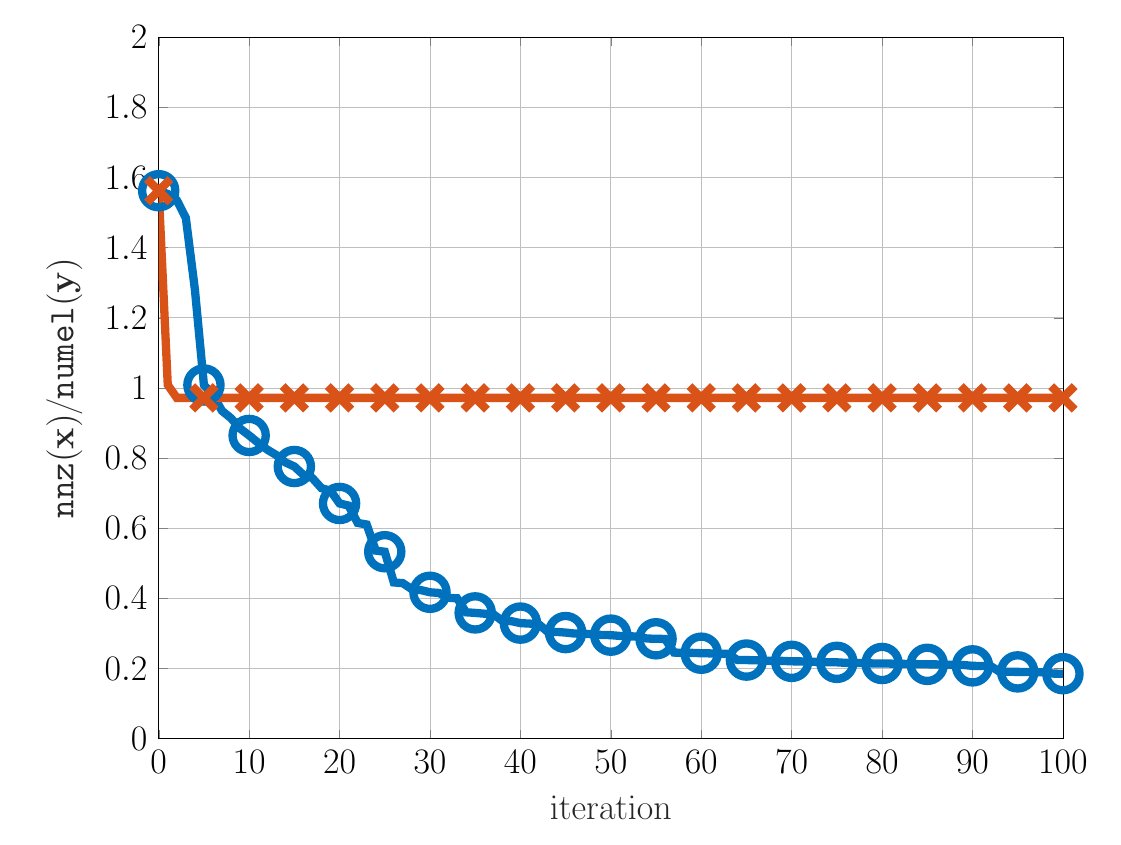}}
    \subfloat[Optimal Step Sizes \label{fig:superresolution_alpha}]{\includegraphics[width=0.32\linewidth]{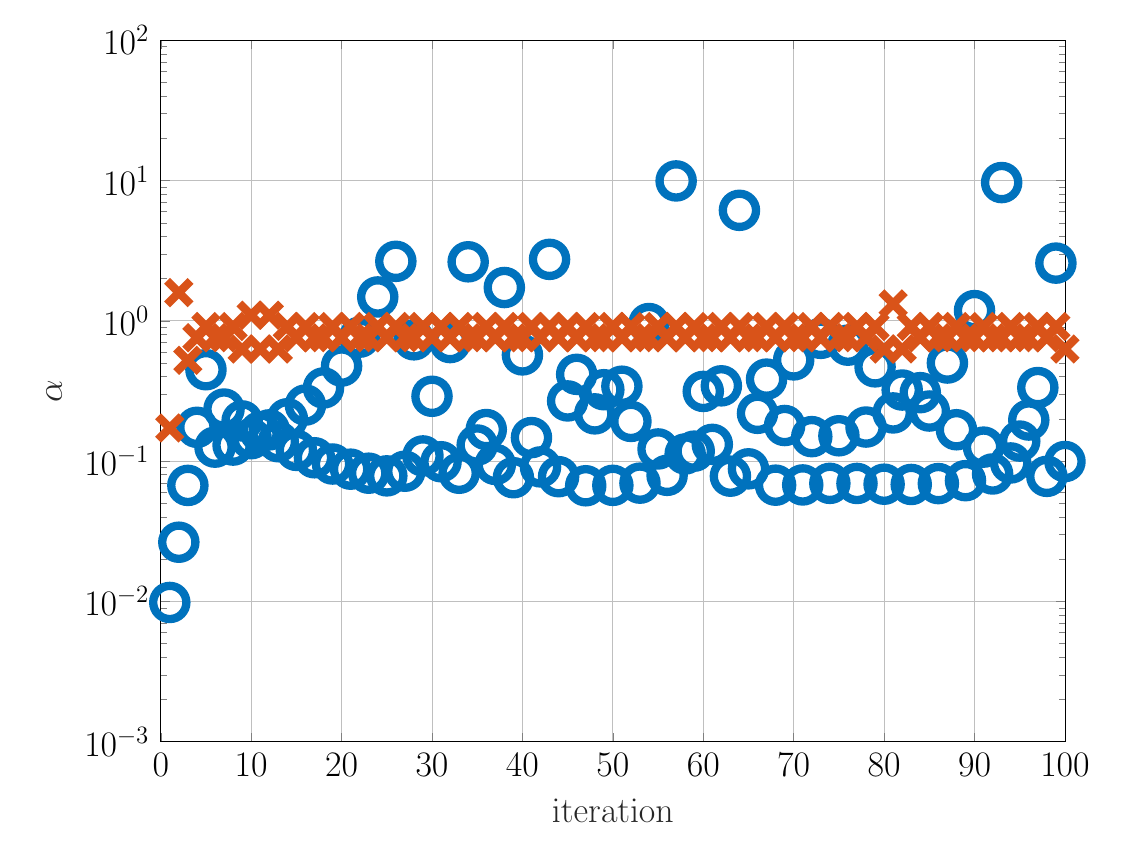}}

    \caption{Convergence behavior for the superresolution problem. The relative residual norm (left) converges at similar rates for both  algorithms.  The proxy for relative sparsity (middle) reaches a lower level to \texttt{spMRNSD}.   On average, the step sizes (right) for \texttt{spNNGD} are larger and more consistent.}
    \label{fig:superresolution_convergence}
\end{figure}

In our experiments, we consider an resized $N\times N = 512 \times 512$ builtin {\sc Matlab} image of the moon.\footnote{\texttt{img = imread('moon.tif');}} 
We form $\theta = 10$ low-resolution images, downsampled by a factor of eight. 
Thus, the low resolution images $\bfB_i$ are of size $64 \times 64$ and the restriction operator is a matrix $\bfR_+\in \Rbb^{(N/8)^2 \times N^2}$. 
Poetically, we reconstruct the moon using the earth dictionary (\Cref{fig:dictionary_patches}).  
We use $100$ iterations and a patch neighbor smoothing regularizer with $\mu = 10^{-2}$.  
We report the results of the superresolution experiment in~\Cref{fig:superresolution_recon} and~\Cref{fig:superresolution_convergence}.

\begin{figure}
    \centering
    \begin{tikzpicture}
        \def\w{1.8}
        \def\s{0.05}
        \foreach \i in {1,...,5}{
            \node (n0) at (\i*\w+\i*\s,0) {\includegraphics[width=\w cm]{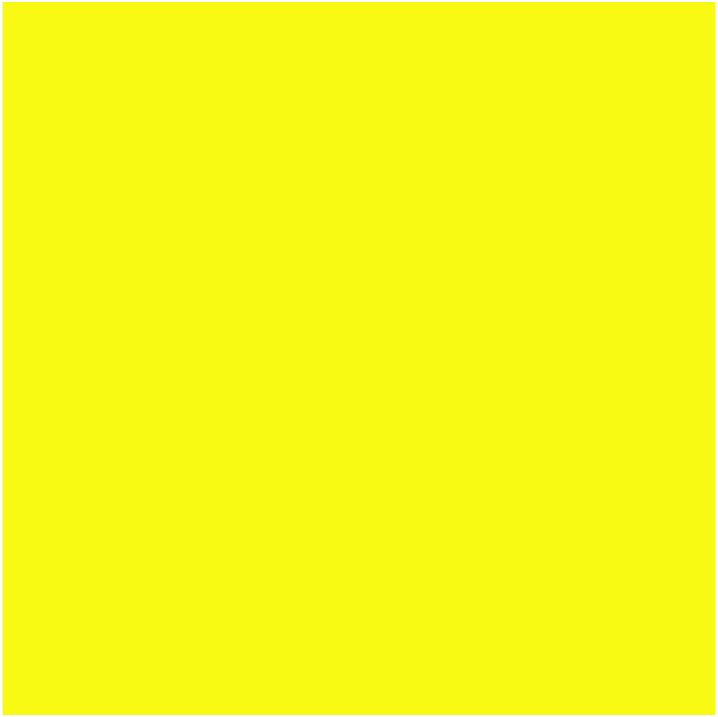}};
            \node (m0) at (\i*\w+\i*\s,-\w-\s) {\includegraphics[width=\w cm]{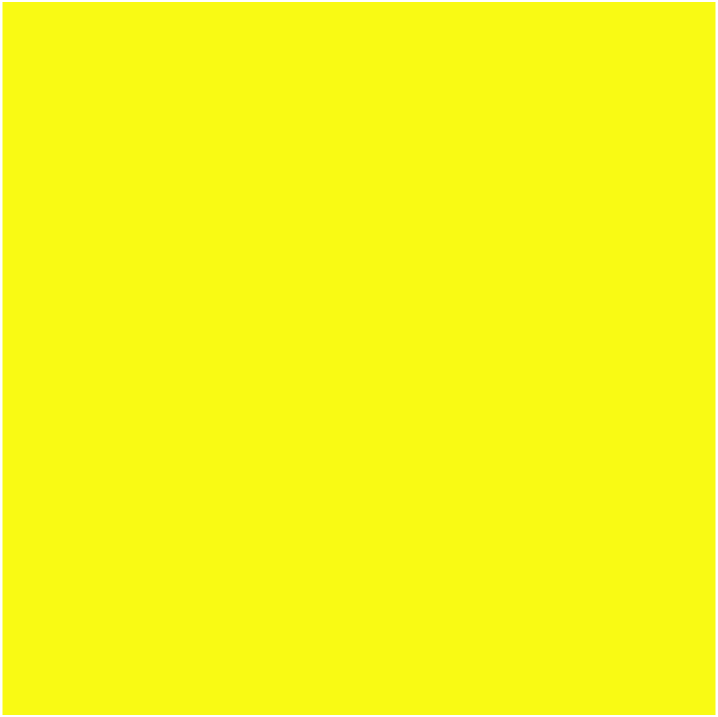}};
            \ifthenelse{\i = 1}
            {
            \node[above=0.0cm of m0.west, anchor=south, rotate=90]{\texttt{spNNGD}};
            
            \node[above=0.0cm of n0.west, anchor=south, rotate=90]{\texttt{spMRNSD}};}
            {}
        }

        \pgfmathsetmacro{\h}{2*\w + \s}
        \node[right=0.0cm of n0.north east, anchor=north west] (parula) 
        {\includegraphics[height=\h cm]{parula.png}};
        \node[below right=0.0cm of parula.north east, anchor=north west] {$1$};
        \node[right=0.0cm of parula.south east, anchor=south west] {$0$};
        
        \def\i{1}
        \def\j{5}
        \draw[->, line width=2pt] (\i*\w+\i*\s,0.5*\w+0.25) -- node[above, midway] {iteration} (\j*\w+\j*\s,0.5*\w+0.25);
        
    \end{tikzpicture}
    \caption{Superresolution approximations for first five iterations.  Approximations from both algorithms contain the lower right corner due to the boundary conditions of the operator.  We show the approximations on the same color scale as the original image, which can be misleading in the initial iteration.  The first approximation contains values larger than $1$, hence all of the pixels appear to be yellow.}
    \label{fig:superresolution_approx_convergence}
\end{figure}

In~\Cref{fig:superresolution_recon}, we observe that both algorithms produce similar quality approximations, and both have the same artifact in the bottom right corner.  We conjecture that the cause of this artifact is twofold.   First, the location of this error may be due to the operator $\bfC$ which rotates the images counterclockwise about the upper left corner.  The second is that the initial guess plays a crucial role in the approximation, and we see this artifact appear in early iterations (see~\Cref{fig:superresolution_approx_convergence}). In this work, we opted for sacrificing some approximation quality to promote more sparsity, and hence focus on the results with the lower right corner error. 

We note in~\Cref{fig:superresolution_recon}  \texttt{spNNGD} does not produce any sparsity.  We conjecture that sparsity is harder to obtain in superresolution partially because the right-hand sides already have many nonzeros (e.g., $\frac{\texttt{nnz}(\bfb_1)}{\texttt{numel}(\bfb_1)} \approx 56\%$ and overall $\frac{\texttt{nnz}(\bfb)}{\texttt{numel}(\bfb)} \approx 75\%$) and partially because we use small patches to reconstruct a higher resolution image. Per the latter point, as shown in~\cite{NewmanKilmer2020}, the ratio between patch size and image size can significantly affect the sparsity pattern.   We have seen in previous experiments that \texttt{spNNGD} settles on a sparsity pattern in early iterations.  Fixating on one sparsity pattern while maintaining a good approximation leads to a non-sparse solution. In contrast, \texttt{spMRNSD} gradually promotes sparsity, and performs better along this metric.

The first iteration of~\Cref{fig:superresolution_approx_convergence} exhibits similar behavior to the previous experiments with the two algorithms converging at similar rates in terms of relative error and \texttt{spNNGD} promoting sparsity faster. However, we quickly see \texttt{spNNGD} stagnate with minuscule step sizes, eventually halting progress at about $40$ iterations.  As expected, the early step sizes in~\Cref{fig:superresolution_alpha} are larger for \texttt{spNNGD}.  Interestingly, both step sizes start to decay at similar rates for iterations $20$ through $30$, which is the iteration at which \texttt{spMNRSD} stops becoming sparser in~\Cref{fig:superresolution_sparsity}.


\section{Conclusions}\label{sec:conclusion}

We have demonstrated success in solving inverse problems using patch dictionary representations to obtain non-negative, sparse solutions in image deblurring (\Cref{sec:deblur}), image completion (\Cref{sec:indicator}), computed tomography (\Cref{sec:tomography}), and superresolution (\Cref{sec:superresolution}).  We proposed and compared two different modifications to the modified residual norm steepest descent (MRNSD) algorithm to produce solutions with the desired structure. Motivated by~\cite{NewmanKilmer2020}, we first introduced a sparsity-promoting MRNSD (\texttt{spMRNSD}) by adding $\ell_1$ regularization and iterating with a soft-thresholding operator.  Our key contribution was using the correct optimal step size based on the soft-thresholding operator update.  We compared this soft-thresholding approach to a new gradient descent approach with a sparsity-promoting, non-negative mapping in \texttt{spNNGD}.   We observed that typically \texttt{spNNGD} converges faster to a sparse solution than \texttt{spMRNSD} and takes larger step sizes.  However, the sparsity pattern can be fixed after a few iterations, which can limit improvements along that metric.  In comparison, \texttt{spMRNSD} gradually increases sparsity via the soft-thresholding operator. The performance of each algorithm can vary significantly based on the task, depending on the given operator, the algorithmic hyperparameters, and the dictionary used for the representation. In general, it seems the gradual approach of \texttt{spMRNSD} offers more promise than the sparsity-promoting, non-negative mapping of \texttt{spNNGD}, but the advantages of the latter could lead to interesting future work.

The complex relationship between all of the parameters led to our comparison of the two algorithms and motivates several exciting extensions.  For \texttt{spMRNSD}, we can we can explore other implicit mappings.  Currently, we use the exponential mapping to enable simple steps in the solution space.  We could perform similar implicit mappings and solution space steps with other mappings. For example, with $\bfx = \bfz^2$, our search in solution space (without forming $\bfz$ explicitly) would use search directions $\bfs = \diag{\tfrac{1}{2}\sqrt{\bfx}} \bfg$. For \texttt{spNNGD}, we can consider adaptive approaches to modify the hyperparameters as we iterate. For example, if we observe stagnation of the sparsity pattern, we may gradually shift the mapping to promote more sparsity. Additionally, we could consider accelerating the selection of the step size via linearization rather than solving a one-dimensional optimization problem explicitly.   Combining the two algorithms, we could consider initializing with a few iterations of \texttt{spNNGD} to reach a sparser solution, and then refine with \texttt{spMRNSD}.

\bigskip
\noindent {\bf Acknowledgement:} This work was partially supported by the National Science Foundation (NSF) under grants [DMS-2309751] and [DMS-2152661] for Newman and Chung, respectively. Any opinions, findings, conclusions, or
recommendations expressed in this material are those of the authors and do not necessarily reflect the views of the National Science Foundation.

\printbibliography

\end{document}